\documentclass[10pt]{amsart}
\usepackage{graphicx,color}

\newtheorem{theorem}{Theorem}[section]
\newtheorem{lemma}[theorem]{Lemma}

\theoremstyle{definition}

\theoremstyle{remark}
\newtheorem{remark}[theorem]{Remark}

\numberwithin{equation}{section}

\keywords{Schr{\"o}dinger operator, lattice, inverse scattering\\
\noindent 
\ \ \ \ This work is supported by Grant-in-Aid for Scientific Research (S) 15H05740, (B) 16H0394, (C) 17K05303 and Grant-in-Aid for Young Scientists (B) 16K17630, Japan Society for the Promotion of Science. }

\title[Inverse Scattering on Perturbed Lattices]
{Inverse scattering for Schr{\"o}dinger Operators on Perturbed Lattices}
\author[K. Ando]{Kazunori ANDO}
\address[K. Ando]{Department of Electrical and Electronic Engineering and Computer Science, Ehime University, Bunkyo-cho 3, Matsuyama, Ehime, 790-8577, Japan}
\email{ando@cs.ehime-u.ac.jp}
\author[H. Isozaki]{Hiroshi ISOZAKI}
\address[H. Isozaki]{Professor Emeritus,
University of Tsukuba, Tennoudai 1-1-1, Tsukuba, Ibaraki, 305-8571, Japan}
\email{isozakih@math.tsukuba.ac.jp}
\author[H. Morioka]{Hisashi MORIOKA}
\address[H. Morioka]{Faculty of Science and Engineering,
Doshisha University, Tataramiyakodani 1-3, Kyotanabe, Kyoto, 610-0394, Japan}
\email{hmorioka@mail.doshisha.ac.jp}

\date{\today}

\begin{document}
\baselineskip 14pt

\maketitle

\begin{abstract}
We study the inverse scattering for Schr{\"o}dinger operators on locally perturbed periodic lattices. We show that the  associated scattering matrix is equivalent to the Dirichlet-to-Neumann map for a boundary value problem on a finite part of the graph, and reconstruct scalar potentials as well as the graph structure from the knowledge of the S-matrix. In particular, we give a procedure for probing  defects in  hexagonal lattices (graphene). 
\end{abstract}


\section{Introduction}


\subsection{Inverse scattering for the continuous model}
The aim of this paper is to investigate inverse problems of scattering for Schr{\"o}dinger operators on locally perturbed periodic lattices. 
For the sake of comparison, we begin with recalling the progress of multi-dimensional inverse scattering theory for the continuous model, made in the last several decades. In ${\bf R}^d$ with $d\geq 2$, consider the Schr{\"o}dinger equation
\begin{equation}
(- \Delta + V(x))u = \lambda u, \quad x \in {\bf R}^d,
\label{S1SchroedingerEq}
\end{equation}
where $V(x)$ is a real-valued compactly supported potential. 
Given a beam of quantum mechanical particles with energy $\lambda > 0$ and incident direction $\omega \in S^{d-1}$, the scattering state is described by
 a solution $u(x)$ of the equation (\ref{S1SchroedingerEq}) satisfying
\begin{equation}
u(x) \simeq e^{i\sqrt{\lambda}\omega\cdot x} + \frac{e^{i\sqrt{\lambda}r}}{r^{(d-1)/2}}a(\lambda;\theta,\omega), \quad
{\rm as} \quad r = |x| \to \infty,
\label{S1FarfieldCont}
\end{equation}
where $\theta = x/r$. The first term of the right-hand side corresponds to the plane wave coming from the direction $\omega$, and the second term
represents the spherical wave scattered to the direction $\theta$. 
The function $a(\lambda;\theta,\omega)$ is called the {\it scattering amplitude}, 
and $|a(\lambda;\theta,\omega)|^2$ is the number of particles scattered to the direction $\theta$. Therefore, it is directly related to the physical experiment.
Let $S(\lambda) = I - 2\pi iA(\lambda)$, where $A(\lambda)$ is the integral operator with kernel $C(\lambda)a(\lambda;\theta,\omega), C(\lambda) = 2^{-1/2}\lambda^{(d-2)/4}(2\pi)^{-d/2}$. Then,  $S(\lambda)$ is a unitary operator on $L^2(S^{d-1})$,  called 
 (Heisenberg's) {\it S-matrix}. The goal of inverse scattering is to reconstruct $V(x)$ from the S-matrix. There is also a time-dependent picture of the scattering theory. Let $H_0 = - \Delta$, $H = H_0 + V$, both of which are self-adjoint on $L^2({\bf R}^d)$. Then, the wave operators
\begin{equation}
W_{\pm} = {\mathop{\rm s-lim}_{t\to\pm\infty}}\, e^{itH}e^{-itH_0}
\nonumber
\end{equation}
exist and are partial isometries with initial set $L^2({\bf R}^d)$ and final set $\mathcal H_{ac}(H)$ = the absolutely continuous subspace for $H$. This implies that for any $f \in \mathcal H_{ac}(H)$, there exist $f_{\pm} \in L^2({\bf R}^d)$ such that 
$$
\|e^{-itH}f - e^{-itH_0}f_{\pm}\|_{L^2({\bf R}^d)} \to 0 \quad {\rm as} \quad t \to \pm \infty.
$$
The scattering operator
\begin{equation}
S = \big(W_+\big)^{\ast}W_-
\nonumber
\end{equation}
is then unitary on $L^2({\bf R}^d)$, and we have $Sf_- = f_+$. By the conjugation by the Fourier transformation
\begin{equation}
\big(\mathcal F_0f)(\lambda) = (2\pi)^{-d/2}\int_{{\bf R}^d}e^{-i\sqrt{\lambda}\omega\cdot x}f(x)dx,
\nonumber
\end{equation}
$S$ is represented as
\begin{equation}
\big(\mathcal F_0S\mathcal F_0^{\ast}f\big)(\lambda,\omega) = \big(S(\lambda)f(\lambda,\cdot)\big)(\omega), 
\label{S1Smatrixformula}
\end{equation}
for $f \in L^2\big((0,\infty);L^2(S^{d-1});\frac{1}{2}\lambda^{(d-2)/2}d\lambda\big)$, 
where $S(\lambda)$ is the S-matrix.

There are three methods for the reconstruction of the potential from $S(\lambda)$. The first one is the high-energy Born approximation due to 
Faddeev \cite{Fa56}:
\begin{equation}
\lim_{\lambda \to \infty}A(\lambda;\theta_{\lambda},\omega_{\lambda}) = C_d\widetilde V(\xi),
\end{equation}
where $\widetilde V(\xi)$ is the Fourier transform of $V$, $C_d$ is a constant and $\theta_{\lambda}, \omega_{\lambda} \in S^{d-1}$ are suitably chosen so that
$\sqrt{\lambda}\big(\theta_{\lambda} - \omega_{\lambda}\big) \to \xi$.
The second method is the multi-dimensional Gel'fand-Levitan theory, again due to Faddeev \cite{Fa76}, which opened a breakthrough, although some parts are formal, to the characterization of the S-matrix and the reconstruction of the potential. The key tool was the  new Green function of  Laplacian introduced in \cite{Fa66}. 
The third method was given by Sylvester-Uhlmann \cite{SyUh87}, Nachman \cite{Na88}, Khenkin-Novikov \cite{KhNo87}, \cite{No88}, which is based on the $\overline{\partial}$-theory, a complex analytic view point for Faddeev's Green function. Let us stress here 
that Sylvester-Uhlmann found Faddeev's Green function independently of Faddeev's approach in studying inverse boundary value problems, which is  another stream of inverse problem initiated by Calder{\'o}n \cite{Ca80}.  The associated exponentially growing solution for the
 Schr{\"o}dinger equation and its analogue are now used in various inverse boundary value problems. 
 
 We explain the details of this third method. It is essential here that the perturbation $V(x)$ is compactly supported. Assuming that the support of $V(x)$ lies in a bounded 
 domain $D_{int} \subset {\bf R}^d$, we consider the boundary value problem
 \begin{equation}
 \left\{
 \begin{split}
 & \big(- \Delta + V(x) - \lambda)u = 0 \quad {\rm in} \quad D_{int}, \\
 & u = f \quad {\rm on} \quad \partial D_{int}. 
 \end{split}
 \right.
 \nonumber
 \end{equation}
 The mapping
 \begin{equation}
 \Lambda : f \to \frac{\partial u}{\partial \nu}\Big|_{\partial D_{int}},
 \nonumber
 \end{equation}
 $\nu$ being the unit normal to $\partial D_{int}$, 
 is called the Dirichlet-to-Neumann map, or simply {\it D-N map}. For any fixed energy $\lambda > 0$, one can show that the scattering amplitude $A(\lambda)$ determines the D-N map and vice versa, if 
 $\lambda$ is not the Dirichlet eigenvalue for the domain $D_{int}$. Using  Faddeev's Green function or exponentially growing solution, one can then reconstruct the potential $V(x)$ from the D-N map.
 
 Let us also recall here that the inverse boundary value problem raised by Calder{\'o}n deals with the following equation appearing in electrical impedance tomography
   \begin{equation}
 \left\{
 \begin{split}
 & \nabla\cdot \big(\gamma(x)\nabla u\big) = 0 \quad {\rm in} \quad D_{int}, \\
 & u = f \quad {\rm on} \quad \partial D_{int},
 \end{split}
 \right.
 \nonumber
 \end{equation}
 where $\gamma(x) = \left(\gamma_{ij}(x)\right)$ is a positive definite matrix, representing the electric conductivity of the body in question.
 The D-N map is defined 
  as the operator
 \begin{equation}
 \Lambda_{\gamma}f = \gamma(x)\frac{\partial u}{\partial \nu}\Big|_{\partial D_{int}}.
\label{S1DNmapCont}
 \end{equation}
For further details of the inverse scattering theory and inverse boundary value problems, see e.g.  review articles \cite{ChaColPaiRun}, \cite{Is03}, \cite{Uhl09}.


\subsection{Inverse scattering on the perturbed periodic lattice} 
In this paper, we  consider  periodic lattices whose finite parts are perturbed by  potentials or some deformation, i.e. addition or removal of edges and vertices. 
 Since the perturbation is finite dimensional,  the wave operators and the scattering operator are introduced in the same way as in the continuous model.  The basic spectral properties of the associated Schr{\"o}dinger operator were investigated in our previous work \cite{AndIsoMor}. 
 To study the inverse scattering, we adopt the above third approach.

 A problem arises in the first step where
 we derive the relation between the S-matrix and the D-N map on a finite domain. The method in the continuous case depends largely on the asymptotic expansion of the form (\ref{S1FarfieldCont}), 
 which follows from the asymptotic expansion of the resolvent $(- \Delta - \lambda \mp i0)^{-1}$ 
 at space infinity. However, for the lattice Hamiltonians, we cannot expect it.
 In fact, the usual way to derive this sort of expansion is to apply the stationary phase method to an integral on the Fermi surface. It requires that the Gaussian curvature does not vanish, which can be expected only on restricted regions of the energy.
 In \cite{AndIsoMor},  to study the spectral properties of the lattice Hamiltonian, we passed it on the flat torus ${\bf R}^d/(2\pi{\bf Z})^d$, and instead of the spatial asymptotics of the resolvent, we studied the singularity expansion of the resolvent of the transformed Hamiltonian on the torus. This makes it possible to obtain an analogue of the expansion (\ref{S1FarfieldCont}) in terms of the singularities of the resolvent and to derive the desired relation between the S-matrix and the D-N map in the bounded domain. This S-matrix coincides with $S(\lambda)$ appearing in the time-dependent picture (\ref{S1Smatrixformula}), and is equal to the one defined through the spatial asymptotics  when the Gaussian curvature of the Fermi surface does not vanish.
 Thus,  the forward problem
 can be treated in  a unified framework encompassing the examples such as  square, triangular, hexagonal, diamond, kagome, subdivision lattices, as well as ladder and graphite.
 
We are then led to a boundary value problem on a finite graph for the Schr{\"o}dinger operator $- \widehat{\Delta} + \widehat V$ or the conductivity operator. Let us consider the latter :
\begin{equation}
\left\{
\begin{split}
& \widehat\Delta_{\gamma}\widehat u := \sum_{w \in \mathcal N_v}\gamma(e_{vw})\big(\widehat u(w)- \widehat u(v)\big) = 0, \quad v \in \stackrel{\circ}{\mathcal V_{int}}, \\
& \widehat u(v) = \widehat f(v) , \quad v \in \partial{\mathcal V_{int}},
\end{split}
\right.
\label{S1InvBVPDiscret}
\end{equation}
where $\gamma(e_{vw})>0$ is a conductance of the edge $e_{vw}$ with end points $v, w$. Precise definitions will be explained in \S 2 and \S 7.
The D-N map for (\ref{S1InvBVPDiscret}) is defined in a manner similar to  (\ref{S1DNmapCont}). 
A remarkable fact is that the inverse problem for the {\it network problem} (\ref{S1InvBVPDiscret}) has already been solved in a satisfactory way. One knows
\begin{itemize}
\item uniqueness of the map $\gamma \to \Lambda_{\gamma}$
\item characterization of the D-N map $\Lambda_{\gamma}$
\item algorithm for the reconstruction of $\gamma$ from $\Lambda_{\gamma}$
\item stability of the map $\gamma \to \Lambda_{\gamma}$
\item reconstruction procedure of the  graph from $\Lambda_{\gamma}$
\end{itemize}
by the works of Curtis, Ingerman, Mooers, Morrow, and Colin de Verdi{\'e}re, Gitler, Vertigan (see 
\cite{CuMo90}, \cite{CuMo91}, \cite{CuMoMo94}, \cite{Col94}, \cite{ColGit96}, \cite{CuInMo98}, \cite{CuMo00}). These results enable us to recover the perturbation term (conductance or scalar potential) and also the 
graph structure. We can then solve the inverse problem starting from the scattering matrix.  


\subsection{Main results}
The main assumptions are (A-1) $\sim$ (A-4), (B-1) $\sim$ (B-4) in \S 2. 
The principal results of this paper are as follows.

\begin{itemize}
\item Theorem \ref{ThSmatrixDNmapEquiv} proves that the S-matrix and the D-N map  determine each other. 

\item In \S 6, we show a reconstruction algorithm for the scalar potential from the D-N map of the finite hexagonal lattice.

\item In Subsection 7.2, we discuss how the resistor network is reconstructed from the S-matrix up to some equivalence.

\item Theorem \ref{Theorem7.5DefectProbe} guarantees that in principle it is possible to probe the defects in the periodic structure from the knowledge of the S-matrix.

\item Theorem \ref{ProbeTheoremHoneycombDefectslambda} gives an algorithm to detect the location of defects forming a finite number of holes of the shape of convex polygons in the hexagonal lattice. 

\end{itemize}

Until the end of \S 5, we deal with a general class of lattices satisfying the assumptions (A-1) $\sim$ (A-4) and (B-1) $\sim$ (B-4). 
 As will be seen from our argument, inverse scattering for the resistor network can be formulated and discussed on square, triangular, $d$-dimensional diamond lattices ($ d\geq 2$), ladder of $d$-dimensional square lattices and graphite. 
To find the location of defects  on the hexagonal lattice, in Subsection \ref{subsectionprovingwaves}, we use a special type of solution to the Schr{\"o}dinger equation, which vanishes in a half space in ${\bf Z}^d$ and growing in the opposite half space. 
This is an analogue of exponentially growing solutions for Schr{\"o}dinger operators in the continuous model. Note that  Ikehata \cite{Ikehata99}, \cite{Ikehata04} 
developed the enclosure method to find locations of inclusions by using exponentially growing solutions for the case of  continuous model. 
 Our detection procedure depends largely on the geometric structure of the lattice and should be checked separately for each lattice. Hence we formulate Theorems \ref{Theorem7.5DefectProbe} and \ref{ProbeTheoremHoneycombDefectslambda} only for the hexagonal lattice. 
The square and triangular lattices are dealt with similarly by our theory. 
However, the inverse scattering by defects
 for the higher dimensional diamond lattice,  ladder, graphite, subdivision and kagome lattice is still an open problem, although the forward problem is settled.


\subsection{Plan of the paper}
In \S 2, we recall basic facts on the spectral properties of periodic lattices proved in \cite{AndIsoMor}. The results are extended in \S 3 to the boundary value problem in an exterior domain. In \S 4, the S-matrix and the D-N map in the interior domain are shown to be equivalent. Our S-matrix is derived from the singularity expansion of solutions to the Helmholtz equation. In some energy region, it coincides with the usual S-matrix obtained from the asymptotic expansion at infinity of solutions to the Schr{\"o}dinger equation in the lattice space. This is  proven in \S 5. 
The remaining  sections are devoted to the reconstruction procedure. In \S 6, we reconstruct the scalar potential from the D-N map. In \S 7, we study the reconstruction of the graph structure as a network problem. 
Picking up the example of hexagonal lattice, we also study the probing problem for the location of defects from the S-matrix. 


\subsection{Related works}
There is an extensive literature on the mathematical theory of graphs and their spectra.
We cite here only the articles which have close relations to this paper, but are not mentioned above. For a general survey, see e.g. \cite{MoWo} and the references therein.

For the foundations of the properties of graph Laplacian, see  \cite{Chung} and \cite{Col98}. A general approach to the spectral properties of periodic systems in terms of Mourre's commutator analysis is given in  \cite{GeNi98}. The Floquet-Bloch theory for periodic differential operators is generalized to  more general covering graphs in \cite{Suna90}, \cite{KoOnoSuna}.  Random walk is often used to study the structure of the graph, see e.g. \cite{Do84} and \cite{KoShiSu98}.  Determination of spectra, spectral gap, and (non) existence of eigenvalues are basic issues for periodic, or more generally, covering graphs, and many works are now presented, e.g. \cite{HiShi04}, \cite{HiNo09}, \cite{HSSS12}, \cite{KoSa13}, \cite{S99}, \cite{Su13}.

The inverse scattering for the multi-dimensional discrete Schr{\"o}dinger operator was first reported in \cite{Es}. In \cite{IsKo}, it was proven by using the complex Born approximation of the scattering amplitude. The extension to the hexagonal lattice is done in  \cite{Ando13}.  
The Rellich type theorem for the uniqueness of solutions to the Helmholtz equation, proved in \cite{IsMo1}, plays an essential role in this paper. The Hilbert Nullstellensatz is used in the proof and this idea goes back to Shaban-Vainberg 
\cite{Sha}.
The long-range scattering is discussed in \cite{Na14}. 

For the recent issues on discretization of Riemannian manifolds and their spectral properties, see \cite{BuIvKu13} and the references therein. 
	The monograph \cite{Yaf} contains an exposition of the spectral theory due to \cite{AgHo76}, over  which leans the method of this paper.

In physical literatures, the 2-dimensional Dirac operator is usually adopted as a mathematical model for the graphene (see e.g. \cite{GoGuVo}, \cite{NeGuPeNov}, 
\cite{CuSie}). 
Therefore, our discrete Laplacian on the hexagonal lattice is regarded as a discretization of this Dirac operator. For the mathematical model of carbon nano-tube, see 
\cite{KoKut}, \cite{KuPo}.
An experimental 
result for the defects in graphite is seen in  \cite{KCSSSHSOTTN12}.


\subsection{Basic notation}
For $f \in \mathcal S'({\bf R}^d)$, 
$\widetilde f(\xi)$ denotes its Fourier transform 
\begin{equation}
\widetilde f(\xi) = (2\pi)^{-d/2}\int_{{\bf R}^d}e^{-ix\cdot\xi}
f(x)dx, \quad \xi \in {\bf R}^d,
\label{S1Fouriertransf}
\end{equation}
while for $f(x) \in \mathcal D'({\bf T}^d)$, $\widehat f(n)$ denotes its Fourier coefficients
\begin{equation}
\widehat f(n) = (2\pi)^{-d/2}\int_{{\bf T}^d}e^{-ix\cdot n}f(x)dx, \quad n \in {\bf Z}^d.
\label{S1Fouriercoeffi}
\end{equation}
We also use $\widehat f = \big(\widehat f(n)\big)_{n\in{\bf Z}^d}$ to denote a function on ${\bf Z}^d$, and by $\mathcal U$ the operator
\begin{equation}
\mathcal U : \mathcal D'({\bf Z}^d) \ni \big(\widehat f(n)\big)_{n\in{\bf Z}^d} \to 
f(x) = (2\pi)^{-d/2}\sum_{n\in{\bf Z}^d}\widehat f(n)e^{in\cdot x} \in 
\mathcal D'({\bf T}^d).
\label{S1Fourierseries}
\end{equation}
For Banach spaces $X$ and $Y$, ${\bf B}(X;Y)$ denotes the set of all bounded operators from $X$ to $Y$. For a self-adjoint operator $A$, $\sigma(A), \sigma_p(A), \sigma_d(A), \sigma_e(A)$  denote its spectrum, point spectrum, discrete spectrum 
 and essential spectrum, respectively. $\mathcal H_{ac}(A)$ is the absolutely continuous subspace for $A$, and $\mathcal H_{p}(A)$ is the closure of the linear hull of eigenvectors of $A$.
For an interval $I \subset {\bf R}$ and a Hilbert space ${\bf h}$, $L^2(I,{\bf h},\rho(\lambda)d\lambda)$ denotes the set of all ${\bf  h}$-valued $L^2$-functions on $I$ with respect to the measure $\rho(\lambda)d\lambda$. $S^m = S^m_{1,0}$ denotes the standard H{\"o}rmander class of symbols for pseudo-differential operators ($\Psi$DO), i.e. $|\partial_x^{\alpha}\partial_{\xi}^{\beta}p(x,\xi)| \leq C_{\alpha\beta}(1 + |\xi|)^{m-\beta}$ (\cite{HoVol3}).


\section{Basic properties of graph} \label{Basic}


\subsection{Vertices and edges}
Our object is an infinite, simple (i.e. without self-loop and multiple edge) graph $\Gamma = \{\mathcal V,\mathcal E\}$, where $\mathcal V$ is a vertex set, and $\mathcal E$ is an edge set. For two vertices $v$ and $w$, $v \sim w$ means that they are the end points of an edge $e \in \mathcal E$. We denote it $e = e(v,w)$, and also
$$
o(e) = v, \quad t(e) = w.
$$
However, we do not assume the orientation for the edge. 
The graph $\Gamma$ is assumed to be {\it connected}, i.e.  for any $v, w \in \mathcal V$, there exist $v_1, \cdots, v_m \in \mathcal V$ such that $v=v_1,  v_m = w$ and
$v_i \sim v_{i+1}$, $1 \leq i \leq m-1$. For $v \in \mathcal V$,  we put
\begin{equation}
{\mathcal N}_v = \{w \in \mathcal V \, ; \, v \sim w\},
\label{C1S1Adjacentpoints}
\end{equation}
and call it the set of points {\it adjacent} to $v$.
The {\it degree} of $v \in \mathcal V$ is then defined by
\begin{equation}
{\rm deg}\, (v) = {\sharp}\, {\mathcal N}_v = {\sharp}\,\{e \in \mathcal E\, ; \, o(e) = v\},
\nonumber
\end{equation}
which is assumed to be finite for all $v \in \mathcal V$.
Let $\ell^2(\mathcal V)$ be the set of ${\bf C}$-valued functions $f = \left(f(v)\right)_{v\in\mathcal V}$ on $\mathcal V$ satisfying
\begin{equation}
\|f\|^2 := {\mathop\sum_{v\in\mathcal V}}|f(v)|^2 \,{\rm deg}\,(v) < \infty,
\nonumber
\end{equation}
which is a Hilbert space
equipped with the inner product
\begin{equation}
(f,g)_{deg} = {\mathop\sum_{v\in\mathcal V}}f(v)\overline{g(v)}\,{\rm deg}\,(v).
\label{Innerprductgeneral}
\end{equation}
The Laplacian $\widehat{\Delta}_{\Gamma}$ on the graph $\Gamma = \{\mathcal V, \mathcal E\}$ is defined by
\begin{equation}
\big(\widehat{\Delta}_{\Gamma}\widehat f\big)(v) = \frac{1}{\deg\,(v)}\sum_{w\sim v}\widehat f(w),
\label{S2LaplacianGeneral}
\end{equation}
which is self-adjoint on $\ell^2(\mathcal V)$.

A subset $\Omega \subset \mathcal V$ is \textit{connected} if, for any $v,w\in \Omega $, there exist $v=v_1, v_2 , \cdots , v_m =w \in \Omega $ such that $ v_i \sim v_{i+1} $, $ 1\leq i \leq m-1 $. For $v \in \mathcal V$, $v \sim \Omega$ means that there exists  $w \in \Omega$ such that $v \sim w$.
For a connected subset $\Omega \subset \mathcal{V} $, we define
\begin{equation}
\Omega ' = \{ v\not\in \Omega \ ; \ v \sim \Omega \},
\label{S2_def_Omegaprime}
\end{equation}
and put $D= \Omega \cup \Omega ' $.
For this set $D$, we put 
\begin{equation}
 \stackrel{\circ}D \, = \Omega ,\quad \partial D = \Omega ' . 
 \label{S2_def_dint_partiald}
 \end{equation}
We call $\stackrel{\circ}D$ the {\it interior} of $D$ and $\partial D$ the {\it boundary} of $D$.

We define
\begin{equation}
{\rm deg}_D(v) =\left\{ \begin{split}
& \sharp\{w \in D \, ; \,  v \sim w\}, \quad v \in \, \stackrel{\circ}D, \\
& \sharp \{ w\in \,  \stackrel{\circ}D \ ; \ v \sim w \} , \quad v \in \partial D .
 \end{split} \right.
\label{definedegDv}
\end{equation}
The normal derivative at the boundary $\partial D$ is defined by
\begin{equation}
(\partial_{\nu}^D\widehat f)(v) = - \frac{1}{\mathrm{deg} _D (v) }\sum_{w\in\stackrel{\circ}D, w \sim v}\widehat f(w), \quad 
v \in \partial D.
\label{normalder.general}
\end{equation}
Then the following Green's formula holds
\begin{equation}
\big(\widehat{\Delta}_{\Gamma}\widehat f,\widehat g\big)_{\ell^2(\stackrel{\circ}D)}
- \big(\widehat f,\widehat{\Delta}_{\Gamma}\widehat g\big)_{\ell^2(\stackrel{\circ}D)} = 
\big(\partial_{\nu}^D\widehat f,\widehat g\big)_{\ell^2(\partial D)} - 
\big(\widehat f,\partial_{\nu}^D\widehat g\big)_{\ell^2(\partial D)}
\label{GreenGeneral}
\end{equation}
for $\widehat f, \widehat g \in \ell^2(\mathcal V)$ such that 
$\widehat f(v) = \widehat g(v)=0$ if $v \not\in D$. 
Note that the inner product on $\partial D$ is defined by
\begin{equation}
\big(\widehat f,\widehat g\big)_{\ell^2(\partial D)} = \sum_{w\in \partial D}
\widehat f(w)\overline{\widehat g(w)} \mathrm{deg}_D (w) ,
\label{InnerproductpartialD}
\end{equation} %
and the sum in the inner product of the left-hand side of  (\ref{GreenGeneral}) ranges over the points in $\stackrel{\circ}D$.


\subsection{Laplacian on the perturbed periodic graph}\label{subsectionperiodicgraph}
A periodic graph in ${\bf R}^d$ is a triple $\Gamma_0 = \{\mathcal L_0, \mathcal V_0, \mathcal E_0\}$, where  $\mathcal L_0$ is a lattice of rank $d$ in ${\bf R}^d$ with basis ${\bf v}_j, j= 1,\cdots,d$, i.e. 
\begin{equation}
\mathcal L_0 = \big\{{\bf v}(n)\, ; \, n \in 
{\bf Z}^d\big\}, \quad
{\bf v}(n) = \sum_{j=1}^dn_j{\bf v}_j, \quad n =(n_1,\cdots,n_d) \in {\bf Z}^d,
\nonumber
\end{equation}
and the vertex set is defined by
\begin{equation}
\mathcal V_0 = {\mathop\cup_{j=1}^s}\big(p_j + \mathcal L_0\big),
\nonumber
\end{equation}
and where $p_j$, $j = 1,\cdots,s$, are the points in ${\bf R}^d$ satisfying
\begin{equation}
p_i - p_j \not\in \mathcal L_0, \quad {\rm if}\quad i\neq j.
\label{C1pi-pj}
\end{equation}
By (\ref{C1pi-pj}), there exists a bijection $\mathcal V_0 \ni a \to (j(a),n(a)) \in \{1,\cdots,s\}\times{\bf Z}^d$ such that
\begin{equation}
a = p_{j(a)} + {\bf v}(n(a)).
\label{a=pja+vna}
\end{equation}
In the following, we often identify $a$ with $(j(a),n(a))$. 
The group ${\bf Z}^d$ acts on ${\mathcal V}_0$ as follows :
\begin{equation}
{\bf Z}^d\times{\mathcal V}_0 \ni (m,a) \to m\cdot a := p_{j(a)}+{\bf v}(m+n(a)) \in {\mathcal V}_0.
\label{Zdaction}
\end{equation}
The edge set  $\mathcal E_0 \subset \mathcal V_0 \times \mathcal V_0$ is assumed to satisfy
\begin{equation}
\mathcal E_0 \ni (a,b) \Longrightarrow  
(m\cdot a, m\cdot b) \in \mathcal E_0,
\quad \forall m \in {\bf Z}^d.
\nonumber
\end{equation}
Then ${\rm deg}\,(p_j + {\bf v}(n))$ depends only on $j$, and is denoted by ${\rm deg}_0(j)$ :
\begin{equation}
{\rm deg}_0(j) ={\rm  deg}\,(p_j + {\bf v}(n)).
\label{S2Definedeg(j)}
\end{equation}
Any function $\widehat f$ on $\mathcal V_0$ is written as
$\widehat f(n) = (\widehat f_1(n),\cdots,\widehat f_s(n)), \  n \in {\bf Z}^d$, where $\widehat f_j(n)$ is 
 a function on $p_j + \mathcal L_0$.  Hence $\ell^2(\mathcal V_0)$ is a Hilbert space equipped with the inner product
\begin{equation}
(\widehat f,\widehat g)_{\ell^2(\mathcal V_0)} = 
\sum_{j=1}^s(\widehat f_j,\widehat g_j)_{{\rm deg}_0(j)}.
\label{wholespaceinnerproduct}
\end{equation}
We then define a unitary operator $\mathcal U_{\mathcal L_0} : \ell^2(\mathcal V_0) \to L^2({\bf T}^d)^s$ by
\begin{equation}
\big(\mathcal U_{\mathcal L_0}\widehat f\big)_j = (2\pi)^{-d/2}
\sqrt{{\rm deg}_0(j)}\sum_{n\in {\bf Z}^d}\widehat f_j(n)e^{in\cdot x},
\label{S1UdDefine}
\end{equation}
where $L^2({\bf T}^d)^s$ is equipped with the inner product
\begin{equation}
(f,g)_{L^2({\bf T}^d)^s} = \sum_{j=1}^s\int_{{\bf T}^d}f_j(x)\overline{g_j(x)}dx.
\label{S2L2bfTdinnerproduct}
\end{equation}

Recall that the shift operator $S_j$ acts on a sequence $\big(\widehat f(n)\big)_{n\in{\bf Z}^d}$ as follows : 
$$
\big(S_j\widehat f\big)(n) = \widehat f(n + {\bf e}_j),
$$
where ${\bf e}_1 = (1,0,\cdots,0), \cdots, {\bf e}_d= (0,\cdots,0,1)$. Then we have 
\begin{equation}
\mathcal U_{\mathcal L_0} S_j = e^{-ix_j}\mathcal U_{\mathcal L_0}.
\label{S2Ud0SjUd0}
\end{equation}

The Laplacian $\widehat \Delta_{\Gamma_0}$ on the graph $\Gamma_0$ is defined by the  formula
\begin{equation}
\begin{split}
  (\widehat\Delta_{\Gamma_0} \widehat f)(n) & = (\widehat g_1(n),\cdots,\widehat g_s(n)),\\
\widehat g_i(n) & = \frac{1}{{\rm deg}_0(i)}\sum_{b\sim p_i + {\bf v}(n)}\widehat f_{j(b)}(n(b)),
\end{split}
\label{S2DefineLaplacian}
\end{equation}
where $b = p_{j(b)} + {\bf v}(n(b))$. Recalling (\ref{a=pja+vna}), we can rewrite it as 
\begin{equation}
\widehat g_i(n) = \frac{1}{{\rm deg}_0(i)}\sum_{(i,n)\sim(j,n')}\widehat f_{j}(n').
\label{S2DefineLaplacian'}
\end{equation}
Passing to the Fourier series, (\ref{S2DefineLaplacian}) has the following form :
\begin{equation}
\mathcal U_{\mathcal L_0}(- \widehat\Delta_{\Gamma_0}) (\mathcal U_{\mathcal L_0})^{-1} f = H_0(x)f(x), \quad 
f \in L^2({\bf T}^d)^s,
\nonumber
\end{equation}
where $H_0(x)$ is an $s\times s$ Hermitian matrix whose entries are trigonometric functions. Let $D$ be the $s\times s$ diagonal matrix whose $(j,j)$ entry is $\sqrt{{\rm deg}_0(j)}$. Then $\mathcal U_{\mathcal L_0} = D\mathcal U$, where $\mathcal U$ means the operator $(\widehat f_1,\cdots,\widehat f_s) \to (\mathcal U\widehat f_1,\cdots,\mathcal U\widehat f_s)$ (see (\ref{S1Fourierseries})), 
 hence
\begin{equation}
H_0(x) = DH_0^0(x)D^{-1}, \quad H_0^0(x) = \mathcal U(- \widehat\Delta_{\Gamma_0}){\mathcal U}^{-1},
\label{S2H0(x)rewritten}
\end{equation}
and $H_0^0(x)$ is computed by (\ref{S1Fourierseries}) and (\ref{S2Ud0SjUd0}).

Let $\mathcal H_0 = L^2\left({\bf T}^d\right)^s$  equipped with the inner product (\ref{S2L2bfTdinnerproduct}).
 Then, the operator of multiplication by $H_0(x)$ is a bounded self-adjoint operator on $\mathcal H_0$, which is denoted by $H_0$.
Let $\lambda_1(x) \leq \lambda_2(x) \leq \cdots \leq \lambda_s(x)$ be the eigenvalues of $H_0(x)$, and 
\begin{equation}
M_{\lambda, j} = \{x \in {\bf T}^d\, ; \, \lambda_j(x) = \lambda\}.
\label{C4S1MlambdajDefine}
\end{equation}
Then we have
\begin{equation}
p(x,\lambda) := \det\big(H_0(x)-\lambda\big) = \prod_{j=1}^s(\lambda_j(x) - \lambda),
\label{C4S1pxlambda=prodpj}
\end{equation}
\begin{equation}
M_{\lambda}:= \{x \in {\bf T}^d\, ; \, p(x,\lambda)=0\} = {\mathop\cup_{j=1}^s}M_{\lambda,j}.
\label{C4S1Mlambda=cupMlambdaj}
\end{equation}

Let 
\begin{equation}
{\bf T}^d_{\bf C} = {\bf C}^d/(2\pi{\bf Z})^d, \quad
M_{\lambda}^{\bf C} = \{z \in {\bf T}^d_{\bf C}\, ; \, p(z,\lambda)=0\},
\end{equation}
\begin{equation}
M^{\bf C}_{\lambda,reg} = \{z \in M_{\lambda}^{\bf C}\, ; \, \nabla_z p(z,\lambda) \neq 0\},
\end{equation}
\begin{equation}
M^{\bf C}_{\lambda,sng} = \{z \in M_{\lambda}^{\bf C}\, ; \, \nabla_z p(z,\lambda) = 0\}.
\end{equation}


\subsection{Assumptions}
The following assumptions are imposed on the free system.

\medskip
\noindent
{\bf(A-1)}  {\it There exists a subset $\mathcal T_1 \subset \sigma(H_0)$ such that for} $\lambda \in \sigma(H_0)\setminus \mathcal T_1$, 

\smallskip
{\bf(A-1-1)} $M_{\lambda,sng}^{\bf C}$ {\it is discrete}.

\smallskip
{\bf(A-1-2)} {\it Each connected component of $M_{\lambda,reg}^{\bf C}$ intersects with ${\bf T}^d$ and the intersection is a $(d-1)$-dimensional real analytic submanifold of ${\bf T}^d$.}

\medskip
\noindent
{\bf(A-2)} {\it  There exists a finite set $\mathcal T_0 \subset \sigma (H_0)$ such that}
$$
M_{\lambda,i}\cap M_{\lambda,j} = \emptyset, \ \ if \ \ i \neq j, \ \ 
\lambda \in \sigma(H_0)\setminus \mathcal T_0.
$$

\medskip
\noindent
{\bf(A-3)} \ {\it $\nabla_xp(x,\lambda) \neq 0$ , \ on \ $M_{\lambda}$, \ $\lambda \in \sigma(H_0)\setminus \mathcal T_0$.}

\medskip
\noindent
{\bf(A-4)}  {\it The unique continuation property holds for $\widehat H_0$ in $\mathcal V_0$. Namely, 
any $\widehat u$ satisfying $(- \widehat\Delta_{\Gamma_0}-\lambda)\widehat u = 0$ on $\mathcal V_{0}$ except for a finite number of points, where $\lambda$ is a constant,  vanishes identically on $\mathcal V_{0}$.}

\medskip
For the square, triangular, hexagonal, Kagome, diamond lattices and the subdivision of square lattice, $\mathcal T_1$ is a finite set. However, for the  ladder and graphite, $\mathcal T_1$ fills closed intervals.  See \cite{AndIsoMor}, \S 5. 

\medskip
We consider a connected graph $\Gamma = \{\mathcal V, \mathcal E\}$, which is a local perturbation of the periodic lattice $\Gamma_0 = \{\mathcal L_0, \mathcal V_0, \mathcal E_0\}$ having the properties described above. 
We impose the following assumptions on $\Gamma$.

\bigskip
\noindent 
{\bf(B-1)} 
{\it  There exist two subsets $\mathcal V_{int}, \mathcal V_{ext} \subset \mathcal V$ having the following properties :}

\smallskip
{\bf(B-1-1)} $\ \mathcal V = \mathcal V_{int}\cup \mathcal V_{ext}$.

\smallskip
{\bf(B-1-2)} $\ \mathcal V_{int}\cap\mathcal V_{ext} = \partial \mathcal V_{int} = 
\partial \mathcal V_{ext}$.

\smallskip
{\bf(B-1-3)} \ {\it $\ \mathcal V_{int}$, $\mathcal V_{ext}$ are connected.}

\smallskip
{\bf(B-1-4)} \ $\sharp \mathcal V_{int} < \infty$.

\medskip
\noindent
{\bf(B-2)} \ {\it The unique continuation  property holds on ${\mathcal V_{ext}}$.} 

\medskip
\noindent
{\bf(B-3)} \ {\it There exist a subset ${\mathcal V^{(0)}_{ext}} \subset \mathcal V_0$ such that 
$\sharp\left(\mathcal V_0 \setminus {\mathcal V^{(0)}_{ext}}\right) < \infty$ and a bijection 
$\mathcal V_{ext} \to {\mathcal V^{(0)}_{ext}}$ which preserves the edge relation.}

\medskip
Because of (B-3), we identify $\mathcal V_{ext}$ with $\mathcal V_{ext}^{(0)}$ and denote the point in $\mathcal V_{ext}$ as $(j,n)$.

\begin{figure}[htbp]
 \begin{minipage}{0.49\hsize}
  \begin{center}
\includegraphics[width=60mm, bb=0 0 276 237]{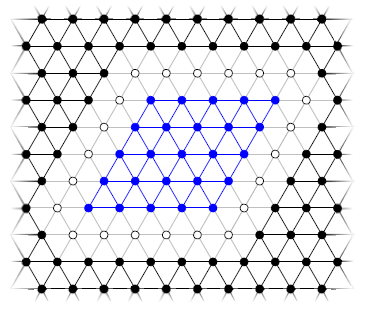} 
 \end{center}
  \caption{Boundary of a domain in the triangular lattice}
  \label{S3_tri_hex_boundary}
 \end{minipage}
 \begin{minipage}{0.49\hsize}
  \begin{center}
\includegraphics[width=55mm, bb=0 0 239 237]{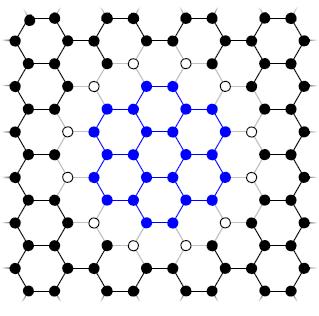}
  \end{center}
  \caption{Boundary of a domain in the hexagonal lattice}
  \label{BoundaryHexagonal}
 \end{minipage}
\end{figure}

\begin{figure}[hbtp]
\centering
\includegraphics[width=70mm, bb=0 0 311 311]{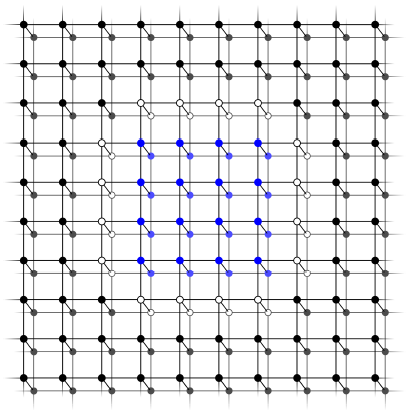}
\caption{Boundary of a domain in the two dimensional square ladder}
\label{BoundarySquareLadder}
\end{figure}

Typical examples of the decomposition $\mathcal V = \mathcal V_{int}\cup \mathcal V_{ext}$ are given in Figures \ref{S3_tri_hex_boundary}, \ref{BoundaryHexagonal} and \ref{BoundarySquareLadder}, where $\Sigma = \partial \mathcal V_{int} = \partial \mathcal V_{ext}$ is the set of the white dots, and $\stackrel{\circ}{\mathcal V_{int}}$, $\stackrel{\circ}{\mathcal V_{ext}}$ are the regions inside $\Sigma$, outside $\Sigma$, respectively.
%
%
%


\begin{lemma}\label{Lemma2.1}
Let $\Sigma = \partial \mathcal V_{int} = \partial \mathcal V_{ext}$.  \\
\noindent
(1) $\ \mathcal V$ is written as a disjoint union : $\ \mathcal V = \stackrel{\circ}{\mathcal V_{int}} \cup \, \Sigma \, \cup \stackrel{\circ}{\mathcal V_{ext}}$. \\\noindent
(2) \ For any $v \in \Sigma$, $v \sim \stackrel{\circ}{\mathcal V_{int}}$ and $v \sim \stackrel{\circ}{\mathcal V_{ext}}$ hold.\\
\noindent
(3)  Any path starting from $\stackrel{\circ}{\mathcal V_{int}}$ and ending in $\stackrel{\circ}{\mathcal V_{ext}}$ passes through $\Sigma$.
\end{lemma}

Proof. By (B-1-2), $\partial \mathcal V_{ext} \subset \mathcal V_{ext}$, hence
$\mathcal V_{ext} = \stackrel{\circ}{\mathcal V_{ext}} \cup \, \partial \mathcal V_{ext}$. Similarly, $\mathcal V_{int} = \stackrel{\circ}{\mathcal V_{int}} \cup \, \partial \mathcal V_{int}$. This and (B-1-1) imply (1). Since $\Sigma = \partial \mathcal V_{int} = \partial \mathcal V_{ext}$, (2) follows. Suppose there exist $v_i \in \stackrel{\circ}{\mathcal V_{int}}$, and $v_e \in \stackrel{\circ}{\mathcal V_{ext}}$ such that $v_i \sim v_e$. Then $v_e \in \partial{\mathcal V_{int}}$. This is in contradiction to (1). \qed
 
 \bigskip
The Hilbert space $\ell^2(\mathcal V)$ then admits an orthogonal decomposition
$$
\ell^2(\mathcal V) = \ell^2(\stackrel{\circ}{\mathcal V_{ext})}\oplus \, \ell^2({\mathcal V_{int}}).
$$
Let $\widehat P_{ext}$ be the associated orthogonal projection :
\begin{equation}
\widehat P_{ext} : \ell^2(\mathcal V) \to \ell^2(\stackrel{\circ}{\mathcal V_{ext}}).
\label{Pextdefine}
\end{equation}
Let ${\widehat \Delta}_{\Gamma}$ be the Laplacian on the graph $\Gamma$. We assume that the perturbation $\widehat V$ has the following property. 

\medskip
\noindent
{\bf(B-4)} \ {\it $\widehat V$ is bounded self-adjoint on $\ell^2(\mathcal V)$ and has  support in $\stackrel{\circ}{\mathcal V_{int}}$, i.e. $\widehat V\widehat u = 0$ on  $\mathcal V_{ext}$, $\forall \widehat u \in \ell^2(\mathcal V)$.}

\bigskip
In \cite{AndIsoMor}, the exterior domain $\mathcal V_{ext}$ was defined in a slightly different, more restricted form. However, all the arguments there work well for the above $\mathcal V_{ext}$   
under the above assumptions (B-1) $\sim$ (B-4).


\subsection{Function spaces}

In \cite{AndIsoMor}, for the periodic graph $\Gamma_0 = \{\mathcal L_0, \mathcal V_0, \mathcal E_0\}$, the spaces $\ell^{2}, \ell^{2,\sigma}, \widehat{\mathcal B}$, $\widehat{\mathcal B}^{\ast}$, $\widehat{\mathcal B}^{\ast}_0$ were defined as the spaces equipped with the following norms : 
\begin{equation}
\|\widehat f\|^2_{\ell^2(\mathcal V_0)} = \sum_{n\in \mathcal V_0} |\widehat f(n)|^2,
\label{l2norm}
\end{equation}
where (see (\ref{wholespaceinnerproduct}))
\begin{equation}
|\widehat f(n)|^2 = \sum_{j=1}^s|\widehat f_j(n)|^2{\rm deg}_0(j),
\label{|f(n)|define}
\end{equation}
\begin{equation}
\|\widehat f\|^2_{\ell^{2,\sigma}(\mathcal V_0)} = \sum_{n\in \mathcal V_0}(1 + |n|^2)^{\sigma} |\widehat f(n)|^2, \quad \sigma \in {\bf R},
\label{l2norm}
\end{equation}
\begin{equation}
\|\widehat f\|_{\widehat B(\mathcal V_0)}^2 = \sum_{j=0}^{\infty}r_j^{1/2}
\Big(\sum_{r_{j-1}\leq |n|<r_j, n \in \mathcal V_0}|\widehat f(n)|^2\Big)^{1/2},
\label{Bspace}
\end{equation}
where $r_{-1}=0, r_j = 2^j \ (j \geq 0)$, 
\begin{equation}
\|\widehat f\|_{\widehat{\mathcal B}^{\ast}(\mathcal V_0)}^2 = \sup_{R>1}\frac{1}{R}\sum_{|n|<R,n\in \mathcal V_0}|\widehat f(n)|^2,
\label{Bastspace}
\end{equation}
\begin{equation}
\widehat{\mathcal B}^{\ast}_0(\mathcal V_0) \ni \widehat f \Longleftrightarrow 
\lim_{R\to\infty}\frac{1}{R}\sum_{|n|<R,n\in \mathcal V_0}|\widehat f(n)|^2=0.
\label{Bast0}
\end{equation}
For the perturbed graph $\Gamma$, these spaces  are defined as above, replacing $\mathcal V_0$ by $\stackrel{\circ}{\mathcal V_{ext}}$ and adding the norm of $\ell^2(\mathcal V_{int})$. They are denoted by $\ell^2(\mathcal V)$, etc, or sometimes $\ell^2$ without fear of confusion.


\subsection{Continuous spectrum and embedded eigenvalues}
\label{SpectralPropertiesWholespace}
We now define the perturbed Hamiltonian $\widehat H$ by
\begin{equation}
\widehat H = - \widehat\Delta_{\Gamma} + \widehat V.
\label{DefinewidehatH}
\end{equation}
 Let us review the  spectral properties of $\widehat H$.

  
\begin{lemma}\label{sigmae(H)=sigma(H0)}
(Theorem 7.1, Lemma 7.2 in \cite{AndIsoMor}). \\
\noindent
(1) \ $\sigma_e(\widehat H) =  \sigma(\widehat H_0)$. \\
\noindent
(2) The eigenvalues of $\widehat H$ in $\sigma_e(\widehat H)\setminus \mathcal T_1$ is finite with finite multiplicities. \\
\noindent
(3) There is no eigenvalue in $\sigma_e(\widehat H)\setminus \mathcal T_1$, provided $\widehat H$ has the unique continuation property in $\mathcal V_{int}$.
\end{lemma}

The assertions (2) and (3) of Lemma \ref{sigmae(H)=sigma(H0)} are based on the following Rellich type theorem:


\begin{theorem}\label{Rellichtypetheorem}
(Theorem 5.1 in \cite{AndIsoMor}). 
Assume (A-1) and $\lambda \in \sigma_e(\widehat H)\setminus \mathcal T_1$. 
If $\widehat u \in \widehat{\mathcal B}^{\ast}_0(\mathcal V_0)$ satisfies
$$
(\widehat H_0 - \lambda)\widehat u(n) = 0, \quad |n| > R_0,
$$
for some $R_0 > 0$, then there exists $R > R_0$ such that $\widehat u(n) = 0$ for $|n| > R$.
\end{theorem}

Theorem \ref{Rellichtypetheorem}, for which the assumption (A-1) is essential, plays also an important role in the inverse  scattering procedure to be developed in \S 4. Note, however, by the well-known perturbation theory for the continuous spectrum by Agmon, Kato-Kuroda (see \cite{Agmon75}, \cite{Kuroda}), one can prove the discreteness of embedded eigenvalues, for which we can avoid (A-1), and construct  the spectral representation and S-matrix outside the embedded eigenvalues. This was already done in \S 7 of \cite{AndIsoMor}. 
In this paper, we always assume (A-1).


\subsection{Radiation condition}
The well-known radiation condition of Sommerfeld is extended to the discrete Schr{\"o}dinger operator in the following way.

For $u  \in \mathcal S'({\bf R}^d)$, the wave front set $WF^{\ast}(u)$ is defined as follows. 
For $(x_0,\omega) \in {\bf R}^d\times S^{d-1}$, $(x_0,\omega) \not\in WF^{\ast}(u)$, if there exist $0 < \delta < 1$ and $\chi \in C_0^{\infty}({\bf R}^d)$ such that $\chi(x_0)=1$ and
\begin{equation}
\lim_{R\to\infty}\frac{1}{R}\int_{|\xi|<R}|C_{\omega,\delta}(\xi)(\widetilde{ \chi u})(\xi)|^2d\xi = 0,
\end{equation}
where $C_{\omega,\delta}(\xi)$ is the characteristic function of the cone $\{\xi \in {\bf R}^d\, ; \, \omega\cdot\xi > \delta|\xi|\}$.

We consider distributions  on the torus ${\bf T}^d = {\bf R}^d/(2\pi{\bf Z})^d$. 
By using the Fourier series, the counter parts of 
the spaces $\ell^2(\mathcal V_0)$, 
$\ell^{2,\sigma}(\mathcal V_0), \widehat{\mathcal B}(\mathcal V_0), \widehat{\mathcal B}^{\ast}(\mathcal V_0)$ and  $\widehat{\mathcal B}^{\ast}_0(\mathcal V_0)$ are naturally defined on ${\bf T}^d$, which are denoted by $L^2({\bf T}^d), H^{\sigma}({\bf T}^d), \mathcal B({\bf T}^d), \mathcal B^{\ast}({\bf T}^d)$ and $\mathcal B^{\ast}_0({\bf T}^d)$, respectively. As above, we often omit ${\bf T}^d$.

Let $H_0(x)$ be the matrix in (\ref{S2H0(x)rewritten}), and $\lambda_j(x), j = 1,\cdots, s$,  be its eigenvalues. By (A-2) and (A-3), if  $\lambda \in \sigma_e(\widehat H)\setminus \mathcal T_0$, they are simple, and non-characteristic, i.e. $\nabla_x\lambda_j(x) \neq 0$ on $M_{\lambda,j}$. Let $P_j(x)$ be the eigenprojection associated with $\lambda_j(x)$. Then, in a small neighborhood of 
$M_{\lambda,j}$, $P_j(x)$ is smooth with respect to $x$. 
Suppose $u \in \mathcal B^{\ast}$ satisfies the equation
\begin{equation}
(H_0(x) - \lambda)u = f \in \mathcal B, \quad {\rm on} \quad {\bf T}^d. 
\label{H0X-lambdau0fonTd}
\end{equation}
Then, outside $M_{\lambda}$, $u$ is in $\mathcal B$. Therefore, when we talk about $WF^{\ast}(u)$, we have only to localize it in a small neighborhood of $M_{\lambda}$.
Now, the solution $u$ of the equation (\ref{H0X-lambdau0fonTd}) is said to satisfy the {\it outgoing radiation condition} if
\begin{equation}
WF^{\ast}(P_ju) \subset \{(x,\omega_x)\, ;\, x \in M_{\lambda,j}\}, \quad 1 \leq j \leq s,
\end{equation}
where $\omega_x$ is the unit normal of $M_{\lambda,j}$ at $x$ such that $\omega_x\cdot\nabla\lambda_j(x)<0$. Strictly speaking, one must multiply a cut-off function near $M_{\lambda,j}$ to $u$, which is omitted for the sake of simplicity. Similarly, $u$ is said to satisfy the {\it incoming radiation condition}, if
\begin{equation}
WF^{\ast}(P_ju) \subset \{(x,-\omega_x)\, ;\, x \in M_{\lambda,j}\}, \quad 1 \leq j \leq s.
\end{equation}

We return to the equation on the perturbed lattice
\begin{equation}
(\widehat H - \lambda)\widehat u = \widehat f \quad {\rm on} \quad \mathcal V.
\label{EqonLattice}
\end{equation}
We say that $\widehat u$ satisfies the outgoing (incoming) radiation condition if $\mathcal U\widehat P_{ext}\widehat u$ is outgoing (incoming),
where $\mathcal U$ and $\widehat P_{ext}$ are defined by (\ref{S1Fourierseries}) and (\ref{Pextdefine}).

Let $\widehat R(z) = (\widehat H - z)^{-1}$, and put
\begin{equation}
\mathcal T = \mathcal T_0 \cup \mathcal T_1 \cup \mathcal \sigma_p(\widehat H).
\label{Thresholds}
\end{equation}


\begin{lemma}\label{RadCondUniquelemma}
(Lemma 7.6 in \cite{AndIsoMor}).
If $\lambda \in \sigma_e(\widehat H)\setminus \mathcal T$, 
the solution of the equation (\ref{EqonLattice}) satisfying the outgoing or incoming radiation condition is unique.
\end{lemma}


\begin{theorem}\label{ThLAP}
(Theorem 7.7 in \cite{AndIsoMor}).
Take any compact set $I \subset \sigma_e(\widehat H)\setminus\mathcal T$,  and $\lambda \in I$. Then for any $\widehat f, \widehat g \in \widehat{\mathcal B}$, there exists a limit
$$
\lim_{\epsilon \to 0}(\widehat R(\lambda \pm i0)\widehat f,\widehat g) = 
(\widehat R(\lambda \pm i0)\widehat f,\widehat g).
$$
Moreover, there exists a constant $C > 0$ such that
$$
\|\widehat R(\lambda \pm i0)\widehat f\|_{\widehat{\mathcal B}^{\ast}} \leq C\|\widehat f\|_{\widehat{\mathcal B}}, \quad \lambda \in I.
$$
For $\widehat f \in \widehat{\mathcal B}$, $\widehat R(\lambda + i0)\widehat f$ satisfies the outgoing radiation condition, and $\widehat R(\lambda - i0)\widehat f$ satisfies the incoming radiation condition. Moreover, letting
\begin{equation}
\widehat Q_1(z) = (\widehat H_0 - z)\widehat P_{ext}\widehat R(z),
\label{DrfineQ1(z)}
\end{equation}
\begin{equation}
Q_1(\lambda \pm i0) = \mathcal U_{\mathcal L_0}\widehat Q_1(\lambda \pm i0),
\label{Q1(lambdapmi0)define}
\end{equation}
and $u_{\pm} = \mathcal U_{\mathcal L_0}\widehat P_{ext}\widehat R(\lambda \pm i0)\widehat f$, we have
\begin{equation}
P_ju_{\pm} \mp \frac{1}{\lambda_j(x) - \lambda \mp i0}\otimes\left(P_jQ_1(\lambda \pm i0)\widehat f\right)\Big|_{M_{\lambda,j}} \in \mathcal B^{\ast}_0.
\label{Pju-outinis iBast0}
\end{equation}
\end{theorem}


\subsection{Spectral representation}
For the case of $\widehat H_0$, the spectral representation means the diagonalization of the matrix $H_0(x)$. We first prepare its representation space.
Take an eigenvector $a_j(x) \in {\bf C}^s$ of $H_0(x)$ satisfying $H_0(x)a_j(x) = \lambda_j(x)a_j(x)$, $|a_j(x)|=1$.  Let ${\bf h}_{\lambda,j}$ be the Hilbert space of ${\bf C}$-valued functions on $M_{\lambda,j}$ equipped with the inner product
$$
(\phi,\psi) = \int_{M_{\lambda,j}}\phi(x)\overline{\psi(x)}\frac{dM_{\lambda,j}}{|\nabla \lambda_j(x)|}.
$$
Put
\begin{equation}
I_j = \{\lambda_j(x) \, ; \, x \in {\bf T}^d\}\setminus \mathcal T,
\end{equation}
\begin{equation}
I = {\mathop\cup_{j=1}^s}I_j = \sigma(H_0)\setminus \mathcal T,
\end{equation}
\begin{equation}
{\bf H}_j = L^2(I_j,{\bf h}_{\lambda,j}a_j,d\lambda).
\end{equation}
We define ${\bf h}_{\lambda,j}$ and ${\bf H}_j$ to be $\{0\}$ for $\lambda \in I\setminus I_j$, and put
\begin{equation}
{\bf h}_{\lambda} = {\bf h}_{\lambda,1}a_1\oplus \cdots\oplus {\bf h}_{\lambda,s}a_s,
\end{equation}
\begin{equation}
{\bf H} = {\bf H}_{1}\oplus \cdots\oplus {\bf H}_{s} = L^2(I,{\bf h}_{\lambda},d\lambda).
\end{equation}
For $f \in \mathcal B({\bf T}^d)$, we put
\begin{equation}
\mathcal F_{0,j}(\lambda) f = 
\left\{
\begin{split}
& P_j(x)f(x)\Big|_{M_{\lambda,j}}, \quad {\rm if} \quad \lambda \in I_j, \\
& 0, \quad {\rm otherwise},
\end{split}
\right.
\label{Pjf(x)define}
\end{equation}
\begin{equation}
\mathcal F_0(\lambda)f = (\mathcal F_{0,1}(\lambda)f,\cdots,\mathcal F_{0,s}(\lambda)f).
\end{equation}
\begin{equation}
\widehat{\mathcal F}_0(\lambda) = \mathcal F_0(\lambda)\mathcal U_{\mathcal L_0}.
\label{DefineF0labda=F0lambdaUL0}
\end{equation}

For the perturbed lattice, we define
\begin{equation}
\widehat{\mathcal F}^{(\pm)}(\lambda) = 
\widehat{\mathcal F}_0(\lambda)\widehat Q_1(\lambda \pm i0).
\label{DefineFpmlambda}
\end{equation}
Note that this is denoted by $\widehat{\mathcal F}_{\pm}(\lambda)$ in \cite{AndIsoMor}.
Then for any compact set $J \subset \sigma_e(\widehat H)\setminus{\mathcal T}$, there exists a constant $C > 0$ such that
\begin{equation}
\|\widehat{\mathcal F}^{(\pm)}(\lambda)\widehat f\|_{{\bf h}_{\lambda}} \leq C\|\widehat f\|_{\widehat{\mathcal B}}, \quad \lambda \in J.
\end{equation}
We define
\begin{equation}
\big(\widehat{\mathcal F}^{(\pm)}\widehat f\big)(\lambda) = 
\widehat{\mathcal F}^{(\pm)}(\lambda)\widehat f.
\end{equation}
Similarly,
\begin{equation}
\big(\widehat{\mathcal F}_{0}\widehat f\big)(\lambda) = \widehat{\mathcal F}_{0}(\lambda)\widehat f, \quad 
\big({\mathcal F}_{0} f\big)(\lambda) = {\mathcal F}_{0}(\lambda)f.
\end{equation}
Let $\widehat E(\cdot)$ be the resolution of the identity for $\widehat H$.


\begin{theorem}
(Theorem 7.11 in \cite{AndIsoMor}).\\
\noindent
(1) $\ \widehat{\mathcal F}^{(\pm)}$ is uniquely extended to a partial isometry with initial set $\mathcal H_{ac}(\widehat H) = \widehat E(I)\ell^2(\mathcal V)$ and final set $\bf H$. \\
\noindent
(2) $\ \big(\widehat{\mathcal F}^{(\pm)}\widehat H\widehat f\big)(\lambda) = \lambda\big(\widehat{\mathcal F}^{(\pm)}\widehat f\big)(\lambda), \quad 
\lambda \in \sigma_e(\widehat H)\setminus{\mathcal T}, \quad \widehat f \in \ell^2 (\mathcal{V} )$. \\
\noindent
(3) \ For $\lambda \in \sigma_e(\widehat H)\setminus{\mathcal T}$, $\widehat{\mathcal F}^{(\pm)}(\lambda)^{\ast} \in {\bf B}({\bf h}_{\lambda}\, ; \, \widehat{\mathcal B}^{\ast})$,  and $\big(\widehat H - \lambda\big)\widehat{\mathcal F}^{(\pm)}(\lambda)^{\ast}\phi = 0$ for $\phi \in {\bf h}_{\lambda}$.
\end{theorem}

The following theorem shows that the spectral representation $\widehat{\mathcal F}^{(\pm)}(\lambda)$ appears in the singularity expansion of the resolvent 
$\widehat R(\lambda \pm i0)$. Let
\begin{equation}
\widehat{\mathcal F}^{(\pm)}(\lambda)\widehat f = \left(\widehat{\mathcal F}^{(\pm)}_1(\lambda)\widehat f, \cdots, \widehat{\mathcal F}^{(\pm)}_s(\lambda)\widehat f\right).
\end{equation}


\begin{theorem}\label{ResolvSingExpand}
(Theorem 7.7 in \cite{AndIsoMor}). 
For $\widehat f \in \widehat{\mathcal B}$, we have
$$
\mathcal U_{\mathcal L_0}\widehat R(\lambda \pm i0)\widehat f \mp
\sum_{j=1}^s\frac{1}{\lambda_j(x) - \lambda \mp i0}\otimes
\widehat{\mathcal F}^{(\pm)}_j(\lambda)\widehat f \in 
\mathcal B^{\ast}_0.
$$
\end{theorem}

Note that by (\ref{Q1(lambdapmi0)define}), (\ref{Pjf(x)define}) and (\ref{DefineFpmlambda}), 
\begin{equation}
\widehat{\mathcal F}_j^{(\pm)}(\lambda)\widehat f = P_jQ_1(\lambda\pm i0)\widehat f\Big|_{\Gamma_{\lambda,j}}.
\end{equation}


\subsection{S-matrix}
The wave operators are defined by the following strong limit
\begin{equation}
\widehat W_{\pm} = {\mathop{\rm s-lim}_{t\to\pm\infty}}\, e^{it\widehat H}\widehat P_{ext}e^{-it\widehat H_0}\widehat P_{ac}(\widehat H_0),
\end{equation}
where $\widehat P_{ac}(\widehat H_0)$ is the projection onto the absolutely  continuous subspace for $\widehat H_0$. The scattering operator is then defined by
\begin{equation}
\widehat S = \big(\widehat W_+\big)^{\ast}\widehat W_-,
\end{equation}
which is unitary on $\ell^2(\mathcal V_0)$. We consider its Fourier transform $S = \widehat{\mathcal F}_0\widehat S\big(\widehat{\mathcal F}_0\big)^{\ast}$.  Letting
\begin{equation}
\widehat K_2 = \widehat H\widehat P_{ext} - \widehat P_{ext}\widehat H_0,
\label{S2DefineK2}
\end{equation}
\begin{equation}
A(\lambda) = \widehat{\mathcal F}^{(+)}(\lambda)\widehat K_2\widehat{\mathcal F}_0(\lambda)^{\ast} = \widehat{\mathcal  F}_0(\lambda)\widehat Q_1(\lambda + i0)\widehat K_2\widehat F_0(\lambda)^{\ast},
\label{DefineAlambda}
\end{equation}
we define the {\it S-matrix} by
\begin{equation}
S(\lambda) = 1 - 2\pi i A(\lambda).
\label{S2DefineSmatrix}
\end{equation}


\begin{theorem}
(Theorem 7.13 in \cite{AndIsoMor})
$S(\lambda)$ is unitary on ${\bf h}_{\lambda}$ and
\begin{equation}
\left(Sf\right)(\lambda) = S(\lambda)f(\lambda), \quad f \in {\bf H}.
\end{equation}
\end{theorem}

This S-matrix appears in the singularity expansion of solutions to the Helmholtz equation in the following way. Define the operator $A_{\pm}(\lambda) : {\bf h}_{\lambda} \to 
\mathcal B^{\ast}$ by
\begin{equation}
A_{\pm}(\lambda) = \frac{1}{2\pi i}\sum_{j=1}^s\frac{1}{\lambda_j(x) - \lambda \mp i0}\otimes P_j(x)\Big|_{x\in M_{\lambda,j}}.
\end{equation}


\begin{theorem}\label{SmatrixTheorem}
(Theorem 7.15 in \cite{AndIsoMor}) \\
\noindent
(1) $\ \{\widehat u \in \widehat{ \mathcal B}^{\ast}\, ; \, (\widehat H - \lambda)\widehat u =0\} = 
\widehat{\mathcal F}^{(-)}(\lambda)^{\ast}{\bf h}_{\lambda}$. \\
\noindent
(2) For any $\phi^{in} \in {\bf h}_{\lambda}$, there exist unique $\phi^{out} \in {\bf h}_{\lambda}$ and $\widehat u  \in \widehat{\mathcal B}^{\ast}$ satisfying
\begin{equation}
\big(\widehat H - \lambda\big)\widehat u =0,
\label{H-lambdau=0}
\end{equation}
\begin{equation}
\mathcal U_{\mathcal{L}_0} \widehat P_{ext}\widehat u + A_-(\lambda)\phi^{in} - A_+(\lambda)\phi^{out} \in 
\mathcal B_0^{\ast}.
\label{UPextuasymp}
\end{equation}
Moreover,
\begin{equation}
S(\lambda)\phi^{in} = \phi^{out}.
\end{equation}
\end{theorem}


\section{Exterior problem}


\subsection{Laplacian in the exterior domain}
\label{LaplaceExtDomain}

\label{FuncSpaceExt}

In this section, we study 
 the exterior Dirichlet  problem
\begin{equation}
\left\{
\begin{split}
& (- \widehat\Delta_{\Gamma_0} - z)\widehat u = \widehat f \quad {\rm in}  \quad \stackrel{\circ}{\mathcal V_{ext}}, \\
& \widehat u = 0 \quad {\rm on} \quad \partial \mathcal V_{ext}.
\end{split}
\right.
\label{ExtDirProb}
\end{equation}
In the exterior domain ${\mathcal V_{ext}}$, the spaces $\ell^2({\mathcal V_{ext}})$, $\ell^{2,\sigma}({\mathcal V_{ext}})$, $\widehat{\mathcal B}({\mathcal V_{ext}})$, $\widehat{\mathcal B}^{\ast}({\mathcal V_{ext}})$, $\widehat{\mathcal B}^{\ast}_0({\mathcal V_{ext}})$ are defined in the same way as in the case of periodic lattice 
$\Gamma_0 = \{\mathcal L_0, \mathcal V_0,\mathcal E_0\}$.
Let $D$ be a connected subset of ${\mathcal V_{ext}}$. 
Then for any function $\widehat u = (\widehat u_1,\cdots,\widehat u_s)$ on $D$, its normal derivative at the boundary of $D$ defined by (\ref{normalder.general}) is rewritten as 
\begin{equation}
\big(\partial_{\nu}^D \widehat u\big)_i(n) =- \frac{1}{{\rm deg}_D(i,n)}
\sum_{(j,n') \in \stackrel{\circ}D, (j,n')\sim(i,n)}\widehat u_j(n'), 
\quad (i,n) \in \partial D,
\end{equation}
where 
\begin{equation}
{\rm deg}_D(i,n) = \left\{ \begin{split}
& \sharp\{(j,n')\, ; \, (j,n') \in D, \ (j,n') \sim (i,n)\} , \quad (i,n)\in \, \stackrel{\circ}D , \\
& \sharp\{(j,n')\, ; \, (j,n') \in \, \stackrel{\circ}D, \ (j,n') \sim (i,n)\} ,\quad (i,n) \in \partial D . 
\end{split}
\right. 
\end{equation} 
 By (\ref{normalder.general}) and (\ref{GreenGeneral}),  the following Green's formula holds
\begin{equation}
\big(\widehat\Delta_{\Gamma_0}\widehat u,\widehat v\big)_{\ell^2(\stackrel{\circ}D)} -
\big(\widehat u,\widehat\Delta_{\Gamma_0}\widehat v\big)_{\ell^2(\stackrel{\circ}D)} 
= \left(\partial_{\nu}^D\widehat u,\widehat v\right)_{\ell^2(\partial D)} - \left(\widehat u,\partial_{\nu}^D\widehat v\right)_{\ell^2(\partial D)}.
\label{GreensFormula}
\end{equation}
In particular, this holds for $D = \mathcal V_{ext}$.

  We define a subspace of $\ell^2(\mathcal V_{ext})$ by 
\begin{equation}
\ell^2_0(\mathcal V_{ext}) = \{
\widehat f \in \ell^2(\mathcal V_{ext})\, ; \, \ \widehat f = 0 \ {\rm on} \ \partial\mathcal V_{ext}\},
\label{ell20Vext}
\end{equation}
and let $\widehat P_{0,ext}$ be the associated orthogonal projection
\begin{equation}
\widehat P_{0,ext}\, : \; \ell^2(\mathcal V_{ext}) \to \ell^2_0(\mathcal V_{ext}).
\label{ell20Projection}
\end{equation}
Note that $\ell^2_0(\mathcal V_{ext})$ is naturally isomorphic to $\ell^2(\stackrel{\circ}{\mathcal V_{ext}})$.
By (\ref{GreensFormula}), $ - \widehat P_{0,ext}\widehat\Delta_{\Gamma}\widehat P_{0,ext}$ is self-adjoint on $\ell^2(\mathcal V_{ext})$. Here, we extend any function $\widehat f \in \ell^2(\mathcal V_{ext})$ to be 0 outside $\mathcal V_{ext}$ so that $\widehat\Delta_{\Gamma}$ can be applied to $\widehat f$.
We take $\ell^2_0(\mathcal V_{ext})$ as the total Hilbert space and define
\begin{equation}
\widehat H_{ext} = - \widehat P_{0,ext}\widehat\Delta_{\Gamma_0}\widehat P_{0,ext}\Big|_{\ell^2_0(\mathcal V_{ext})},
\end{equation}
which is self-adjoint on $\ell^2_0(\mathcal V_{ext})$.
 Note  that
\begin{equation}
\widehat H_{ext}\widehat u = - \widehat \Delta_{\Gamma_0}\widehat u \quad 
{\rm on} \quad \stackrel{\circ}{\mathcal V_{ext}}, \quad \forall \widehat u \in \ell^2_0(\mathcal V_{ext}).
\label{HextnoP0}
\end{equation}
In fact, by definition, for $\widehat u \in \ell^2_0(\mathcal V_{ext})$
$$
\widehat H_{ext}\widehat u = - \widehat P_{0,ext}\widehat \Delta_{\Gamma_0}\widehat P_{0,ext}\widehat u = - \widehat\Delta_{\Gamma_0}\widehat u + (1 - \widehat P_{0,ext})\widehat{\Delta}_{\Gamma_0}\widehat u, 
$$
and the 2nd term of the right-hand side vanishes on $\stackrel{\circ}{\mathcal V_{ext}}$.


\begin{lemma}
(1) $\ \sigma(\widehat H_{ext}) = \sigma(\widehat H_0) = 
\sigma_e(\widehat H) = \sigma_e(\widehat H_{ext})$. \\
\noindent
(2) $\ \sigma_p(\widehat H_{ext})\cap \big(\sigma(\widehat H_{ext})\setminus
(\mathcal T_0\cup\mathcal T_1)\big) = \emptyset$.
\end{lemma}

Proof. The assertion (1) is proven by the standard method of singular sequences. The assertion (2) follows from Theorem \ref{Rellichtypetheorem} and the assumption (B-2). \qed

\medskip

 We show that Lemma \ref{RadCondUniquelemma} also holds for the exterior Dirichlet problem. The radiation condition is naturally extended to  solutions $\widehat u \in \widehat{\mathcal B}^{\ast}({\mathcal V_{ext}})$ of the equation
\begin{equation}
(- \widehat\Delta_{\Gamma_0} - \lambda)\widehat u = 0 \quad {\rm on} \quad \stackrel{\circ}{\mathcal V_{ext}}
\label{ExtHomogEq}
\end{equation}
by extending to be 0 outside $\mathcal V_{ext}$.


\begin{lemma}\label{RadCondUniqueExt}
Let $\lambda \in \sigma_e(\widehat H_{ext})\setminus\big(\mathcal T_0 \cup \mathcal T_1\big)$. If $\widehat u \in \widehat{\mathcal B}^{\ast}({\mathcal V_{ext}})$ satisfies the equation (\ref{ExtHomogEq}), the boundary condition $\widehat u = 0$ on $\partial{\mathcal V_{ext}}$ and the radiation condition, then $\widehat u$  vanishes identically  on $\mathcal V_{ext}$.
\end{lemma}

Proof. We consider the case that $\widehat u$ satisfies the outgoing radiation condition. Take $R$ large enough, and split $\widehat u$ as $\widehat u = \widehat u_e + \widehat u_0$, where
\begin{equation}
\widehat u_e = 
\left\{
\begin{split}
& \widehat u,  \quad \forall (j,n) \ {\rm s. t. } \  |n|>R, \\
& 0, \quad {\rm otherwise}.
\end{split}
\right.
\nonumber
\end{equation}
We first show
\begin{equation}
{\rm Im}\big((- \widehat\Delta_{\Gamma_0} - \lambda)\widehat u_e,\widehat u_e) =0.
\label{Inh-lambdaueue=0}
\end{equation}
In fact, since $(- \widehat\Delta_{\Gamma_0}-\lambda)\widehat u=0$  holds on $\stackrel{\circ}{\mathcal V_{ext}}$, we have by using Green's formula and the fact that $\widehat u$, $\widehat u_e$ and $\widehat u_0$ vanish on $\partial{\mathcal V}_{ext}$
\begin{equation}
\begin{split}
\big((- \widehat\Delta_{\Gamma_0}- \lambda)\widehat u_e,\widehat u_e\big) &= \big((- \widehat\Delta_{\Gamma_0}-\lambda)(\widehat u - \widehat u_0),\widehat u_e\big) = 
- \big((- \widehat\Delta_{\Gamma_0}-\lambda)\widehat u_0,\widehat u_e\big)\\
& = - \big(\widehat u_0,(- \widehat\Delta_{\Gamma_0}-\lambda)\widehat u_e\big) = -
\big(\widehat u_0,(- \widehat\Delta_{\Gamma_0}-\lambda)(\widehat u - \widehat u_0)\big)\\
&= 
\big(\widehat u_0,(\widehat H_{ext}-\lambda)\widehat u_0\big).
\end{split}
\end{equation}
The imaginary part of the right-hand side vanishes, since $\widehat H_{ext}$ is self-adjoint, and $\widehat u_0 \in D(\widehat H_{ext})$.

We now define $\widehat v \in \ell^2(\mathcal V_0)$ by the 0-extension of $\widehat u_e$ on whole $\mathcal V_0$. Then, we have
\begin{equation}
(- \widehat\Delta_{\Gamma_0}-\lambda)\widehat v = \widehat f,
\label{Eqwidehatv}
\end{equation}
where $\widehat f$ is compactly supported. 
Passing to the Fourier series, $v = \mathcal U_{\mathcal L_0}\widehat v$ satisfies
\begin{equation}
(H_0(x) - \lambda)v = f \quad {\rm on} \quad {\bf T}^d,
\end{equation}
where $f(x)$ is a trigonometric polynomial. Since $v$ is outgoing, by Lemma 6.2 of \cite{AndIsoMor}, we have
$$
P_j(x)v(x) = \frac{P_j(x)f(x)}{\lambda_j(x)-\lambda - i0},
$$
and also
\begin{equation}
{\rm Im}\, \big(v,f) = \pi \|f\big|_{M_{\lambda}}\|^2_{L^2(M_{\lambda})},
\label{tracef=0onMlambda}
\end{equation}
which vanishes by virtue of (\ref{Inh-lambdaueue=0}). Take $x^{(0)}\in M_{\lambda}$ and $\chi(x) \in C^{\infty}({\bf T}^d)$ such that $\chi(x^{(0)})=1$ and $\chi(x)=0$ outside a small neighborhood of $x^{(0)}$. We multiply the equation $(H_0(x)-\lambda)v(x) = f(x)$ by the cofactor matrix of $H_0(x)-\lambda$, and also $\chi(x)$.
Letting $w(x) = \chi(x)v(x)$, $g(x) = \chi(x)^{co}(H_0(x)-\lambda)f(x)$, we have
$$
p(x,\lambda)w(x) = g(x), \quad p(x,\lambda) = \det (H_0(x)- \lambda).
$$
Since $p(x,\lambda)$ is simple characteristic on $M_{\lambda}$, we can make a change of variables $x \to y$ taking $y_1 = p(x,\lambda)$. We write $w(x(y))$, $g(x(y))$ as $w(y)$, $g(y)$ for the sake of simplicity. Since $w(y)$ is outgoing, 
by Lemma 6.2 of \cite{AndIsoMor}, it is written as
$$
w(y) = \frac{g(y)}{y_1 - i0}.
$$
Passing to the Fourier transform, we then have 
\begin{equation}
\big\|\widetilde w(\xi_1,\cdot) - i\theta(-\xi_1)\int_{-\infty}^{\infty}
\widetilde g(\eta_1,\cdot)d\eta_1\big\|_{L^2({\bf R}^{d-1})}\to 0, \quad
{\rm as } \quad |\xi_1| \to \infty,
\label{S3wxi1cdot-itetato0}
\end{equation}
where $\theta$ is the Heaviside function (see the proof of \cite{AndIsoMor}, Lemma 4.5). By virtue of (\ref{tracef=0onMlambda}), 
$g\big|_{M_{\lambda}}=0$ holds. Therefore $g(0,y')=0$, hence $\int_{-\infty}^{\infty}\widetilde g(\eta_1,\eta')d\eta_1=0$. We have by (\ref{S3wxi1cdot-itetato0}), $\|\widetilde w(\xi_1,\cdot)\|_{L^2({\bf R}^{d-1})} \to 0$ as $\xi_1 \to \pm \infty$.
Therefore, $w(x)$ is both outgoing and incoming, hence $w(x) \in \mathcal B_0^{\ast}$.

We have thus seen that $\widehat u$ is a $\widehat{\mathcal B}_0^{\ast}$-solution to the  $(- \widehat\Delta_{\Gamma_0} - \lambda)\widehat u=0$, hence vanishes identically on $\mathcal V_{ext}$ by virtue of Theorem \ref{Rellichtypetheorem}. 
\qed

\medskip
Once we have proven Lemma \ref{RadCondUniqueExt}, the following Theorem \ref{LAPExteriordomain} can be derived  in  the same way as in the whole space \cite{AndIsoMor}, 
as was done in Theorem 6.3 of \cite{IsMo2} for the square lattice.
We do not repeat the details.

We put 
\begin{equation}
\widehat R_{ext}(z) = (\widehat H_{ext}-z)^{-1}.
\label{ExtResolvent}
\end{equation}
\begin{equation}
\mathcal T_e = \mathcal T_0 \cup \mathcal T_1.
\end{equation}


\begin{theorem}
\label{LAPExteriordomain}
 Take any compact set $I \subset \sigma_e(\widehat H_{ext})\setminus{\mathcal T}_e$, and $\lambda \in I$. Then for any $\widehat f, \widehat g \in \mathcal B^{\ast}$, there exists a limit 
\begin{equation}
\lim_{\epsilon\to 0}(\widehat R_{ext}(\lambda \pm i0)\widehat f,\widehat g) := 
(\widehat R_{ext}(\lambda \pm i0)\widehat f,\widehat g).
\end{equation}
Moreover, there exists a constant $C > 0$ such that
\begin{equation}
\|\widehat R_{ext}(\lambda \pm i0)\widehat f\|_{{\mathcal B}^{\ast}} \leq C\|\widehat f\|_{\mathcal B}, \quad \forall \lambda \in I.
\end{equation}
\end{theorem}


\subsection{Exterior and interior D-N maps}
Take $\lambda \in \sigma_e(\widehat H_{ext})\setminus \mathcal T_e$, and consider the solution $\widehat u^{(\pm)}_{ext} \in  \widehat{\mathcal B}^{\ast}$ of the following equation
\begin{equation}
\left\{
\begin{split}
(- \widehat\Delta_{\Gamma_0} - \lambda)\widehat u^{(\pm)}_{ext} = 0 \quad 
{\rm in} \quad \stackrel{\circ}{\mathcal V_{ext}}, \\
\widehat u^{(\pm)}_{ext} = \widehat f \quad {\rm on} \quad \partial\mathcal V_{ext},
\end{split}
\right.
\label{ExtDirichletProb}
\end{equation}
satisfying the radiation condition (outgoing for $u^{(+)}_{ext}$ and incoming for $u^{(-)}_{ext}$).


\begin{lemma}
For any $\lambda  \in \sigma_e(\widehat H_{ext})\setminus \mathcal T_e$, there exists a unique solution $\widehat u^{(\pm)}_{ext}$ of the exterior Dirichlet problem (\ref{ExtDirichletProb}) satisfying the radiation condition.
\end{lemma}

Proof. The uniqueness follows from Lemma \ref{RadCondUniqueExt}. To prove the existence, we extend $\widehat f$ to be 0 outside $\partial\mathcal V_{ext}$ and put
\begin{equation}
\widehat u^{(\pm)}_{ext} = \widehat f - \widehat P_0\widehat R_{ext}(\lambda \pm i0)(- \widehat\Delta_{\Gamma_0}-\lambda)\widehat f,
\label{defineupmext}
\end{equation}
where $\widehat P_0 = \widehat P_{0,ext}$. 
Then $\widehat u^{(\pm)}_{ext}=\widehat f$ on $\partial\mathcal V_{ext}$. Letting $\widehat w^{(\pm)} = \widehat P_0\widehat R_{ext}(\lambda \pm i0)(- \widehat\Delta_{\Gamma_0}-\lambda)\widehat f$, we have
$$
- \widehat\Delta_{\Gamma_0}\widehat P_0\widehat w^{(\pm)}= -  \widehat P_0\widehat\Delta_{\Gamma_0}\widehat P_0\widehat w^{(\pm)}  = 
\widehat H_{ext}\widehat w^{(\pm)} \quad 
{\rm in} \quad \stackrel{\circ}{\mathcal V_{ext}}.
$$
Here, we note that
\begin{equation}
\begin{split}
\widehat H_{ext}\widehat w^{(\pm)} &= \widehat H_{ext}\widehat R_{ext}(\lambda \pm i0)(- \widehat\Delta_{\Gamma_0}-\lambda)\widehat f\\
&= (- \widehat\Delta_{\Gamma_0}-\lambda)\widehat f + 
\lambda \widehat R_{ext}(\lambda \pm i0)(- \widehat\Delta_{\Gamma_0}-\lambda)\widehat f.
\end{split}
\nonumber
\end{equation}
Hence
\begin{equation}
\begin{split}
(- \widehat\Delta_{\Gamma_0}-\lambda)\widehat u^{(\pm)}_{ext} & = (- \widehat\Delta_{\Gamma_0}-\lambda)\widehat f  - 
(\widehat H_{ext}-\lambda)\widehat w^{(\pm)} \\
& = (1 - \widehat P_0) (- \widehat\Delta_{\Gamma_0}-\lambda)\widehat f = 0, \quad {\rm in} \quad 
\stackrel{\circ}{\mathcal V_{ext}},
\end{split}
\end{equation}
which proves the lemma. \qed

\medskip
We define the {\it exterior D-N map} $\Lambda^{(\pm)}_{ext}(\lambda)$ by
\begin{equation}
\Lambda^{(\pm)}_{ext}(\lambda)\widehat f = - \partial_{\nu}^{\mathcal V_{ext}}\widehat u^{(\pm)}_{ext}\Big|_{\partial{\mathcal V}_{ext}}.
\label{ExtDNmap}
\end{equation}

\medskip
By the assumption (B-1-3), $\mathcal V_{int}$ is a connected subgraph of $\mathcal V$, hence has the Laplacian, which is denoted by $\widehat{\Delta}_{int}$.
 We define a subspace of $\ell^2(\mathcal V_{int})$ by
\begin{equation}
\ell^2_0({\mathcal V}_{int}) =\{
 \widehat f \in \ell^2(\mathcal V_{int})\, ; \, \widehat f = 0, \ {\rm on} \ 
 \partial{\mathcal V}_{int}\},
\end{equation}
and let $\widehat P_{0,int}$ be the associated orthogonal projection
\begin{equation}
\widehat P_{0,int} : \ell^2(\mathcal V_{int}) \to \ell^2_0(\mathcal V_{int}).
\end{equation}
We define the interior Schr\"{o}dinger operator 
\begin{equation}
\widehat H_{int} = - \widehat P_{0,int}\widehat\Delta_{int}\widehat P_{0,int} + \widehat V,
\label{intDirichletop}
\end{equation}
on $\mathcal V_{int}$ with Dirichlet boundary condition on $\partial\mathcal V_{int}$.  Note that by the assumption (B-4), we have
\begin{equation}
\widehat V = \widehat P_{0,int}\widehat V\widehat P_{0,int} \quad 
{\rm on} \quad \mathcal V_{int}.
\label{AssumpV}
\end{equation}
$\widehat H_{int}$  is a finite dimensional operator, hence has a finite discrete spectrum. Then the {\it interior D-N map} $\Lambda_{int}(\lambda)$ is defined by
\begin{equation}
\Lambda_{int}(\lambda)\widehat f =  \partial_{\nu}^{\mathcal V_{int}}\widehat u_{int}\Big|_{\partial{\mathcal V}_{int}},
\label{IntDNmap}
\end{equation}
where $\lambda \not\in \sigma(\widehat H_{int})$, and $\widehat u_{int}$ is a unique solution to the equation
\begin{equation}
\left\{
\begin{split}
(- \widehat\Delta_{int} + \widehat V -\lambda)\widehat u_{int} = 0 \quad {\rm in} \quad \stackrel{\circ}{\mathcal V_{int}} ,\\
\widehat u_{int} = \widehat f \quad {\rm on} \quad \partial{\mathcal V}_{int}.
\end{split}
\right.
\label{IntBVP}
\end{equation}
Note that the uniqueness of $\widehat u_{int}$ follows from $\lambda \not\in \sigma(\widehat H_{int})$ and the existence is shown by putting
$$
\widehat u_{int} = \widehat f - \widehat P_{0,int}\widehat R_{int}(\lambda)(- \widehat\Delta_{int}-\lambda)\widehat f,
$$
where $\widehat f$ is extended to be 0 outside $\partial\mathcal V_{int}$, and 
$\widehat R_{int}(\lambda) = (\widehat H_{int} - \lambda)^{-1}$.

As in Lemma \ref{Lemma2.1}, we put
$$
\Sigma = \partial \mathcal V_{int} = \partial \mathcal V_{ext},
$$
 and define an operator $\widehat S_{\Sigma} \in {\bf B}(\ell^2(\Sigma);\ell^2(\mathcal V))$ by
\begin{equation}
\big(\widehat S_{\Sigma}\widehat f\big)(a) = 
 \frac{1}{{\rm deg}_{\mathcal V} (a)}\sum_{b \sim a, b \in \Sigma}\widehat f(b),
\end{equation}
where 
$$
{\rm deg}_{\mathcal V}(a) = {\rm deg}\,(a) = \sharp\{c \in \mathcal V\, ; \, c \sim a\}
$$
is the degree on $\mathcal V$. For a subset $A$ in $\mathcal V$, let $\chi_A$ be the characteristic function of $A$. 
In the following, we use $\chi_{\Sigma}$ to mean both of the operator of restriction 
$$
\chi_{\Sigma} : \ell^{2}_{loc}(\mathcal V) \ni \widehat f \to \widehat f\Big|_{\Sigma},
$$
 $\ell^{2}_{loc}(\mathcal V)$ being the set of locally bounded sequences, 
and the operator of extension
$$
\chi_{\Sigma} : \ell^{2}(\Sigma) \ni \widehat f \to \left\{
\begin{split}
& \widehat f, \quad {\rm on} \quad \Sigma, \\
& 0, \quad {\rm otherwise},
\end{split}
\right.
$$
without fear of confusion.  Then, we have for $\widehat f \in \ell^2(\Sigma)$
\begin{equation}
\widehat{\Delta}_{\Gamma}\chi_{\Sigma}\widehat f = 
\widehat S_{\Sigma}\widehat f.
\end{equation}
We also introduce multiplication operators by
\begin{equation}
\big(\mathcal M_{int}\widehat f\big)(a) = \frac{{\rm deg}_{\mathcal V_{int}}(a)}{{\rm deg}_{\mathcal V} (a)}\widehat f(a),
\label{S3DefineMint}
\end{equation}
\begin{equation}
\big(\mathcal M_{ext}\widehat f\big)(a) = \frac{{\rm deg}_{\mathcal V_{ext}}(a)}{{\rm deg}_{\mathcal V}(a)}\widehat f(a).
\end{equation}


\begin{lemma}\label{Lemmauintuextchif}
Let $\widehat u^{(\pm)}_{ext}$ and $\widehat u_{int}$ be the solutions of (\ref{ExtDirichletProb}) and (\ref{IntBVP}), and put
\begin{equation}
\widehat u^{(\pm)} = \chi_{\stackrel{\circ}{\mathcal V_{int}}}\widehat u_{int} + \chi_{\stackrel{\circ}{\mathcal V_{ext}}}\widehat u^{(\pm)}_{ext} + \chi_{\Sigma}\widehat f.
\label{uintuextchisigmaf}
\end{equation}
Then we have
\begin{equation}
\widehat u^{(\pm)} = \widehat R(\lambda \pm i0)\chi_{\Sigma}B^{(\pm)}_{\Sigma}(\lambda)\widehat f,
\label{u(pm)=R(lambdapmi0)chiB}
\end{equation}
where
\begin{equation}
B^{(\pm)}_{\Sigma}(\lambda) = \mathcal M_{int}\Lambda_{int}(\lambda) - 
\mathcal M_{ext}\Lambda^{(\pm)}_{ext}(\lambda) - \widehat S_{\Sigma} - \lambda\chi_{\Sigma}.
\label{DefineBpmSigmalambda}
\end{equation}
In particular, 
\begin{equation}
\widehat u^{(\pm)}_{ext} = \widehat R(\lambda \pm i0)\chi_{\Sigma}B^{(\pm)}_{\Sigma}(\lambda)\widehat f \quad {\rm on} \quad \stackrel{\circ}{\mathcal V_{ext}} ,
\end{equation}
\begin{equation}
\widehat f = \widehat R(\lambda \pm i0)\chi_{\Sigma}B^{(\pm)}_{\Sigma}(\lambda)\widehat f \quad {\rm on} \quad \Sigma.
\label{1=RchiSigmaBSIgma}
\end{equation}
\end{lemma}

Proof. Since $\widehat u^{(\pm)} = \widehat u_{int}$ in $\mathcal V_{int}$, and   $\widehat u^{(\pm)} = \widehat u_{ext}^{(\pm)}$ in $\mathcal V_{ext}$, we have\begin{equation}
( - \widehat{\Delta}_{\Gamma} + \widehat V - \lambda)\widehat u^{(\pm)} = 0 \quad {\rm in} \quad \stackrel{\circ}{\mathcal V_{int}} \cup \stackrel{\circ}{\mathcal V_{ext}}.
\label{upmsatusfyVextVint}
\end{equation}
For $a \in \Sigma$, we have
$$
\big(\widehat{\Delta}_{\Gamma}\, \chi_{\stackrel{\circ}{\mathcal V_{int}}}\widehat u_{int}\big)(a) = - \big(\mathcal M_{int} \partial_{\nu}^{\mathcal V_{int}}\widehat u_{int}\big)(a) = 
- \big(\mathcal M_{int}\Lambda_{int}(\lambda)\widehat f\big)(a),
$$
$$
\big(\widehat{\Delta}_{\Gamma}\, \chi_{\stackrel{\circ}{\mathcal V_{ext}}}\widehat u_{ext}^{(\pm)}\big)(a) = - \big(\mathcal M_{ext}\partial_{\nu}^{\mathcal V_{ext}}\widehat u_{ext}^{(\pm)}\big)(a) = 
\big(\mathcal M_{ext}\Lambda_{ext}^{(\pm)}(\lambda)\widehat f\big)(a).
$$
Therefore, we have in view of (\ref{upmsatusfyVextVint}) 
\begin{equation}
\begin{split}
& \big(- \widehat{\Delta}_{\Gamma} + \widehat V - \lambda\big)\big(\chi_{\stackrel{\circ}{\mathcal V_{int}}}\widehat u_{int} +  \chi_{\stackrel{\circ}{\mathcal V_{ext}}}\widehat u_{ext}^{(\pm)} + \chi_{\Sigma}\widehat f\big) \\
& = \chi_{\Sigma}\left(\mathcal M_{int}\Lambda_{int}(\lambda) - \mathcal M_{ext}\Lambda_{ext}^{(\pm)}(\lambda) - \big(\widehat{\Delta}_{\Gamma} + \lambda)\chi_{\Sigma}\right)\widehat f \\
&= \chi_{\Sigma}B^{(\pm)}_{\Sigma}(\lambda)\widehat f.
\end{split}
\end{equation}
Taking account of the radiation condition, we get (\ref{u(pm)=R(lambdapmi0)chiB}). \qed


\begin{lemma}\label{Lambdapmext=Lmabdampext}
For any $\lambda \in \sigma_e(\widehat H_{ext})\setminus\left(\mathcal T_e\cup\sigma(\widehat H_{int})\right)$, and $\widehat f, \widehat g \in \ell^2(\Sigma)$, we have
\begin{equation}
 \big(\Lambda_{int}(\lambda)\widehat f, \widehat g\big)_{\ell^2(\Sigma)} = 
\big(\widehat f,\Lambda_{int}(\lambda)\widehat g\big)_{\ell^2(\Sigma)},
\label{IntDNselfadjoint}
\end{equation}
\begin{equation}
 \big(\Lambda_{ext}^{(\pm)}(\lambda)\widehat f, \widehat g\big)_{\ell^2(\Sigma)} = 
\big(\widehat f,\Lambda_{ext}^{(\mp)}(\lambda)\widehat g\big)_{\ell^2(\Sigma)}.
\label{ExtDNselfadjoint}
\end{equation}
\end{lemma}

Proof. The first equality (\ref{IntDNselfadjoint}) follows from Green's formula (\ref{GreenGeneral}). To show (\ref{ExtDNselfadjoint}), let for $
z \not\in {\bf R}$
$$
\widehat u(z) = \chi_{\Sigma}\widehat f - \widehat R_{ext}(z)\big(\chi_{\stackrel{\circ}{\mathcal V_{ext}}}(-\widehat{\Delta}_{\Gamma}-z)\chi_{\Sigma}\widehat f\big).
$$
 Then, $\widehat u(z)$ is the $\ell^2$-solution to the exterior Dirichlet problem:
\begin{equation}
\left\{
\begin{split}
(- \widehat\Delta_{\Gamma_0}-z)\widehat u = 0 \quad {\rm in} \quad \stackrel{\circ}{\mathcal V_{ext}}, \\
 \widehat u = \widehat f \quad {\rm on} \quad \partial\mathcal V_{ext},
\end{split}
\right.
\nonumber
\end{equation}
and $\widehat u(\lambda + i0)$ ($\widehat u(\lambda -i0)$) satisfies the outgoing (incoming) radiation condition. Similarly, we put
$$
\widehat v(z) = \chi_{\Sigma}\widehat g - \widehat R_{ext}(z)\big(\chi_{\stackrel{\circ}{\mathcal V_{ext}}}(-\widehat{\Delta}_{\Gamma}-z)\chi_{\Sigma}\widehat g\big).
$$
Since $\widehat u(z)$, $\widehat v(z) \in \ell^2$, we have by Green's formula
\begin{equation}
\begin{split}
& \big(\widehat\Delta_{\Gamma_0}\widehat u(\lambda + i\epsilon),\widehat v(\lambda -i\epsilon)\big)_{\ell^2(\stackrel{\circ}{\mathcal V_{ext}})} -\big (\widehat u(\lambda + i\epsilon),\widehat\Delta_{\Gamma_0}\widehat v(\lambda -i\epsilon)\big)_{\ell^2(\stackrel{\circ}{\mathcal V_{ext}})} \\
& = - \big(\partial_{\nu}^{\mathcal V_{ext}}\widehat u(\lambda + i\epsilon),\widehat v(\lambda - i\epsilon)\big)_{\ell^2(\Sigma)} + \big(\widehat u(\lambda + i\epsilon),\partial_{\nu}^{\mathcal V_{ext}}\widehat v(\lambda - i\epsilon)\big)_{\ell^2(\Sigma)}.
\end{split}
\nonumber
\end{equation}
Due to the equation $\widehat \Delta_{\Gamma_0}\widehat u = - z\widehat u$, the left-hand side vanishes.
Letting $\epsilon \to 0$, we get (\ref{ExtDNselfadjoint}). \qed


\section{Scattering amplitude and D-N maps}


\subsection{Imbedding of $\ell^2(\Sigma)$ into ${\bf h}_{\lambda}$}

Let us derive the resolvent equation for $\widehat H_{ext}$.


\begin{lemma}\label{ResoventEquation}
(1) For $\widehat f \in \widehat{\mathcal B}(\mathcal V_0)$,
\begin{equation}
\begin{split}
\widehat R_{ext}(\lambda \pm i0)\chi_{\stackrel{\circ}{\mathcal V}_{ext}}\widehat f  = \left(1
 - \widehat R(\lambda \pm i0)\chi_{\Sigma}B^{(\pm)}_{\Sigma}(\lambda)\chi_{\Sigma}\right)\widehat R_0(\lambda \pm i0)\widehat f,
\end{split}
\nonumber
\end{equation}
in $\stackrel{\circ}{\mathcal V}_{ext}$. \\
\noindent
(2) For $\widehat g \in \widehat{\mathcal B}(\stackrel{\circ}{\mathcal V_{ext}})$,
\begin{equation}
\begin{split}
\widehat R_{ext}(\lambda \pm i0)\widehat g  = \widehat R_0(\lambda \pm i0)\left(1  -\chi_{\Sigma}\big(B^{(\mp)}_{\Sigma}(\lambda)\big)^{\ast}\chi_{\Sigma}\widehat R(\lambda \pm i0)\right)\chi_{\stackrel{\circ}{\mathcal V_{ext}}}\widehat g,
 \end{split}
\nonumber
\end{equation}
in $\stackrel{\circ}{\mathcal V_{ext}}$.
\end{lemma}

Proof. Let $\widehat v_0 = \widehat R(\lambda \pm i0)\chi_{\Sigma}B^{(\pm)}_{\Sigma}(\lambda)\chi_{\Sigma}\widehat R_0(\lambda \pm i0)\widehat f$. 
We replace $\widehat f$ in (\ref{uintuextchisigmaf}) by $\chi_{\Sigma}\widehat R_0(\lambda \pm i0)\widehat f$. Then, by (\ref{u(pm)=R(lambdapmi0)chiB}), we have
\begin{equation}
\begin{split}
& \chi_{{\stackrel{\circ}{\mathcal V}_{int}}}\widehat u_{int} + 
\chi_{{\stackrel{\circ}{\mathcal V}_{ext}}}\widehat u_{ext}^{(\pm)} + 
\chi_{\Sigma}\widehat R_0(\lambda \pm i0)\widehat f \\
& = \widehat R(\lambda \pm i0)\chi_{\Sigma}B^{(\pm)}_{\Sigma}(\lambda)\chi_{\Sigma}\widehat R_0(\lambda \pm i0)\widehat f = \widehat v_0.
\end{split}
\end{equation}
This implies $\widehat v_0 = \widehat u^{(\pm)}_{ext}$ in $\stackrel{\circ}{\mathcal V_{ext}}$, hence
\begin{equation}
\left\{
\begin{split}
& 
(- \widehat\Delta_{\Gamma_0}-\lambda)\widehat v_0 =0 \quad {\rm in} \quad {\stackrel{\circ}{\mathcal V_{ext}}},
\\
& \widehat v_0 = \widehat R_0(\lambda \pm i0)\widehat f \quad{\rm on} \quad \Sigma = \partial{\mathcal V}_{ext}.
\end{split}
\right.
\end{equation}
Let $\widehat w = \widehat v_0 - \widehat R_0(\lambda \pm i0)\widehat f$. Then
$$
\left\{
\begin{split}
& ( - \widehat{\Delta}_{\Gamma_0} - \lambda)\widehat w = - \widehat f \quad {\rm in} \quad \stackrel{\circ}{\mathcal V_{ext}}, \\
& \widehat w = 0 \quad {\rm on} \quad \partial{\mathcal V_{ext}}.
\end{split}
\right.
$$
Taking account of the radiation condition, we then have $\widehat w = - \widehat R_{ext}(\lambda \pm i0) \chi_{\stackrel{\circ}{\mathcal{V}_{ext}}} \widehat f =  \widehat v_0 - \widehat R_0(\lambda \pm i0)\widehat f$, which implies (1).
Taking the adjoint, we obtain 
(2).  \qed

\medskip
We introduce a spectral representation for $\widehat H_{ext}$ by
\begin{equation}
\widehat{\mathcal F}^{(\pm)}_{ext}(\lambda) = 
\widehat{\mathcal F}_0(\lambda)\left(1 - \chi_{\Sigma}\big(B^{(\mp)}_{\Sigma}(\lambda)\big)^{\ast}\chi_{\Sigma}\widehat R(\lambda \pm i0)\right)\chi_{\stackrel{\circ}{\mathcal V_{ext}}}. 
\label{DefineFpmextlambda}
\end{equation}
Lemma \ref{ResoventEquation} (2) implies
$$
\widehat{\mathcal F}^{(\pm)}_{ext}(\lambda) = 
\widehat{\mathcal F}_0(\lambda)(\widehat H_0 - \lambda)\widehat R_{ext}(\lambda \pm i0)\chi_{\stackrel{\circ}{\mathcal V_{ext}}}.
$$
Therefore, $\widehat{\mathcal F}^{(\pm)}_{ext}(\lambda)$ does not depend on the perturbation $\mathcal V_{int}$ and $\widehat V$. 
By (\ref{DefineFpmextlambda}), we have
$$
\widehat{\mathcal F}^{(-)}_{ext}(\lambda)^{\ast} \phi = \chi_{\stackrel{\circ}{\mathcal V_{ext}}}\left(1- \widehat R(\lambda + i0)\chi_{\Sigma}B^{(+)}_{\Sigma}(\lambda)\chi_{\Sigma}\right)\widehat{\mathcal F}_0(\lambda)^{\ast}\phi.
$$
By (\ref{1=RchiSigmaBSIgma}), $\widehat R(\lambda + i0)\chi_{\Sigma}B^{(+)}_{\Sigma}(\lambda)\chi_{\Sigma}=1$ on $\Sigma$, hence it is natural to define
$$
\widehat{\mathcal F}^{(-)}_{ext}(\lambda)^{\ast} \phi = 0, \quad {\rm on} \quad \Sigma.
$$
Then, we have
\begin{equation}
\widehat{\mathcal F}^{(-)}_{ext}(\lambda)^{\ast} \phi =\left(1- \widehat R(\lambda + i0)\chi_{\Sigma}B^{(+)}_{\Sigma}(\lambda)\chi_{\Sigma}\right)\widehat{\mathcal F}_0(\lambda)^{\ast}\phi, \quad 
{\rm in} \quad \mathcal V_{ext}.
\label{S4Fextast-F0ast}
\end{equation}


\begin{lemma}\label{Lemma4.2}
For any $\phi \in {\bf h}_{\lambda}$, $\widehat{\mathcal F}^{(-)}_{ext}(\lambda)^{\ast} \phi$ satisfies the equation
\begin{equation}
\left\{
\begin{split}
& (- \widehat\Delta_{\Gamma_0} - \lambda)\widehat{\mathcal F}^{(-)}_{ext}(\lambda)^{\ast} \phi =0 \quad {\rm in} \quad{\stackrel{\circ}{\mathcal V_{ext}}}, \\
& \widehat{\mathcal F}^{(-)}_{ext}(\lambda)^{\ast} \phi =0 \quad {\rm on} \quad \Sigma,
\end{split}
\right.
\end{equation}
and $\widehat{\mathcal F}^{(-)}_{ext}(\lambda)^{\ast} \phi - {\mathcal F}_0(\lambda)^{\ast} \phi$ is outgoing.
\end{lemma}

Proof. 
By Lemma \ref{Lemmauintuextchif}, $\widehat v = \widehat R(\lambda + i0)\chi_{\Sigma}B^{(+)}_{\Sigma}(\lambda)\chi_{\Sigma}\widehat{\mathcal F}_0(\lambda)^{\ast}\phi$ satisfies
$$
(- \widehat{\Delta}_{\Gamma_0} - \lambda)\widehat v=0 \quad {\rm in} \quad 
\stackrel{\circ}{\mathcal V_{ext}}, \quad 
\widehat v\big|_{\Sigma} = \widehat{\mathcal F}_0(\lambda)^{\ast}\phi.
$$
In view of (\ref{S4Fextast-F0ast}), noting that $\widehat{\mathcal F}_0(\lambda)^{\ast}\phi$ satisfies
$$
(- \widehat{\Delta}_{\Gamma_0}-\lambda)\widehat{\mathcal F}_0(\lambda)^{\ast}\phi = 0 
\quad {\rm in} \quad \stackrel{\circ}{\mathcal V_{ext}},
$$
we obtain the lemma. \qed

\medskip
We put
\begin{equation}
\widehat I^{(\pm)}(\lambda) = \widehat{\mathcal F}^{(\pm)}(\lambda)\chi_{\Sigma}B^{(\pm)}_{\Sigma}(\lambda) : \ell^2(\Sigma) \to {\bf h}_{\lambda}.
\end{equation}
By (\ref{DrfineQ1(z)}), (\ref{DefineF0labda=F0lambdaUL0}) and (\ref{DefineFpmlambda}), we have
$$
\widehat I^{(\pm)}(\lambda) = \mathcal F_0(\lambda)\mathcal U_{\mathcal L_0}(\widehat H_0-\lambda)\widehat P_{ext}\widehat R(\lambda \pm i0)\chi_{\Sigma}B^{(\pm)}_{\Sigma}(\lambda),
$$
which yields by virtue of (\ref{Pextdefine}), (\ref{uintuextchisigmaf}) and 
(\ref{u(pm)=R(lambdapmi0)chiB})
\begin{equation}
\widehat I^{(\pm)}(\lambda)\widehat f = \mathcal F_0(\lambda)\mathcal U_{\mathcal L_0}(\widehat H_0 - \lambda)\widehat P_{ext}\widehat u^{(\pm)}.
\label{Ipmlabda=mathcalFH-lambdauext}
\end{equation} 
This formula shows that $I^{(\pm)}(\lambda)$ depends neither  on $\mathcal V_{int}$ nor on $\widehat V$, i.e. it is independent of the perturbation. 


\begin{lemma}\label{Lemma4.3}
\noindent
(1) $\ \widehat I^{(\pm)}(\lambda) : \ell^2(\Sigma) \to {\bf h}_{\lambda}$ is 1 to 1.  \\
\noindent
(2) $\ \widehat I^{(\pm)}(\lambda)^{\ast}  : {\bf h}_{\lambda} \to \ell^2(\Sigma)$ is onto.
\end{lemma}

Proof. Suppose $\widehat I^{(\pm)}(\lambda)\widehat f=0$, and let $\widehat u^{(\pm)}_{ext}$ be the solution to (\ref{ExtDirichletProb}). 
In view of (\ref{uintuextchisigmaf}) and (\ref{u(pm)=R(lambdapmi0)chiB}), we have
 in $\stackrel{\circ}{\mathcal V_{ext}}$,
$$
\widehat u^{(\pm)}_{ext} = \widehat R(\lambda \pm i0)\widehat g, \quad
\widehat g = \chi_{\Sigma}B^{(\pm)}_{\Sigma}(\lambda)\widehat f.
$$
 Theorem \ref{ResolvSingExpand} implies
 $$
\mathcal U_{\mathcal L_0}\widehat R(\lambda \pm i0)\widehat g  
\mp \sum_{j=1}^s\frac{1}{\lambda_j(x) - \lambda \mp i0}\otimes
\widehat{\mathcal F}_j^{(\pm)}(\lambda)\widehat g \in {\mathcal B}^{\ast}_0.
 $$
Since $\widehat I^{(\pm)}(\lambda)\widehat f = \widehat{\mathcal  F}^{(\pm)}(\lambda)\widehat g$, this implies $\widehat u^{(\pm)}_{ext} \in 
\widehat{\mathcal B}^{\ast}_0$.
 The Rellich type theorem (Theorem \ref{Rellichtypetheorem}) and the unique continuation property (B-2) entails $\widehat u^{(\pm)}_{ext}=0$, which yields $\widehat f=0$. This proves (1), which implies that the range of $\widehat I^{(\pm)}(\lambda)^{\ast}$ is dense.  Since $\ell^2(\Sigma)$ is finite dimensional, (2) follows. \qed


\subsection{Scattering amplitude in the exterior domain}
Similarly to the scattering amplitude (\ref{DefineAlambda}), we define the scattering amplitude in the exterior domain by
\begin{equation}
A_{ext}(\lambda) = \widehat{\mathcal F}^{(+)}(\lambda)\chi_{\Sigma}B^{(+)}_{\Sigma}(\lambda)\chi_{\Sigma}\widehat{\mathcal F}_0(\lambda)^{\ast}.
\label{Aextlambda}
\end{equation}
Using Theorem \ref{ResolvSingExpand} and (\ref{S4Fextast-F0ast}), and extending $\widehat{\mathcal F}^{(-)}_{ext}(\lambda)^{\ast}\phi - 
\widehat{\mathcal F}_0(\lambda)^{\ast}\phi$ to be 0 outside $\mathcal V_{ext}$, we have
\begin{equation}
\begin{split}
\mathcal U_{\mathcal L_0}\left(\widehat{\mathcal F}^{(-)}_{ext}(\lambda)^{\ast}\phi - 
\widehat{\mathcal F}_0(\lambda)^{\ast}\phi\right)& = 
- \mathcal U_{\mathcal L_0} \widehat R(\lambda + i0)\chi_{\Sigma}B^{(+)}_{\Sigma}(\lambda)\chi_{\Sigma}
\widehat{\mathcal F}_0(\lambda)^{\ast}\phi \\
& \simeq - \sum_{j=1}^s\frac{1}{\lambda_j(x)-\lambda - i0}\otimes A_{ext,j}(\lambda)\phi,
\end{split}
\end{equation}
where $A_{ext,j}(\lambda)\phi$ denotes the $j$-th component of $A_{ext}(\lambda)\phi$. This shows that $A_{ext}(\lambda)\phi$ depends only on $\mathcal V_{ext}$.


\subsection{Single layer and double layer potentials}
The operator
\begin{equation}
\widehat R(\lambda \pm i0)\chi_{\Sigma}B^{(\pm)}_{\Sigma}(\lambda) : \ell^2(\Sigma) \to \mathcal B^{\ast}(\mathcal V)
\end{equation}
is an analogue of the double layer potential in the continuous case. Similarly, the operator defined by
\begin{equation}
\ell^2(\Sigma) \ni \widehat f \to M^{(\pm)}_{\Sigma}(\lambda)\widehat f := \left(\widehat R(\lambda \pm i0)\chi_{\Sigma}\widehat f\right)\Big|_{\Sigma} \in \ell^2(\Sigma)
\end{equation}
is an analogue of the single layer potential.


\begin{lemma}\label{Lemma4.4}
$M^{(\pm)}_{\Sigma}(\lambda)B^{(\pm)}_{\Sigma}(\lambda) =1$ on $ \ell^2(\Sigma)$.
\end{lemma}

This follows from (\ref{1=RchiSigmaBSIgma}). In particular, 
$M^{(\pm)}_{\Sigma}(\lambda) = (B^{(\pm)}_{\Sigma}(\lambda))^{-1}$.


\subsection{S-matrix and interior D-N map}
The scattering amplitude $A(\lambda)$ in the whole space (\ref{DefineAlambda}) and the scattering amplitude $A_{ext}(\lambda)$ in the exterior domain (\ref{Aextlambda}) have the following relation.


\begin{theorem}
\label{ThSmatrixDNmapEquiv}
We have
\begin{equation}
A_{ext}(\lambda) - A(\lambda) = \widehat I^{(+)}(\lambda)\left(B^{(+)}_{\Sigma}(\lambda)\right)^{-1}
\widehat I^{(-)}(\lambda)^{\ast}.
\label{Aext-A=I+MI-}
\end{equation}
\end{theorem}

Proof. For $\phi \in {\bf h}_{\lambda}$, let
\begin{equation}
\widehat u = \widehat{\mathcal F}^{(-)}(\lambda)^{\ast} \phi - 
\widehat{\mathcal F}^{(-)}_{ext}(\lambda)^{\ast}\phi.
\end{equation}
By (\ref{Q1(lambdapmi0)define}) and (\ref{S2DefineK2}), we have
\begin{equation}
\begin{split}
\widehat Q_1(\lambda - i0) = (\widehat H_0 - \lambda)\widehat P_{ext}\widehat R(\lambda - i0) = - \widehat K_2^{\ast}\widehat R(\lambda - i0) + \widehat P_{ext}.
\end{split}
\nonumber
\end{equation}
Using (\ref{DefineFpmlambda}) and (\ref{S4Fextast-F0ast}), we  then have
\begin{equation}
\widehat u =\left( (\widehat P_{ext}-1) + \widehat R(\lambda+i0)(\chi_{\Sigma}B_{\Sigma}^{(+)}(\lambda)\chi_{\Sigma}-\widehat K_2)\right)\widehat{\mathcal F}_0(\lambda)^{\ast}\phi.
\end{equation}
Theorem \ref{ResolvSingExpand} then implies
\begin{equation}
\mathcal U_{\mathcal L_0}\widehat u \simeq \sum_{j=1}^d\frac{1}{\lambda_j(x)-\lambda - i0}
\otimes \widehat{\mathcal F}^{(+)}_j(\lambda)(\chi_{\Sigma}B^{(+)}_{\Sigma}(\lambda)\chi_{\Sigma}
- \widehat K_2)\widehat{\mathcal F}_0(\lambda)^{\ast}\phi.
\label{S4uwidahatuexpand1}
\end{equation}
By virtue of Lemma \ref{Lemma4.2}, $\widehat u$ is the outgoing solution of the equation
\begin{equation}
(- \widehat\Delta_{\Gamma_0}-\lambda)\widehat u = 0 \quad {\rm in} \quad \stackrel{\circ}{\mathcal V}_{ext}, \quad
\widehat u\big|_{\Sigma} = \widehat{\mathcal F}^{(-)}(\lambda)^{\ast}\phi.
\end{equation}
Lemma \ref{Lemmauintuextchif} yields
\begin{equation}
\widehat u = \widehat R(\lambda + i0)\chi_{\Sigma}B^{(+)}(\lambda)\widehat{\mathcal F}^{(-)}(\lambda)^{\ast}\phi.
\end{equation}
Again using Theorem \ref{ResolvSingExpand},
\begin{equation}
\mathcal U_{\mathcal L_0} \widehat u \simeq \sum_{j=1}^d\frac{1}{\lambda_j(x)-\lambda -i0}\otimes
\widehat{\mathcal F}_j^{(+)}(\lambda)\chi_{\Sigma}B_{\Sigma}^{(+)}(\lambda)\widehat{\mathcal F}^{(-)}(\lambda)^{\ast}\phi.
\label{S4uwidahatuexpand2}
\end{equation}
Comparing (\ref{S4uwidahatuexpand1}) and (\ref{S4uwidahatuexpand2}), we have
\begin{equation}
\begin{split}
& \widehat{\mathcal F}^{(+)}(\lambda)(\chi_{\Sigma}B^{(+)}_{\Sigma}(\lambda)\chi_{\Sigma}
- \widehat K_2)\widehat{\mathcal F}_0(\lambda)^{\ast}\phi \\
&=\widehat{\mathcal F}^{(+)}(\lambda)\chi_{\Sigma}B_{\Sigma}^{(+)}(\lambda)\widehat{\mathcal F}^{(-)}(\lambda)^{\ast}\phi.
\end{split}
\end{equation}
By (\ref{DefineAlambda}) and (\ref{Aextlambda}), the left-hand side is equal to 
$A_{ext}(\lambda) - A(\lambda)$.

Let us note here that $1 = B^{(-)}_{\Sigma}(\lambda)M^{(-)}_{\Sigma}(\lambda)$ by virtue of Lemma \ref{Lemma4.4}. Since $M^{(\pm)}_{\Sigma}(\lambda) = 
\chi_{\Sigma}\widehat R(\lambda \pm i0)\chi_{\Sigma}$, we have $(M^{(-)}_{\Sigma}(\lambda))^{\ast} = M^{(+)}_{\Sigma}(\lambda)$, which implies 
\begin{equation}
1 = M^{(+)}_{\Sigma}(\lambda)\big(B^{(-)}_{\Sigma}(\lambda)\big)^{\ast}.
\label{1=M+B-ast}
\end{equation}
Inserting (\ref{1=M+B-ast}) 
between  $B^{(+)}_{\Sigma}(\lambda)$ and $\widehat{\mathcal F}^{(-)}(\lambda)^{\ast}\phi$, we obtain
\begin{equation}
\begin{split}
& \widehat{\mathcal F}^{(+)}(\lambda)\chi_{\Sigma}B^{(+)}\chi_{\Sigma}\widehat{\mathcal F}^{(-)}(\lambda)^{\ast} \\
&= \widehat{\mathcal F}^{(+)}(\lambda)\chi_{\Sigma}B^{(+)}_{\Sigma}(\lambda)M^{(+)}_{\Sigma}(\lambda)\big(B^{(-)}_{\Sigma}(\lambda)\big)^{\ast}\chi_{\Sigma}\widehat{\mathcal F}^{(-)}(\lambda)^{\ast} \\
&= \widehat I^{(+)}(\lambda)M^{(+)}_{\Sigma}(\lambda)\widehat I^{(-)}(\lambda)^{\ast}.
\end{split}
\end{equation}
We have thus proven (\ref{Aext-A=I+MI-}). \qed


\subsection{The operator $\widehat J^{(\pm)}(\lambda)$}
To construct $A(\lambda)$ from $B^{(+)}_{\Sigma}(\lambda)$, we need to invert $\widehat I^{(\pm)}(\lambda)$ and its adjoint.
 To compute them, we first construct a solution $\widehat u^{(\pm)}_{ext}$ to the exterior Dirichlet problem satisfying (\ref{ExtDirichletProb}) 
and the radiation condition in the form 
$\widehat R_0(\lambda \pm i0)\widehat \psi$, where $\widehat\psi \in \ell^2(\Sigma)$. Then it is the desired solution if and only if
\begin{equation}
\widehat R_0(\lambda \pm i0)\widehat \psi = \widehat f \quad {\rm on} \quad \Sigma.
\label{S4BoundaryIntEq}
\end{equation}
Suppose $\widehat R_0(\lambda \pm i0)\widehat\phi^{(\pm)}=0$ on $\Sigma$. Then, $\widehat v^{(\pm)} = \widehat R_0(\lambda \pm i0)\widehat\phi^{(\pm)}$ is the solution to the equation (\ref{ExtDirichletProb})  with 0 boundary data. Since $\widehat v^{(\pm)}$ satisfies the radiation condition, by Lemma \ref{RadCondUniqueExt}, it vanishes identically in $\stackrel{\circ}{\mathcal V}_{ext}$ hence on all $\mathcal V_0$. It then follows that $\widehat\phi^{(\pm)} =0$. Therefore, the equation (\ref{S4BoundaryIntEq}) is uniquely solvable for any $\widehat f \in \ell^2(\Sigma)$. Let $\widehat\psi = \widehat r_{\Sigma}^{(\pm)}(\lambda)\widehat f$ be the solution. 
Then, we have
\begin{equation}
\widehat u^{(\pm)}_{ext} = \widehat R_0(\lambda \pm i0)\widehat r_{\Sigma}^{(\pm)}(\lambda)\widehat f, 
\label{S4Solformulauext}
\end{equation}
which is a potential theoretic solution to the boundary value problem (\ref{ExtDirichletProb}).

Let $\widehat g_n, n = 1,\cdots,N$, be a basis of $\ell^2(\Sigma)$ and put
\begin{equation}
\begin{split}
v^{(\pm)}_n & = \widehat I^{(\pm)}(\lambda)\widehat g_n \\
& = 
\mathcal F_0(\lambda)\mathcal U_{\mathcal L_0}(\widehat H_0 -\lambda)\widehat P_{ext}\widehat R_0(\lambda \pm i0)\widehat r_{\Sigma}^{(\pm)}(\lambda)\widehat g_n \in {\bf h}_{\lambda}.
\end{split}
\end{equation}
Let $\mathcal M_{\Sigma}^{(\pm)}$ be the linear hull of $v^{(\pm)}_1,\cdots,v^{(\pm)}_N$. Then, the mapping $\widehat g_n \to v_n^{(\pm)}$ induces a bijection 
$$
\widehat J^{(\pm)}(\lambda) : \ell^2(\Sigma) \ni \sum_{n=1}^Nc_n\widehat g_n \to \sum_{n=1}^Nc_nv^{(\pm)}_n \in \mathcal M_{\Sigma}^{(\pm)}.
$$
In view of  Theorem 
\ref{ThSmatrixDNmapEquiv}, we have the following theorem.


\begin{theorem} \label{LemmaAlambdatoBSigmalambda}
The following formula holds:
\begin{equation}
\big(B^{(+)}_{\Sigma}(\lambda)\big)^{-1} = 
\big(\widehat J^{(+)}(\lambda)\big)^{-1}\big( A_{ext}(\lambda) - A(\lambda)\big)
\big(\widehat J^{(-)}(\lambda)^{\ast}\big)^{-1}.
\end{equation}
\end{theorem}

By virtue of Theorems \ref{ThSmatrixDNmapEquiv} and \ref{LemmaAlambdatoBSigmalambda}, the S-matrix and the D-N map determine each other.

\subsection{Perturbation of S-matrices}
Suppose we are given two interior domains ${\mathcal V}_{int,1}$ and ${\mathcal V}_{int,2}$ such that $\Sigma = \partial{\mathcal V}_{int,1} = \partial{\mathcal V}_{int,2}$ and 
\begin{equation}
{\rm deg}_{{\mathcal V}_{int,1}} = {\rm deg}_{{\mathcal V}_{int,2}} \quad {\rm on} \quad \Sigma.
\label{S4degVint1=degVint2}
\end{equation}
We put the suffix $i = 1,2, $ for the operators $A(\lambda)$, $\Lambda_{int}(\lambda)$, $\mathcal M_{int}$, $\widehat S_{\Sigma}$ and $B^{(+)}_{\Sigma}(\lambda)$ associated with the domain $\mathcal V_{int,i}$. Then, by the condition (\ref{S4degVint1=degVint2}),
\begin{equation}
\mathcal M_{int,1} = \mathcal M_{int,2}, \quad
\widehat S_{\Sigma,1} =\widehat S_{\Sigma,2} \quad {\rm on} \quad \Sigma.
\label{S4i=1,2equal}
\end{equation} 
Then, we have by the resolvent equation, 
\begin{equation}
\begin{split}
& (B^{(+)}_{\Sigma,2}(\lambda))^{-1} - (B^{(+)}_{\Sigma,1}(\lambda))^{-1}\\
& = - (B^{(+)}_{\Sigma,2}(\lambda))^{-1}(\mathcal M_{int}\Lambda_{int,2}(\lambda) - \mathcal M_{int}\Lambda_{int,1}(\lambda))(B^{(+)}_{\Sigma,1}(\lambda))^{-1},
\end{split}
\end{equation}
where $\mathcal M_{int} = \mathcal M_{int,1} = \mathcal M_{int,2}$.
 Theorem \ref{ThSmatrixDNmapEquiv} then implies the following lemma.


\begin{lemma} \label{LemmaS4A2-A1=int2-int1}
The following formula holds:
\begin{equation}
\begin{split}
&A_2(\lambda) - A_1(\lambda) \\
& = \widehat I^{(+)}(\lambda)\big(B^{(+)}_{\Sigma,2}(\lambda)\big)^{-1}
\mathcal M_{int}\big(\Lambda_{int,2}(\lambda) - \Lambda_{int,1}(\lambda)\big)\big(B^{(+)}_{\Sigma,1}(\lambda)\big)^{-1}\widehat I^{(-)}(\lambda)^{\ast}.
\end{split}
\nonumber
\end{equation}
\end{lemma}

Therefore, if we can find a data $\widehat f$ on $\Sigma$ such that $\big(\Lambda_{int,2}(\lambda) - \Lambda_{int,1}(\lambda)\big)\widehat f \neq 0$, we can distinguish between $\mathcal V_{int,1}$ and $\mathcal V_{int,2}$ by the scattering experiment.


\section{Asymptotic behavior of wave functions in the lattice space}
We have defined the S-matrix by using the singularity expansion of the solution to the Schr{\"o}dinger equation. However, in some energy region, we can derive it from the spatial asymptotics at infinity of the lattice space. We prove it here because of its physical importance, although it is not used in the later sections.

Recall that by (\ref{DefineAlambda}) and (\ref{S2DefineSmatrix}),
\begin{equation}
\phi^{out} = S(\lambda)\phi^{in},
\end{equation}
\begin{equation}
\phi^{out} = 
\phi^{in} - 2 \pi i \mathcal F_0(\lambda)\mathcal U_{\mathcal L_0}\widehat Q_1(\lambda + i0)\widehat K_2\widehat{\mathcal F}_0(\lambda)^{\ast}\phi^{in}.
\label{S5_eq_inout}
\end{equation}
We compute $\widehat Q_1(z)$ as follows
\begin{equation}
\begin{split}
\widehat Q_1(z) = (\widehat H_0-z)\widehat P_{ext}\widehat R(z) = \widehat P_{ext} + \widehat K_1\widehat R(z),
\end{split}
\end{equation}
\begin{equation}
\widehat K_1 = \widehat H_0\widehat P_{ext} - \widehat P_{ext}\widehat H.
\end{equation}
Since $\widehat{K}_1 $ and $\widehat{K}_2 $ are finite dimensional operators, (\ref{S5_eq_inout}) implies that $\phi^{out} \in C^{\infty}(M_{\lambda})$ if $\phi^{in} \in C^{\infty}(M_{\lambda})$. 
By Theorem \ref{SmatrixTheorem}, there exists a unique $\widehat u \in \widehat{\mathcal B}^{\ast}$ satisfying (\ref{H-lambdau=0}), (\ref{UPextuasymp}). We observe the behavior of $\widehat u$ modulo  $\widehat{\mathcal B}^{\ast}_0$. In view of Theorem \ref{SmatrixTheorem}, we have only to study 
\begin{equation}
\frac{1}{\lambda_j(x) - \lambda \mp i0}\otimes\big(P_j(x)\phi^{(\pm)}(x)\big)\Big|_{x \in M_{\lambda,j}},
\end{equation}
where $\phi^{(+)} = \phi^{out}, \phi^{(-)} = \phi^{in}$.

\medskip
Here we impose a new assumption which is used only in this section :

\medskip
\noindent
{\bf(C)} {\it There exists $\lambda \in \sigma_e(\widehat H)\setminus{\mathcal T}$ such that 
for any $ 1 \leq j \leq s$, $M_{\lambda,j}$ is strictly convex.}

\medskip
Note that in the Assumption (C), we allow the case in which  $M_{\lambda,j} = \emptyset$ for some $j$. 
Let us  compute the asymptotic expansion of the integral
\begin{equation}
I(k) = \int_{{\bf T}^d}\frac{e^{ix\cdot k}f(x)}{\lambda(x) - \lambda \mp i0}dx,
\end{equation}
assuming that $M_{\lambda} = \{x \in {\bf T}^d\,; \, \lambda(x)=\lambda\}$ is strictly convex. It is well-known that
\begin{equation}
I(k) = \pm i \pi \int_{M_{\lambda}}\frac{e^{ix\cdot k}f(x)}{|\nabla \lambda(x)|}dM_{\lambda} + {\rm p.v}\int_{{\bf T}^d}\frac{e^{ix\cdot k}f(x)}{\lambda (x) - \lambda}dx.
\end{equation}

Let $N(x)$ be the outward unit normal field at $x \in M_{\lambda}$. Since $M_{\lambda}$ is strictly convex, for any $\omega \in S^{d-1}$, there exists a unique $x^{(\pm)}(\lambda,\omega) \in M_{\lambda}$ such that
$$
N(x^{(\pm)}(\lambda,\omega)) = \pm \omega.
$$
Letting $\omega_k = k/|k|$, we have by the stationary phase method
\begin{equation}
I(k) =  \sum_{\pm}|k|^{-(d-1)/2}e^{ik\cdot x^{(\pm)}(\lambda,\omega_k)}a^{(\pm)}(\lambda,\omega_k)f(x^{(\pm)}(\lambda,\omega_k)) + 
O(|k|^{-(d+1)/2}),
\nonumber
\end{equation}
\begin{equation}
a^{(\pm)}(\lambda, \omega_k) = C_{\pm}\frac{K(x^{(\pm)}(\lambda,\omega_k))^{-1/2}}{|\nabla \lambda(x^{(\pm)}(\lambda,\omega_k))|}
\nonumber
\end{equation}
as $|k| \to \infty$,
where $C_{\pm} = \pm i (2\pi )^{(d-1)/2} e^{\mp i (d-1) \pi /4} $ and $K(x)$ is the Gaussian curvature of $M_{\lambda}$ at $x \in M_{\lambda}$ (see Lemma 4.4 of \cite{IsMo2}). 
We replace $\lambda(x)$ by $\lambda_j(x)$, and define $x_j^{(\pm)}(\lambda,\omega)$, $a_j^{(\pm)}(\lambda,\omega)$ and $K_j(x)$ in the same way as above. We can thus reformulate (\ref{UPextuasymp}) into the the following theorem.


\begin{theorem}
Assume (C). If $\widehat u \in \widehat B^{\ast}$ satisfies $(\widehat H - \lambda)\widehat u = 0$ and $\phi^{in}, \phi^{out} \in C^{\infty}(M_{\lambda})$, we have 
the following asymptotic expansion in the lattice space
\begin{equation}
\begin{split}
\widehat u = & - \sum_{j=1}^s |k|^{-(d-1)/2}e^{ik\cdot x^{(-)}_j(\lambda,\omega_k)}a_j^{(-)}(\lambda,\omega_k)\phi_j^{in}(x^{(-)}_j(\lambda,\omega_k)) \\
& + \sum_{j=1}^s |k|^{-(d-1)/2}e^{ik\cdot x^{(+)}_j(\lambda,\omega_k)}a_j^{(+)}(\lambda,\omega_k)\phi_j^{out}(x^{(+)}_j(\lambda,\omega_k)) \\
& + O(|k|^{-(d+1)/2}).
\end{split}
\label{LatticeSpaceAsymp}
\end{equation}
\end{theorem}

The standard way of defining the S-matrix is to use the asymptotic expansion of the form (\ref{LatticeSpaceAsymp}), i.e. the operator
\begin{equation}
\phi^{in}_j(x^{(-)}_j(\lambda,\omega_k))  \to \phi^{out}_j(x^{(+)}_j(\lambda,\omega_k)), \quad \omega_k \in S^{d-1}
\label{GemetricSmatrix}
 \end{equation} 
 is the S-matrix based on the far field pattern of wave functions. 
This coincides with our definition of S-matrix (\ref{S2DefineSmatrix}) up to the parametrization of $M_{\lambda}$ (i.e. $\omega$ or $x^{(\pm)}(\lambda,\omega)$). 
 We omit the proof of this fact. 

 \medskip
 Let us check the assumption (C) in our case. For all of the examples given in \cite{AndIsoMor}, $p(x,\lambda) = \det (H_0(x) - \lambda)$ is written as
 $$
 p(x,\lambda) = f(a_d(x),\lambda), \quad {\rm or} \quad f(b_d(x),\lambda),
 $$
 where $f(z,\lambda)$ is a polynomial of two variables $z, \lambda$, and 
$$
a_d(x) = \sum_{j=1}^d\cos x_j, \quad b_d(x) = \sum_{j=1}^d\cos x_j + \sum_{1\leq j < k\leq d}\cos(x_j - x_k).
$$
We factorize $p(x,\lambda)$ as: 
 $$
 p(x,\lambda) = c_0(\lambda)
\left\{
\begin{split}
&\prod_{i=1}^m(a_d(x) - c_i(\lambda))^{\mu_i}, \quad (A) \\
& \prod_{i=1}^m(b_d(x) - c_i(\lambda))^{\mu_i}, \quad (B)
\end{split}
\right.
 $$
 where $c_i(\lambda) \neq c_j(\lambda)$ if $i \neq j$.
 
\begin{figure}[htbp]
 \begin{minipage}{0.49\hsize}
  \begin{center}
\includegraphics[width=50mm, bb=0 0 217 222]{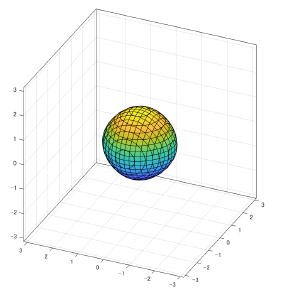}
  \end{center}
  \caption{$a_3 (x)= c \in I'_3$} 
  \label{S5FermiSurfacef(x)square}
 \end{minipage}
 \begin{minipage}{0.49\hsize}
  \begin{center}
\includegraphics[width=50mm, bb=0 0 228 235]{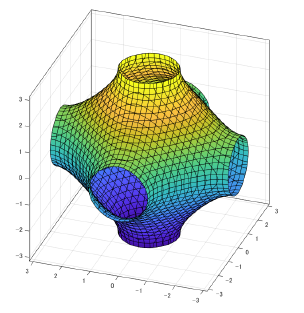}
  \end{center}
  \caption{$a_3 (x)= c \in (-3,3) \setminus I'_3$}
  \label{S5FermiSurfacef(x2)square}
 \end{minipage}
\end{figure}

 \medskip
 \noindent
 {\it The case (A)}. In this case, $M_{\lambda,j} = \{x \in {\bf T}^d\, ;\, \, a_d(x)=c_i(\lambda)\}$ for some $i$. Let 
 \begin{equation}
 I'_d = \left\{
 \begin{split}
 & (-2, 0) \cup (0,2), \quad {\rm for} \quad d=2, \\
 & (-d,-d+1)\cup(d-1,d), \quad {\rm for} \quad d\geq3.
 \end{split}
 \right.
 \end{equation}
 By Lemma 2.1 of \cite{AndIsoMor}, we have
 \begin{equation}
 a_d({\bf T}^d) = [-d,d],
 \end{equation}
 and by Lemma 4.3 of \cite{IsMo2}, for $c \in I_d'$, the surface $\{x \in {\bf T}^d\, ; \, a_d(x) = c\}$ is strictly convex (see Figures \ref{S5FermiSurfacef(x)square} and \ref{S5FermiSurfacef(x2)square}). 
 In view of the formulas given in  \S 3 of \cite{AndIsoMor}, we thus have:
 
 \begin{itemize}
 \item For the square lattice, 
$$
p(x,\lambda) = - \frac{1}{d}\left(a_d(x) + \lambda d\right).
$$
  Hence, $M_{\lambda}$ is strictly convex for 
  $\lambda \in (-1,0)\cup (0,1)$ when $d=2$, and for $\lambda \in (-1, -1 + 1/d)\cup (1- 1/d,1)$ when $d \geq 3$.
 
\item For the subdivision of $d$-dim. square lattice, 
$$
p(x,\lambda) = - \frac{(-\lambda)^{d-1}}{2d}\left(a_d(x) - 2d \lambda^2 + d\right).
$$
Therefore, when $d=2$, $M_{\lambda}$ is strictly convex, if 
$$
\lambda \in \big(-1, 1\big)\setminus\Big\{ \pm \sqrt{\frac{1}{2}}, 0\Big\},
$$
and when $d \geq 3$,  if
$$
\lambda \in \Big(-1,-\sqrt{1- \frac{1}{2d}}\Big)\cup\Big(-\sqrt{\frac{1}{2d}},0\Big)
\cup\Big(0,\sqrt{\frac{1}{2d}}\Big)\cup\Big(\sqrt{1- \frac{1}{2d}},1\Big).
$$

 \item For the ladder of $d$-dim. square lattice, 
$$
p(x,\lambda) = \Big(\frac{2}{2d+1}\Big)^2\Big(a_d(x) + \frac{(2d+1)\lambda +1}{2}\Big)\Big(a_d(x) + \frac{(2d+1)\lambda -1}{2}\Big).
$$
Therefore, when $d=2$, $M_{\lambda}$ is strictly convex if
$$
\lambda \in \big(-1,1\big)\setminus\big\{\pm \frac{1}{5}\big\},
$$
and when $d \geq 3$ if
 $$
 \lambda \in \Big(-1, - \frac{2d-3}{2d+1}\Big)\cup\Big(\frac{2d-3}{2d+1},1\Big)\setminus\big\{\pm \frac{2d-1}{2d+1}\big\}.
 $$
 \end{itemize}
 
 \medskip
 \noindent
 {\it The case (B)}. In this case, $M_{\lambda,j} = \{x \in {\bf T}^d\, ;\, \, b_d(x)=c_i(\lambda)\}$ for some $i$.  
 For the sake of simplicity, we consider only the case $d=2$. 
 By Lemma 2.2 of \cite{AndIsoMor}, we have $b_2 ({\bf T}^2) = [-3/2,3]$. 
\begin{figure}[htbp]
 \begin{minipage}{0.49\hsize}
  \begin{center}
 \includegraphics[width=65mm, bb=0 0 269 207]{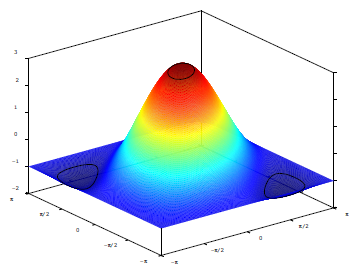}
  \end{center}
  \caption{$x_3 = b_2(x_1, x_2)$}
  \label{S5FermiSurfacef(x)}
 \end{minipage}
 \begin{minipage}{0.49\hsize}
  \begin{center}
   \includegraphics[width=65mm, bb=0 0 219 212]{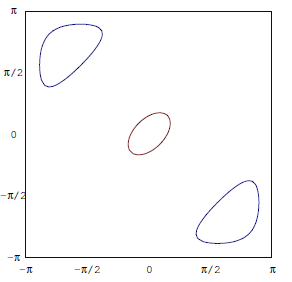}
  \end{center}
  \caption{$b_2(x_1, x_2) = \kappa$}
  \label{S5FermiSurfacef(x2)}
 \end{minipage}
\end{figure}
%
%
%
%
%
 
 Put $C_{\kappa} = \{ x \in {\bf T}^2\, ; \, b_2 (x) = \kappa\}$. 
 Taking note of the inequality
 $\cos x_1 + \cos x_2 + \cos(x_1-x_2) \leq 3$, and noting that the equality occurs only when $x_1 = x_2 = 0$, we have
 $C_3 = \{(0,0)\}$. Therefore, if $3 - \epsilon < \kappa < 3$,  $\epsilon > 0$ being chosen sufficiently small, $C_{\kappa}$ is a regular closed curve enclosing 
 $(0,0)$. By the Taylor expansion, 
 $$
b_2 (x_1,x_2) = \frac{3}{2} - \frac{1}{2}(x_1^2 - x_1x_2 + x_2^2) + O(|x|^3), 
$$
 which does not vanish on $C_{\kappa}$ for $3 -\epsilon < \kappa < 3$. Therefore, $C_{\kappa}$ is strictly convex if $3 - \epsilon < \kappa < 3$.
 
 We also have $\cos x_1 + \cos x_2 + \cos(x_1-x_2) \geq - 3/2$, and the equality occurs only when 
 $(x_1,x_2) = (4\pi/3,2\pi/3), (2\pi/3,4\pi/3)$. Therefore, 
 $$
 C_{-3/2} =  \{(4\pi/3,2\pi/3), (2\pi/3,4\pi/3)\}.
 $$
 Letting $\xi = (x_1 -4\pi/3, x_2-2\pi/3)$ or $(x_1 -2\pi/3, x_2-4\pi/3)$, we have
 $$
b_2 (x_1,x_2) =  - \frac{3}{2} + \frac{1}{2} ( x_1^2 - x_1 x_2 + x_2^2) + O(|x|^3).
$$
 Therefore, $C_{\kappa}$ is strictly convex if $ - 3/2 < \kappa < - 3/2 + \epsilon$.

 In view of 
 \S 3 of \cite{AndIsoMor}, we obtain :
 
 \begin{itemize}
 \item For the triangular lattice, 
$$
p(x,\lambda) = - \frac{1}{3}\big(b_2(x) + 3\lambda\big).
$$
Therefore, $M_{\lambda}$ is strictly convex for $-1< \lambda < -1 + \epsilon$, and $1/2 - \epsilon < \lambda < 1/2$.
 
 \item For the hexagonal  lattice, 
$$
p(x,\lambda) = - \frac{2}{9}\big(b_2(x) - \frac{9\lambda^2-3}{2}\big).
$$
Therefore,
$M_{\lambda}$ is strictly convex for $-1< \lambda < -1 + \epsilon$, $\lambda \in (-\epsilon,\epsilon)\setminus\{0\}$, and 
  $1 - \epsilon < \lambda < 1$.
  
  \item For the Kagome lattice, 
$$
p(x,\lambda) = \frac{1}{8}\big(\lambda - \frac{1}{2}\big)\big(b_2(x) - 8\lambda^2 - 4\lambda +1\big).
$$
Therefore, $M_{\lambda}$ is strictly convex for $-1 < \lambda < -1+\epsilon$, $\lambda \in (-1/4-\epsilon,-1/4+\epsilon)\setminus\{-1/4\}$ and $1/2-\epsilon < \lambda < 1/2$. 
  
   \item For the graphite, 
$$
p(x,\lambda) = \frac{1}{64}\big(b_2(x) - (8\lambda^2 + 4\lambda -1)\big)\big(b_2(x) - (8\lambda^2 - 4\lambda -1)\big).
$$
Therefore,
$M_{\lambda}$ is strictly convex for $- 1 < \lambda < -1 + \epsilon$, 
$\lambda \in (\pm 1/2-\epsilon, \pm 1/2 + \epsilon)\setminus\{\pm 1/2\}$, $(\pm 1/4-\epsilon,\pm 1/4+\epsilon)\setminus\{\pm 1/4\}$ and $1 - \epsilon < \lambda < 1$.
 \end{itemize}

%
%
%
%
%

\section{Reconstruction of scalar potentials}
\subsection{Parallelogram in the hexagonal lattice}
In this section, we  reconstruct a scalar potential from the D-N map on a bounded domain following the works \cite{CuMo90}, \cite{CuMoMo94}, \cite{Ober00}, \cite{IsMo2}.
This method depends strongly on the geometry of the lattice. Therefore, we explain it for the 2-dimensional
hexagonal lattice. First let us recall its structure. We identify ${\bf R}^2$ with ${\bf C}$, and put
$$
\omega = e^{\pi i/3}.
$$
 For $n = n_1 + i n_2 \in {\bf Z}[i]= {\bf Z} + i{\bf Z}$, 
let
$$
\mathcal L_0 = \left\{{\bf v}(n)\, ; \, n \in {\bf Z}[i]\right\}, \quad
{\bf v}(n) = n_1{\bf v}_1 + n_2{\bf v}_2,
$$
$$
{\bf v}_1 = 1 + \omega = (3 + \sqrt3i)/2, \quad {\bf v}_2 = \omega(1+ \omega) =  \sqrt3 i,
$$
$$
p_1 = \omega^{-1} = \omega^5, \quad p_2 = 1,
$$
and define the vertex set $\mathcal V_0$ by
$$
\mathcal V_0 = \mathcal V_{01} \cup \mathcal V_{02}, \quad \mathcal V_{0i} = p_i + \mathcal L_0.
$$
The adjacent points of $a_1 \in \mathcal V_{01}$ and $a_2 \in \mathcal V_{02}$ are defined by
\begin{equation}
\begin{split}
\mathcal N_{a_1} &= \{z \in {\bf C}\, ; \, |a_1 - z|=1\}\cap \mathcal V_{02} \\
&= \left\{a_1 + \omega, a_1 + \omega^{3}, a_1 + \omega^5\right\},
\end{split}
\nonumber
\end{equation}
\begin{equation}
\begin{split}
\mathcal N_{a_2} &= \{z \in {\bf C}\, ; \, |a_2 - z|=1\}\cap \mathcal V_{01} \\
&= \left\{a_2 + 1, a_2 + \omega^{2}, a_2 + \omega^4\right\}.
\end{split}
\nonumber
\end{equation}
Let $\mathcal D_0$ be the fundamental domain by the ${\bf Z}^2$-action (\ref{Zdaction}) on $\mathcal V_0$. It is a hexagon
 with 6 vertices 
 $\omega^k,\ 0 \leq k \leq 5$,
 with center at the origin. Take $D_N = \{n \in {\bf Z}[i]\, ; \, 0 \leq n_1 \leq N, \  0 \leq n_2 \leq N\}$, where $N$ is chosen large enough,  and put
 $$
 \mathcal D_N = {\mathop\cup_{n \in D_N}}\Big( \mathcal D_0 + {\bf v}(n)\Big).
 $$
 This is a parallelogram in the hexagonal lattice.
 
\begin{figure}[hbtp]
\includegraphics[width=9cm, bb=0 0 410 267]{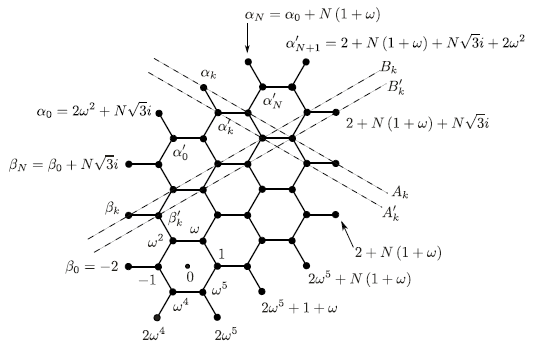}
\caption{Hexagonal parallelogram ($N = 2$)}
\label{S6HexaParallel}
\end{figure}
 
The interior angle of each vertex on the periphery of $\mathcal D_N$ is either $2\pi/3$ or $4\pi/3$. Let $\mathcal A$ be the set of the vertices with interior angle $2\pi/3$. 
 We regard $\mathcal D_N$ to be a subgraph of the original graph $\Gamma_0 = \{\mathcal L_0, \mathcal V_0, \mathcal E_0\}$,   and for  each $z \in \mathcal A$,  let $ e_{z,\zeta} \in \mathcal E_0$ be the outward edge emanating from $z$,  and  $\zeta = t(e_{z,\zeta})$ its terminal point. (See Figure \ref{S6HexaParallel} for the case $N = 2$.) 
Let $\Omega$ be the set of vertices of the resulting graph.
The boundary $\partial\Omega = \{t(e_{z,\zeta})\, ; \, z \in \mathcal A\}$ is  divided into 4 parts,   called  top, 
 bottom, right,  left sides, which are denoted by 
$(\partial\Omega)_T, (\partial\Omega)_B, (\partial\Omega)_R, (\partial\Omega)_L$, i.e. 
\begin{gather*}
\begin{split}
 (\partial\Omega)_T =& \{ \alpha_0 , \cdots , \alpha_N \}, \\
 (\partial\Omega)_B = & \{2\omega^5 + k(1 + \omega)\,  ; \, 0 \leq k \leq N\}, \\
 (\partial\Omega)_R = & \{ 2+ N(1+\omega ) + k\sqrt{3} i \, ; \, 1 \leq k \leq N \} \cup \{ 2+N(1+\omega ) +N\sqrt{3} i + 2 \omega ^2  \} , \\
 (\partial\Omega)_L =& \{2\omega^4\}\cup\{\beta_0,\cdots,\beta_N\} , 
\end{split}
\end{gather*}
where $ \alpha_k = \beta _N + 2\omega + k (1+\omega ) $ and $ \beta _k = -2 + k\sqrt{3} i $ for $ 0\leq k \leq N$.


\subsection{Matrix representation}
Our argument is close to that for the resistor network.
 To facilitate the comparison, we change the definition of the Laplacian as follows. 
\begin{equation}
\big(\widehat{\Delta}'_{\Omega}\widehat u)(v) = \frac{1}{{\rm deg}_{\Omega} (v)}\sum_{w\sim v}\left(\widehat u(v) - \widehat u(w)\right) = \widehat u(v) - \big(\widehat{\Delta}_{\Omega}\widehat u)(v).
\label{S6NewLaplacian}
\end{equation}
Putting
\begin{equation}
\widehat Q = \widehat V - \lambda -1,
\label{S6DefinewidehatQ}
\end{equation}
we consider the Dirichlet problem
\begin{equation}
\left\{
\begin{split}
& \big(\widehat{\Delta}'_{\Omega} + \widehat Q\big)\widehat u = 0 \quad {\rm in} \quad \stackrel{\circ}\Omega, \\
& \widehat u = \widehat f \quad {\rm on} \quad \partial\Omega.
\end{split}
\right.
\label{S6BVP}
\end{equation}
We assume the unique solvability of this equation. Then, the D-N map is defined by
\begin{equation}
\big(\Lambda_{\widehat Q}\widehat f\big)(v) = \partial_{\nu}^{\Omega}\widehat u(v).
\end{equation}

Given 4 vertices $z^{(0)} \in \, \stackrel{\circ}\Omega$, $z^{(0)} + \omega^k$, $z^{(0)} + \omega^{k+2}$, $z^{(0)}+\omega^{k+4} \in \Omega$, 
where $k=0$ or 1,  let us call $z^{(0)}$ the central point and $z^{(0)} + \omega^{k+2j}$ the peripheral point.
If $\widehat u$ satisfies $(\widehat \Delta'_{\Omega} + \widehat Q)\widehat u = 0$,  we have
\begin{equation}
\widehat u(z^{(0)}) - \frac{1}{3}\sum_{j=0}^2\widehat u(z^{(0)} + \omega^{k+2j}) + \widehat Q(z^{(0)})\widehat u(z^{(0)}) = 0.
\label{S6Equation4points}
\end{equation}
Therefore, we can compute the value at a peripheral point $\widehat u(z^{(0)} + \omega^{k + 2\ell})$ using the values at the central point $\widehat u(z^{(0)})$,  the other peripheral points $\widehat u(z^{(0)} + \omega^{k + 2j})$, $j \neq \ell$, and 
$\widehat Q(z^{(0)})$. Moreover, if we know the values of $\widehat u(z^{(0)})$ and $\widehat u(z^{(0)} + \omega^{k+2j})$, we can compute the potential 
$\widehat Q(z^{(0)})$ as long as $\widehat u(z^{(0)}) \neq 0$.

We split the vertex set of $\Omega$ into two parts:
$$
{\mathcal V}_0 = \, \stackrel{\circ}{\Omega} \, =  \{z^{(1)}, \cdots,z^{(\nu)}\}, \quad 
\mathcal V_1 = \partial\Omega = \{z^{(\nu+1)}, \cdots, z^{(\nu+\mu)}\}.
$$
Let us define matrices ${\bf D} = \left(d_{ij}\right)$, ${\bf A} = \left(a_{ij}\right)$ by
\begin{eqnarray*}
d_{ij} &  =& \left\{
\begin{split}
& 1, \quad {\rm if} \quad i = j, \\
& 0, \quad {\rm if} \quad  i \neq j,
\end{split}
\right. \\
a_{ij}  & = &
\left\{
\begin{split}
 \big( \mathrm{deg} _{\Omega} (z^{(i)} ) \big)^{-1} &, \quad {\rm if} \quad z^{(i)} \sim z^{(j)} \quad{\rm for} \quad z^{(i)} \in \, \stackrel{\circ}\Omega \quad {\rm or} \quad z^{(j)} \in \, \stackrel{\circ}\Omega, \\
 0& , \quad   {\rm if} \quad z^{(i)} \not\sim z^{(j)} \quad {\rm or} \quad  z^{(i)}, \ z^{(j)} \in \partial\Omega.
\end{split}
\right.
\nonumber
\end{eqnarray*}      
We define a $(\nu + \mu)\times (\nu + \mu)$ matrix ${\bf H}_0$ by
$$
{\bf H}_0 = {\bf D} - {\bf A} ,
$$ 
which corresponds to (\ref{S6NewLaplacian}). Since $\widehat Q$ is a scalar potential supported in
$\stackrel{\circ}\Omega$, it is identified with the diagonal matrix ${\bf Q} = \left(q_{ij}\right)$, where
\begin{equation}
q_{ij} = 
\left\{
\begin{split}
 \widehat Q(z^{(i)}), & \quad {\rm if} \quad i=j \leq \nu, \\
 0, & \quad {\rm if} \quad i \neq j \quad   {\rm or} \quad \ i =j \geq \nu + 1.
\end{split}
\right.
\nonumber
\end{equation}
We put
\begin{equation}
{\bf H} = {\bf H}_0 + {\bf Q}.
\end{equation}
In the following, $\widehat u(\mathcal V_i)$ denotes a vector in ${\bf C}^{\sharp\mathcal V_i}$,
and ${\bf H}(\mathcal V_i;\mathcal V_j)$ denotes a $\sharp\mathcal V_i \times \sharp\mathcal V_j$-submatrix of ${\bf H}$. 
For the solution $\widehat u$ to the boundary value problem (\ref{S6BVP}), noting $\mathrm{deg}_{\Omega} (v)=1 $ for $v\in \partial \Omega $, we rewrite the D-N map by
\begin{equation}
\big( \Lambda_{\widehat Q}\widehat f \big) (v) =  \sum_{w\in \stackrel{\circ}\Omega , w\sim v} \big( \widehat{f} (v)-\widehat{u} (w) \big) = \widehat{f} (v) + \big( \partial^{\Omega}_{\nu} \widehat{u} \big) (v) , \quad v\in \partial \Omega .
\label{DNnew}
\end{equation} 
Then, (\ref{S6BVP}) together with (\ref{DNnew})  is rewritten as 
the following system of equations
\begin{equation}
\left(
\begin{array}{cc}
{\bf H}(\mathcal V_0;\mathcal V_0) & {\bf H}(\mathcal V_0;\mathcal V_1) \\
{\bf H}(\mathcal V_1;\mathcal V_0) & {\bf H}(\mathcal V_1;\mathcal V_1)
\end{array}
\right)
\left(
\begin{array}{c}
\widehat u(\mathcal V_0) \\
\widehat f(\mathcal V_1)
\end{array}
\right) = 
\left(
\begin{array}{c}
0 \\
\Lambda_{\widehat Q}\widehat f
\end{array}
\right).
\label{S6matrixEq}
\end{equation}

It is easy to see that 
\begin{equation}
 0 \not\in \sigma(\widehat\Delta'_{\Omega}+ \widehat Q) \Longleftrightarrow
\det {\bf H}(\mathcal V_0;\mathcal V_0) \neq 0.
\label{S60notinsgmaequivdetnono0}
\end{equation}
 In fact, assume that 0 is not a Dirichlet eigenvalue of $\widehat\Delta'_{\Omega} + \widehat Q$, and ${\bf H}(\mathcal V_0;\mathcal V_0)\widehat u(\mathcal V_0)= 0$. Then, letting $\widehat u\big|_{\partial\Omega} = \widehat f = 0$, we see that $\widehat u$ satisfies 
 (\ref{S6BVP}). Since $0 \not\in \sigma(\widehat\Delta'_{\Omega} + \widehat Q)$, we have $\widehat u=0$, which implies $\det {\bf H}(\mathcal V_0;\mathcal V_0) \neq 0$. Conversely, suppose $\widehat u$ satisfies (\ref{S6BVP}) with $\widehat f = 0$. Then, (\ref{S6matrixEq}) is satisfied with $\widehat f=0$. 
Hence ${\bf H}(\mathcal V_0;\mathcal V_0)\widehat u(\mathcal V_0)=0$ and $\mathrm{det} {\bf H} ( \mathcal{V}_0 : \mathcal{V}_0 ) \not= 0$ imply $\widehat u(\mathcal V_0)=0$. This proves $0 \not\in \sigma(\widehat\Delta'_{\Omega} + \widehat Q)$.

From now on we assume that 

\medskip
\noindent
{\bf (D)} \ $ 0 \not\in \sigma(\widehat\Delta'_{\Omega}+ \widehat Q)$.

\medskip
\noindent
Hence the D-N map $\Lambda_{\widehat Q}$ has the following matrix representation
\begin{equation}
{\bf{\Lambda}}_{\bf{\widehat Q}} = {\bf{H}}(\mathcal V_1;\mathcal V_1) - {\bf H}(\mathcal V_1;\mathcal V_0){\bf H}(\mathcal V_0;\mathcal V_0)^{-1}{\bf H}(\mathcal V_0;\mathcal V_1).
\end{equation}

 The key to the inverse procedure is the following partial data problem.


\begin{lemma}\label{S6partialDNdata}
(1) Given a partial Dirichlet data $\widehat f$ on $\partial\Omega\setminus(\partial \Omega)_R$, and a partial Neumann data $\widehat g$ on $(\partial\Omega)_L$, there is a unique solution $\widehat u$ on $\stackrel{\circ}\Omega \cup (\partial\Omega)_R$ to the equation
\begin{equation}
\left\{
\begin{split}
& (\widehat\Delta'_{\Omega} + \widehat Q)\widehat u = 0 \quad {\rm in} \quad \stackrel{\circ}\Omega,\\
& \widehat u =\widehat f \quad {\rm on} \quad \partial\Omega\setminus(\partial\Omega)_R, \\
& \partial_{\nu}^{\Omega}\widehat u = \widehat g \quad {\rm on} \quad (\partial\Omega)_L.
\end{split}
\right.
\label{Lemma61Equation}
\end{equation}
(2) For subsets $A, B \subset \partial\Omega$, we denote the associated submatrix of $\bf{\Lambda_{\widehat Q}}$ by ${\bf\Lambda_{\widehat{\bf Q}}}(A;B)$. Then, the submatrix ${\bf \Lambda_{\widehat{\bf Q}}}((\partial\Omega)_{L};(\partial\Omega)_R)$ is non-singular, i.e. 
$$
{\bf\Lambda_{\widehat{\bf Q}}}((\partial\Omega)_{L};(\partial\Omega)_{R}) : (\partial\Omega)_R \to 
(\partial\Omega)_{L}
$$ 
is a bijection. \\
\noindent
(3) Given the D-N map ${\bf\Lambda_{\widehat{\bf Q}}}$, a partial Dirichlet data $\widehat f_2$ on $\partial\Omega\setminus(\partial\Omega)_R$ and a partial Neumann data $\widehat g$ on $(\partial\Omega)_{L}$, there exists a unique $\widehat f$ on $\partial\Omega$ such that $\widehat f = \widehat f_2$ on $\partial\Omega\setminus(\partial\Omega)_R$ and ${\bf\Lambda_{\widehat{\bf Q}}}\widehat f = \widehat g$ on $(\partial\Omega)_{L}$.
\end{lemma}

Proof. (1) Look at Figure \ref{S6HexaParallel}. The values of $\widehat u(x_1 + ix_2)$ at $\omega^4$ and on the line $x_1 = -1$ are computed from the D-N map and the values of $\widehat f$, $\widehat g$.
Uisng the equation (\ref{S6Equation4points}), one can then compute $\widehat u(x_1 +ix_2)$ on $\omega^5 $ and the line $x_1=-1/2$.  (For the line $x_1 = -1/2$, start from $\omega^2$ and go up).   This and the  Dirichlet data $\widehat f(x_1+ix_2)$ at $2\omega^5$ give $\widehat u(x_1 +ix_2)$ on $1$ and the line $x_1= 1/2$.  
Repeating this procedure, we get $\widehat u(z)$ for all $z \in \Omega$. 

 (2)  Suppose $\widehat f =0$ on $\partial\Omega\setminus(\partial\Omega)_R$ and $\Lambda_{\widehat{\bf Q}}\widehat f=0$ on $(\partial\Omega)_L$.
By (1), the solution $\widehat u$ vanishes identically. Hence $\widehat f=0$ on $(\partial\Omega)_R$. This proves the injectivity, hence the surjectivity. 

(3)  We seek $\widehat f$ in the form
$$
({\bf{\Lambda_{\widehat{\bf Q}}}} \widehat f)\big|_{(\partial\Omega)_{L}}= 
{\bf{\Lambda_{\widehat{\bf Q}}}}((\partial\Omega)_{L};(\partial\Omega)_R)\widehat f_1 + {\bf\Lambda_{\widehat{\bf Q}}}((\partial\Omega)_{L};\partial\Omega\setminus(\partial\Omega)_R)\widehat f_2 = \widehat g,
$$
where $\widehat f_1 = \widehat f\big|_{(\partial\Omega)_R}$. 
By (2), we have only to take
$$
\widehat f_1 = \left({\bf\Lambda_{\widehat{\bf Q}}}((\partial\Omega)_{L};
(\partial\Omega)_R)\right)^{-1}\left(\widehat g - {\bf\Lambda_{\widehat{\bf Q}}}((\partial\Omega)_{L};\partial\Omega\setminus(\partial\Omega)_R)\widehat f_2\right). \qed
$$
\begin{figure}[htbp]
 \begin{minipage}{0.40\hsize}
  \begin{center}
\includegraphics[width=48mm, bb=0 0 214 215]{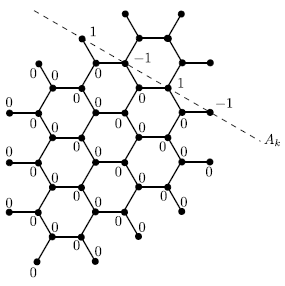}
  \end{center}
  \caption{Line $A_k$}
  \label{S6LineAk}
 \end{minipage}
 \begin{minipage}{0.59\hsize}
  \begin{center}
   \includegraphics[width=80mm, bb=0 0 358 176]{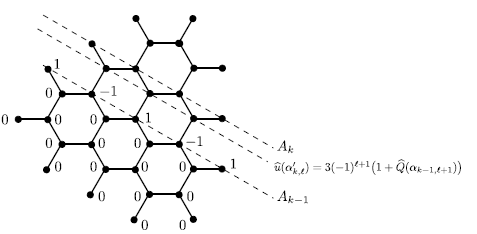}
  \end{center}
  \caption{Line $A_k'$}
  \label{S6LineAkprime}
 \end{minipage}
\end{figure}
%
%
%

\medskip
Now, for $0 \leq k \leq N$, let us consider a {\it diagonal} line $A_k$ :
\begin{equation}
A_k = \{x_1 + ix_2\, ; \, x_1 + \sqrt3x_2 = a_k\},
\label{S6Dinagonalline}
\end{equation}
where $a_k$ is chosen so that $A_k$ passes through
\begin{equation}
\alpha_k =  \alpha_0 + k (1+ \omega )  \in 
(\partial\Omega)_R.
\end{equation} 
The vertices on $A_k\cap\Omega$ are written as
\begin{equation}
\alpha_{k,\ell} = \alpha_k + \ell (1 + \omega^5 ), \quad \ell = 0, 1, 2, \cdots.
\end{equation}
 We also need another diagonal line $A_k'$ between $A_k$ and $A_{k-1}$ : 
\begin{equation}
A_k' = \{x_1 + ix_2\, ; \, x_1 + \sqrt3x_2 = a_k'\},
\label{S6Dinagonallineprime}
\end{equation}
where $a_k'$ is such that $A_k'$ passes through 
\begin{equation}
\alpha_k' = \alpha_k + \omega ^5 .
\end{equation} 
The vertices on $A_k'\cap \Omega$ are written as 
\begin{equation}
\alpha_{k,\ell}'  = \alpha'_k + \ell (1 + \omega^5 ), \quad \ell = 0, 1, 2, \cdots.
\end{equation} 
Finally, we let 
\begin{equation}
A_{N+1}' = \{x_1 +ix_2\, ; \, x_1 + \sqrt3 x_2 = a_{N+1}'\},
\end{equation}
which passes through
\begin{equation}
\alpha_{N+1}'  = \alpha_N + 2 = 2 + N(1 + \omega) + N\sqrt 3i + 2\omega^2.
\end{equation}


\begin{lemma}
\label{alphareconstruct}
(1) Let $A_k \cap \partial\Omega = \{\alpha_{k,0}, \alpha_{k,m}\}$.  
Then there exists a unique solution  $\widehat u$  to the equation
$$
(\widehat\Delta'_{\Omega} + \widehat Q)\widehat u = 0 \quad {\rm in} \quad
\stackrel{\circ}\Omega,
$$
 with partial Dirichlet data $\widehat f$ such that
\begin{equation}
\left\{
\begin{split}
& \widehat f(\alpha_{k,0}) = 1, \\
& \widehat f(\alpha_{k,m}) = (-1)^m, \\
& \widehat f(z) = 0 \quad {\rm for} \quad z \in \partial\Omega\setminus
\left((\partial\Omega)_R\cup\{\alpha_{k,0}\cup\alpha_{k,m}\}\right),
\end{split}
\right.
\label{Lemma6.4Equation}
\end{equation}
and partial Neumann data $\widehat g=0$ on $(\partial\Omega)_{L}$. It satisfies
\begin{equation}
\widehat u(x_1 + ix_2) = 0 \quad {\rm if} \quad  x_1 + \sqrt3x_2 < a_k,
\label{u=0belowx1+3x2<ck}
\end{equation}
and on $x_1 + \sqrt3 x_2 = a_k$, 
\begin{equation}
\widehat u(\alpha_{k,\ell}) = (-1)^{\ell}, \quad \ell = 0, 1, 2, \cdots.
\label{u=pm1onx1+3x2}
\end{equation}
(2) Using the solution $\widehat u$ for the data (\ref{Lemma6.4Equation}) with $k$ replaced by $k-1$, $\widehat Q(\alpha_{k-1,\ell+1})$ is computed as
\begin{equation}
\widehat Q(\alpha_{k-1,\ell+1}) = 
\frac{\widehat u(\alpha'_{k,\ell})}{3(-1)^{\ell+1}} - 1, \quad \ell =0, 1, 2, \cdots.
\label{Qalphakellvalue}
\end{equation}
\end{lemma}

Proof. The uniqueness of $\widehat u$ follows from Lemma \ref{S6partialDNdata}. To prove the existence, we argue as in the proof of Lemma \ref{S6partialDNdata} (1). By the equation (\ref{S6Equation4points}) with central point below $A'_k$ and the condition on $\widehat f$, $\widehat g$, one can compute $\widehat u(x_1 + ix_2)$ successively to obtain (\ref{u=0belowx1+3x2<ck}). 
 Using  again (\ref{S6Equation4points}), putting central point on $A'_k$, we obtain (\ref{u=pm1onx1+3x2}).  
 
We replace $k$ by $k-1$ in the above procedure. Then, $\widehat u(z)$ is computed as 
$$
\widehat u(x_1 + ix_2) = \left\{
\begin{split}
 0 & ,\quad {\rm if} \quad x_1 + \sqrt3 x_2 < a_{k-1}, \\
 (-1)^{\ell} &, \quad {\rm if} \quad x_1 + ix_2 = \alpha_{k-1,\ell}.
\end{split}
\right.
$$
We use (\ref{S6Equation4points}) with central point $\alpha_{k-1,\ell +1} \in A_{k-1}$. Then, 
$$
\big(1 + \widehat Q(\alpha_{k-1,\ell+1})\big)\widehat u(\alpha_{k-1,\ell+1}) = 
\frac{1}{3}\widehat u(\alpha'_{k,\ell}),
$$
which shows
(\ref{Qalphakellvalue}). \qed

\begin{figure}[hbtp]
\includegraphics[width=6cm, bb=0 0 275 246]{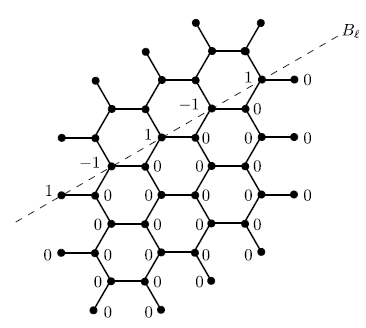}
\caption{Line $B_{\ell}$}
\label{S6LineBell}
\end{figure}

\medskip
Let us exchange the roles of $(\partial\Omega)_R$, $(\partial\Omega)_{L}$ and 
$(\partial\Omega)_T$, $(\partial\Omega)_B$.  For $0 \leq \ell \leq N$, consider a diagonal line $B_{\ell}$ 
\begin{equation}
B_{\ell} = \{x_i + ix_2\, ; \, x_1 - \sqrt3x_2 = b_{\ell}\},
\end{equation}
where $b_{\ell}$ is chosen so that $B_{\ell}$ passes through 
\begin{equation}
\beta_{\ell} = -2 + \ell\sqrt3 i \in (\partial\Omega)_{L}.
\end{equation}
The vertices on $B_{\ell}\cap\Omega$ are written as
\begin{equation}
\beta_{k,\ell} = \beta_{\ell} + k(1 + \omega), \quad k = 0, 1, 2, \cdots.
\end{equation}
Another diagonal line is 
\begin{equation}
B'_{\ell} = \{x_i + ix_2\, ; \, x_1 - \sqrt3x_2 = b'_{\ell}\},
\end{equation}
where $b'_{\ell}$ is chosen so that $B'_{\ell}$ passes through 
\begin{equation}
\beta'_{\ell} = - 1 + \ell\sqrt3 i.
\end{equation}
The vertices on $B'_{\ell}\cap\Omega$ are written as
\begin{equation}
\beta_{k,\ell}' = \beta'_{\ell} + k(1 + \omega), \quad k = 0, 1, 2, \cdots.
\end{equation}
Finally, we put 
\begin{equation}
B'_{N+1} = \{x_1 + ix_2\, ; \, x_1 - \sqrt3 x_2 = b'_{N+1}\},
\end{equation}
which passes through $(\partial\Omega)_T$.

Then, the following lemma is proven in the same way as above.


\begin{lemma}
\label{betareconstruct}
(1) If $B_{\ell}\cap {\partial\Omega} = \{\beta_{0,\ell}\}$, take the Dirichlet data $\widehat f$ such that
\begin{equation}
\left\{
\begin{split}
& \widehat f(\beta_{0,\ell}) = 1, \\
& \widehat f(z) = 0 \quad {\rm for} \quad z \in \partial\Omega\setminus
\left((\partial\Omega)_T\cup\{\beta_{0,\ell}\}\right).
\end{split}
\right.
\label{Lemma6.5equation}
\end{equation}
 If $B_{\ell}\cap {\partial\Omega} = \{\beta_{0,\ell}, \beta_{m,\ell}\}$, take the Dirichlet data $\widehat f$ such that
\begin{equation}
\left\{
\begin{split}
& \widehat f(\beta_{0,\ell}) = 1, \\
& \widehat f(\beta_{m,\ell}) = (-1)^m, \\
& \widehat f(z) = 0 \quad {\rm for} \quad z \in \partial\Omega\setminus
\left((\partial\Omega)_T\cup\{\beta_{0,\ell}, \beta_{m,\ell}\}\right).
\end{split}
\right.
\label{Lemma6.5equation2}
\end{equation}
Then, there exists a unique solution $\widehat u$  to 
$$
(\widehat\Delta'_{\Omega} + \widehat Q)\widehat u = 0 \quad {\rm in} \quad
\stackrel{\circ}\Omega,
$$
 with partial Dirichlet data $\widehat f$ 
and partial Neumann data $\widehat g=0$ on $(\partial\Omega)_{B}$. It satisfies
\begin{equation}
\widehat u(x_1 + ix_2) = 0 \quad {\rm if} \quad  x_1 - \sqrt3x_2 > b_{\ell},
\label{u=0belowx1-3x2<bk}
\end{equation} 
and on $x_1 - \sqrt3 x_2 = b_{\ell}$, 
\begin{equation}
\widehat u(\beta_{k,\ell}) = (-1)^{k}, \quad \beta_{k,\ell} = \beta_{\ell} + k(1 + \omega), \quad k = 0, 1, 2,\cdots.
\label{u=pm1onx1-3x2}
\end{equation}
(2) Using the solution $\widehat u$ for the data (\ref{Lemma6.5equation}) with $\ell$ replaced by $\ell -1$, $\widehat Q(\beta_{k+1,\ell-1})$ is computed as
\begin{equation}
\widehat Q(\beta_{k+1,\ell-1}) = \frac{\widehat u(\beta'_{k,\ell})}{3(-1)^{k+1}} -1, \quad k =0, 1, 2, \cdots.
\label{Qbetakellvalue}
\end{equation}
\end{lemma}


\subsection{Reconstruction algorithm}

We are now in a position to give an {\it algorithm for the reconstruction} of the potential. First let us note that given the boundary data $\widehat f$ and the D-N map, one can compute the values of $\widehat u$ on the points adjacent to $\partial\Omega$.

\medskip
\noindent
(1)  Use Lemma \ref{alphareconstruct} to construct the data $\widehat f$ and the solution $\widehat u$ with $k=N$. Use the equation, and 
\begin{itemize}
\item 
the fact that $\widehat u=0$ in the region $x_1 + \sqrt3 x_2 < a_N$, 
\item the value $\widehat f(\alpha'_{N+1})$, 
\item 
 the values $1, -1, 1$ of $\widehat u$ on the line $A_N$, 
\end{itemize}
compute the value of $\widehat Q$ at $A_N \cap \stackrel{\circ}\Omega = \{\alpha'_{N+1} + \omega^4\}$.
\begin{figure}[hbtp]
\includegraphics[width=13cm, bb=0 0 577 412]{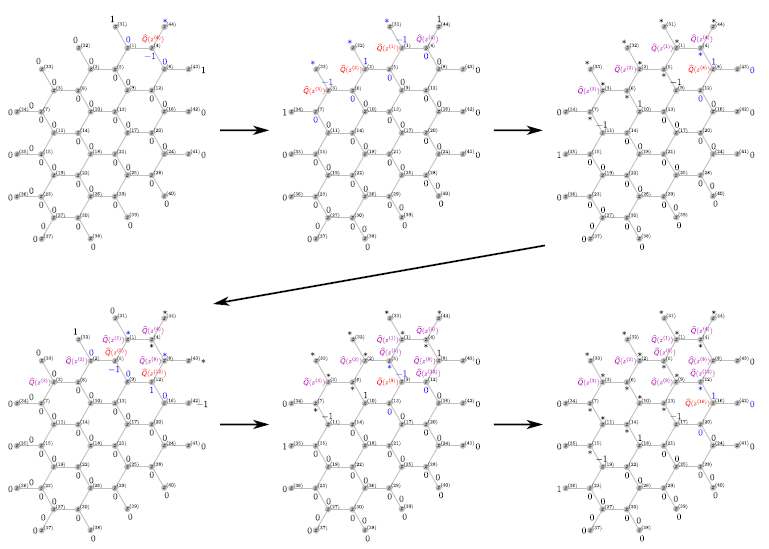}
\caption{Reconstruction of the potential in the hexagonal lattice}
\label{HexagonalPotentialReconstruction}
\end{figure}

\begin{figure}[hbtp]
\includegraphics[width=13cm, bb=0 0 577 362]{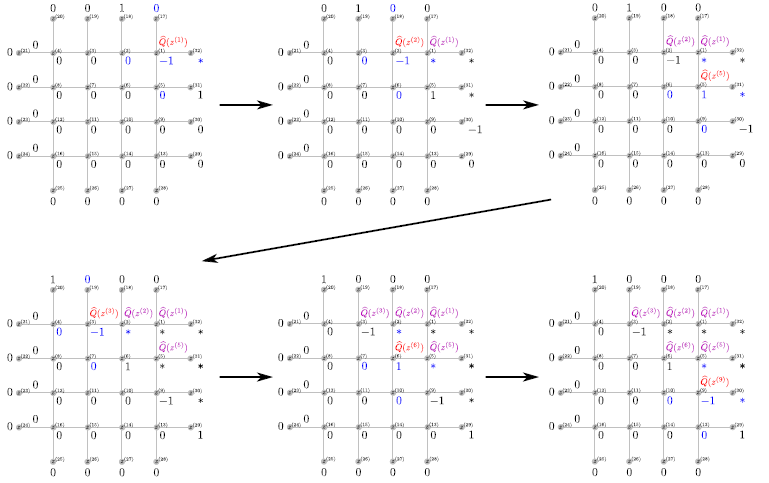}
\caption{Reconstruction of the potential in the square lattice}
\label{SquarePotentialReconstruction}
\end{figure}

\medskip
\noindent
(2) Use Lemma \ref{betareconstruct} to construct the data $\widehat f$ and the solution $\widehat u$ with $\ell =N$. Determine the values of $\widehat Q$ on $B_N$  by the argument similar to the one in (1).

\medskip
\noindent
(3) Assume that all the values of $\widehat Q$ on $\{x_1 + \sqrt3 x_2 > a_k\}$ are computed. Use Lemma \ref{alphareconstruct} to construct the data $\widehat f$ and the solution $\widehat u$ which takes values $1, -1, 1, \cdots$ on $A_k$. Use the equation, $\widehat f$, the D-N map and the values of $\widehat Q$, compute the values $\widehat u$ in the region $\{x_1 + \sqrt3x_2 > a_k\}$. Then, calculate $\widehat Q$ on $A_k$ using the equation.

\medskip
\noindent
(4) Assume that all the values of $\widehat Q$ on $\{x_1 + \sqrt3 x_2 \geq a_k\} \cap \{x_1 - \sqrt{3}x_2 < b_{\ell}\}$ are computed. 
Use Lemma \ref{betareconstruct}  to construct the data $\widehat f$ and the solution $\widehat u$ which  takes values $\pm 1$ on $B_{\ell}$. Use the equation to compute the value of $\widehat Q$ at $A'_k\cap B_{\ell}$. This makes it possible to compute $\widehat Q$ on $A'_k$.

\medskip
\noindent
(5) Repeat the above procedure until $k=0$. 

\medskip
\noindent
(6) Rotate the domain, and compute the values of $\widehat Q$ at the remaining points by the same procedure as above.

\medskip
The above  reconstruction procedure and that for the case of the square lattice are illustrated in the Figures \ref{HexagonalPotentialReconstruction} and \ref{SquarePotentialReconstruction}.


\section{Inverse problems for resistor networks}
In the previous section, we studied the inversion procedure for the scalar potential from the S-matrix. In this section, we consider the inverse problems for the electric conductivity and graph structure from the S-matrix.This should be compared with the perturbation of Riemannian metric or domain for the case of continuous model.


\subsection{Inverse boundary value problem for resistor network}
\label{SubsectionRegistornetwork}
A {\it circular planer graph} $G = \{\mathcal V,\mathcal E\}$ is a graph which is imbedded in a disc $D \subset {\bf R}^2$ so that 
its boundary $\partial\mathcal V$ lies on the circle $\partial D$ and its interior $\stackrel{\circ}{\mathcal V}$ is in the topological interior of $D$.
A conductivity on $G$ is a positive function $\gamma$ on the edge set $\mathcal E$. For $e \in \mathcal E$, the value $\gamma(e)$ is 
called the {\it conductance} of $e$, and its inverse $1/\gamma(e)$ the {\it resistance} of $e$. Equipped with the conductance, the graph $G $ is called the 
{\it resistor network}. The Laplacian $\widehat\Delta_{res}$ is defined by
\begin{equation}
\big(\widehat\Delta_{res}\widehat u\big)(v) = \sum_{w\in\mathcal N_v}\gamma(e_{vw})\Big(\widehat u(w) - \widehat u(v)\Big),
\end{equation}
where $e_{vw}$ is an edge with end points $v, w$. 
Then, for any boundary value $\widehat f$, there exists a unique $\widehat u$ satisfying
\begin{equation}
\left\{
\begin{split}
& \widehat\Delta_{res}\widehat u = 0 \quad {\rm in} \quad \stackrel{\circ}{ \mathcal V}, \\
& \widehat u = \widehat f \quad {\rm on} \quad \partial \mathcal V.
\end{split}
\right.
\label{S5resBVP1}
\end{equation}
The D-N map $\Lambda_{res}$ is defined by
\begin{equation}
\Lambda_{res}\widehat f = - \widehat\Delta_{res}\widehat u \quad {\rm  on } \quad \partial{\mathcal V},
\label{S5DN1}
\end{equation}
where $\widehat u$ is the solution to (\ref{S5resBVP1}). (See \cite{Col94}, \cite{CuMo00}.)


We compare inverse problems for resistor networks with inverse conductivity problems for continuous models.
 For a compact manifold $M$ with boundary $\partial M$, the D-N map defined on $\partial M$ is invariant by any  diffeomorphism on $M$  leaving $\partial M$ invariant. It is shown that in 2-dimensions the D-N map determines the conductivity up to these diffeomorphisms. For higher dimensions, it is proven under the additional assumption of real-analyticity. Also it is known that the scattering relation 
(which is defined through the geodesic in $M$ with end points on $\partial M$) determines  the {\it simple} manifold (\cite{Muho77}, \cite{PesUhl}). Here, the simple manifold is a compact Riemannian manifold with strictly convex boundary,  whose exponential map ${\rm exp}_x : {\rm exp}_x^{-1}(\overline{M}) \to \overline{M}$ is a diffeomorphism.
 
 The above mentioned issues in differential geometry have counter parts in the resistor network. Instead of the diffeomorphism, one uses the notion of {\it criticality} and {\it elementary transformation}.  
The geodesics is also extended on the graph. {\it Connection} is defined  to introduce a mapping between two parts of the boundary.  First let us recall these notions.

A {\it path}  from $v \in \partial{\mathcal V}$ to $w \in \partial{\mathcal V}$ is a sequence of edges : $e_i = e_i(v_i,v_{i+1}), i = 0,\cdots, n-1$, such that 
$v_0 = v$, $v_n= w$,  and $v_i \ (1 \leq i \leq n-1)$ are distinct vertices in $\stackrel{\circ}{\mathcal V}$.  Two sequences of boundary vertices $P = (p_1,\cdots,p_k), Q = (q_1,\cdots,q_k)$ are said to be {\it connected} through $G $ 
if there exists a  path $c_i$ from $p_i$ to $q_i$, $i = 1,\cdots,k$, moreover
$c_i$ and $c_j$ have no vertex in common if $i \neq j$. In this case, the pair $(P,Q)$ is said to be a $k$-{\it connection}.
The graph $G$ is said to be {\it critical} if one removes any edge in $\mathcal E$, there exist  $k$ and a $k$-connection $(P,Q)$ in $G$ which is no longer
 connected in the resulting graph $G '$ after the removal of the edge. Here, to {\it remove} the edge means either to {\it delete} the edge leaving  end points as two vertices, or to {\it contract} the edge letting two end points together into one vertex (see \cite{CuMo00}, p. 16).  An {\it arc} $c(t), 0 \leq t \leq 1$, in $G$ is a union of edges such that $c(t)$ is  continuous, and $c(t) \in \partial\mathcal V$ if and only if $t = 0, 1$.


\begin{figure}[hbtp]
\includegraphics[width=8cm, bb=0 0 358 414]{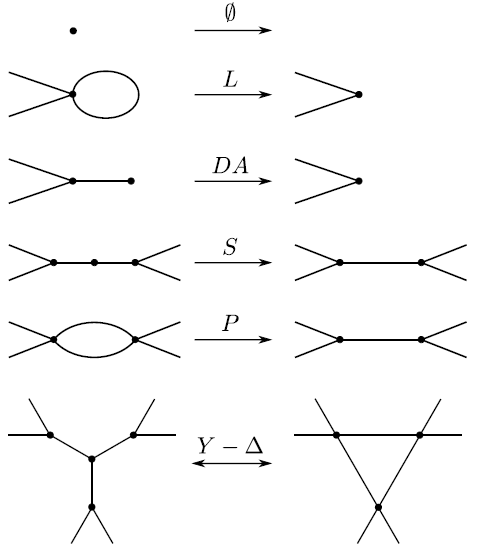}
\caption{Elementary transformations}
\label{ElemenTransf}
\end{figure}


 The elementary transformations of a planar graph consist of the following 6 operations. 
 
\begin{itemize}
\item ($\emptyset$) ({\it Point}) means to remove an isolated vertex. 
 
\item
(L) ({\it Loop}) means to remove a loop.

\item
(DA) ({\it Dead arm}) means to remove an edge with end point of degree 1. 

\item 
(S) ({\it Series}) means to remove a vertex $a$ ($a \not\in \partial{\mathcal V}$) of degree 2 and the two adjacent edges and to join the neighboring vertices $b$, $c$ of $a$ by one edge. The conductance of the edge $e_{bc}$ is defined by
\begin{equation}
\gamma(e_{b,c}) = \left(\gamma(e_{b,a})^{-1} + \gamma(e_{a,c})^{-1}\right)^{-1}.
\end{equation}
\item
(P) ({\it Pararelle}) means to replace two edges joining the same vertices  $a$ and $b$ by one edge. The conductance is defined by
\begin{equation}
\gamma_{a,b} = \gamma'_{a,b} + \gamma''_{a,b},
\end{equation}
where $\gamma'_{a,b}$ and $\gamma''_{a,b}$ are the conductances of the double edges.

\item
($Y - \Delta$) means to replace a star with 3 branches of center $0$ and edges $(0,1)$, $(0,2)$, $(0,3)$ by one triangle of the vertices 1, 2, 3, and remove the center 0. The conductances are defined by
\begin{equation}
\gamma(e_{i,j}) = \frac{\gamma(e_{0,i})\gamma(e_{0,j})}{\gamma(e_{0,1}) + \gamma(e_{0,2}) + \gamma(e_{0,3})}.
\end{equation}
\end{itemize}

The meaning of these graph operations will be clear by Figure \ref{ElemenTransf}. 

\begin{remark}
Note that the elementary transformations increase neither the number of vertices nor that of edges. Hence, it does not increase the number of arcs. 
\end{remark}

Using these notions, it is shown that 


\begin{theorem}
\label{TheoremGraphInv}
(1)  Any circular planar graph is transformed to a critical graph, which is unique within elementary transformations. \\
\noindent
(2)  Any critical circular planar graph is uniquely determined by its D-N map 
up to elementary transformations. \\
\noindent
(3)  Any two circular planar equivalent graphs have the same number of arcs.
 \end{theorem}
 
For the proof, see Corollary 9.4 and Theorem 9.5 in p. 168 of \cite{CuMo00}, and 
Theorem 4 in p. 146 of \cite{ColGit96},
where the reconstruction procedure of the graph is also explained.


\subsection{From the S-matrix to the D-N map for resistor network}
\subsubsection{Conductivity problem}
Let us return to the inverse scattering problem.
Assuming that there is no defect, whose meaning will be given in  \ref{SubsectionDefects}, we can formulate the perturbation of conductivity in the form $\widehat V$ in the assumption { (B-4)}. 
By the arguments in \S 4, the inverse scattering problem is reduced to the inverse boundary value problem in a bounded domain. Note that from the S-matrix of energy $\lambda$, we obtain the D-N map with energy $\lambda$ associated with the equation (\ref{IntBVP}).

\subsubsection{Defect problem} 
We need to pay attention in formulating the inverse boundary value problem in a domain with defect in terms of the resistor network problem. We explain it more precisely.
Suppose that we are given a periodic lattice $\Gamma_0 = \{ \mathcal{L}_0 , \mathcal{V}_0 , \mathcal{E}_0 \} $ as in \S 2 satisfying the assumptions (A-1) $\sim$ (A-4).  We perturb $\Gamma_0$ locally, and let the resulting graph $\Gamma = \{\mathcal V, \mathcal E\}$ satisfy the assumptions (B-1) $\sim $ (B-4).
Let $  \mathcal{V}_{int}$ be the associated interior domain.
Assume that $\mathcal V_{int}$ is a planar graph in the sense of Subsection \ref{SubsectionRegistornetwork}, and denote it by $\mathcal V_{def}$.
To make this perturbation consistent with the previous arguments, we assume that 
$$
\sharp \{ w \in \, \stackrel{\circ}{\mathcal{V}_{def}} \ ; \ w \sim v \}  ,\quad v\in \partial \mathcal{V}_{def} ,
$$
is a constant on $\partial \mathcal{V}_{def} $ which is denoted by $\mu_0$.
In order to apply Theorem \ref{TheoremGraphInv} to our problem, we take 
$$
\gamma (e_{vw} )=1 .
$$
Then we have 
$$
- \widehat{\Delta}_{res} = \mathrm{deg} _{\mathcal{V}_{def}} (\cdot ) ( \widehat{\Delta}_{\Gamma} - 1)  \quad \text{in} \quad \stackrel{\circ}{\mathcal{V}_{def}} .
$$
Hence the equation (\ref{S5resBVP1}) is rewritten as 
\begin{equation}
\left\{
\begin{split}
& ( - \widehat{\Delta}_{\Gamma} +1 ) \widehat{u} =0 \quad \text{in} \quad \stackrel{\circ}{\mathcal{V}_{def}} , \\
& \widehat u = \widehat f \quad {\rm on} \quad \partial \mathcal{V}_{def} ,
\end{split}
\right.
\label{S7BVPwithV}
\end{equation}
Note that $ -\widehat{\Delta}_{\Gamma} $ is self-adjoint on $\ell^2 ( \stackrel{\circ}{\mathcal{V}_{def}} )$ equipped with the inner product
$$
( \widehat{f} , \widehat{g} ) _{\ell^2 (\stackrel{\circ}{\mathcal{V} _{def} } )} = \sum_{v\in \stackrel{\circ}{\mathcal{V} _{def} } }  \widehat{f} (v) \overline{\widehat{g} (v)} \mathrm{deg} _{\mathcal{V}_{def}} (v)   ,
$$
where $\widehat{f} (v)= \widehat{g} (v)=0$ for any $v\in \partial \mathcal{V} _{def} $.  This operator is denoted by $-\widehat{\Delta}_{\Gamma,\mathcal{V} _{def}}$.

\begin{lemma}
We have
$ \ - 1 \not\in \sigma(-\widehat{\Delta}_{\Gamma,\mathcal{V} _{def}})$.
\label{S7_lem_eigenvaluemu0}
\end{lemma}

Proof. 
Suppose  $-1\in \sigma(-\widehat{\Delta}_{\Gamma,\mathcal{V} _{def}})$.
Then there exists a function $\widehat{u}$ satisfying the equation
\begin{equation*}
\left\{
\begin{split}
& \widehat{\Delta}_{res} \widehat{u} =0 \quad \text{in} \quad \stackrel{\circ}{\mathcal{V}_{def}} , \\
& \widehat u = 0 \quad {\rm on} \quad \partial \mathcal{V}_{def} ,
\end{split}
\right.
\end{equation*}
However, the maximum principle for harmonic functions associated with $\widehat{\Delta}_{res}$ implies that $\widehat{u}$ must vanish in $\stackrel{\circ}{\mathcal{V}_{def}}$ (see Theorem 3.2 and Corollary 3.3 in \cite{CuMo00}).
This is a contradiction.
\qed

\medskip
The results in \S 2 $\sim$ \S 5 also hold for this perturbed system. To study $\mathcal V_{def}$, in view of (\ref{S7BVPwithV}), we have to consider the operator $ -\widehat{\Delta}_{\Gamma,\mathcal V_{def}} $ and its D-N map with $\lambda = -1 $.
However, we have $-1 \in \sigma_e (\widehat{H} ) \cap \mathcal{T}$ in some examples of lattices satisfying the above assumptions.
For these reasons, we pick up some cases in which we can compute the D-N map $\Lambda _{int} (\lambda )$ from the S-matrix of $\widehat{H} = -\widehat{\Delta}_{\Gamma} $.

\medskip

\textit{Case 1.}
If $ \lambda \in \sigma_e (\widehat{H} ) \setminus \mathcal{T} $, there is no problem, and we can apply our previous arguments.

\medskip

\textit{Case 2.}
Let $\lambda \not\in \sigma(-\widehat{\Delta}_{\Gamma,\mathcal{V} _{def}}) $. Take any open set $ \mathcal{O} \subset \sigma_e (\widehat{H} ) \setminus \mathcal{T} $. 
If we are given the S-matrix for all energy $\mu \in \mathcal{O} $, we can compute the D-N map for $-\widehat{\Delta}_{\Gamma,\mathcal V_{def}} - \mu $ with $\mu \in \mathcal{O}\setminus\sigma(-\widehat{\Delta}_{\Gamma,\mathcal{V} _{def}})$.
Since the D-N map $\Lambda_{int}(\mu)$ is meromorphic with respect to the energy $\mu$, using the analytic continuation and taking $\mu = \lambda  $, we can compute the D-N map $\Lambda_{int}(\lambda)$ for $\widehat{\Delta}_{res} $.

\medskip
\textit{Case 3.}
In practical applications, it often happens that $ \lambda $ is an end point of $\sigma_e (\widehat{H} )$ as above.
For example, this is the case for the perturbation of the hexagonal lattice.
In view of Lemma \ref{S7_lem_eigenvaluemu0}, the D-N map for $ - \widehat{\Delta} _{\Gamma,\mathcal V_{def}} - \mu$ is continuous with respect to $\mu$ when $\mu $ is close to $\lambda $.
Therefore, choosing a sequence $\mu_j $ convergent to $ \lambda $, one can compute the D-N map $\Lambda_{int} (\lambda )$ from the S-matrices $S(\mu_j ) $ for $j\geq 1$.

\medskip
The above arguments have a general character and work for many lattices satisfying our assumptions (A-1) $\sim$ (A-4) and (B-1) $\sim$ (B-4). 
Thus, when we perturb a bounded part of these lattices by a planer network,
\begin{center}
\textit{we can determine the perturbation as a planer network}
\end{center}
by using Theorem \ref{TheoremGraphInv}.

Main barriers for this fact are 
the Rellich type theorem and the unique continuation property for the associated spectral problems. 
In Theorem 5.10 of \cite{AndIsoMor}, we summarized examples of the lattices having this property.
Therefore, when we perturb a bounded part of square, triangular, and hexagonal lattices by removing a finite number of edges in such a way that the unique continuation property holds in the remaining part (see \cite{IsMo1} and the Figures 11 and 12 in \cite{AndIsoMor}),
we can determine the network.

As is seen from the definition, the elementary transformations are of topological nature, and its physical realization is not an obvious problem.
Therefore we must be careful in the application of above results, which we shall discuss in the next subsection.

 
\subsection{Inverse resistor network problem in the hexagonal lattice}
\label{Subsection7.3}

\subsubsection{Conductivity}
To study the conductivity problem, in the assumptions (B-1) $\sim$ (B-4), we take $\mathcal V_{int}$ to be a sufficiently large planar graph so that the conductivity is constant on $\mathcal V_{ext}$. Then by Theorem \ref{TheoremGraphInv}, the D-N map (\ref{S5DN1}) determines the Laplacian on $\mathcal V_{int}$ as a planar network. Since the S-matrix with energy $\lambda$ determines the associated D-N map, taking note  that the S-matrix and the D-N map are analytic with respect to the energy $\lambda$, we can compute the D-N map (\ref{S5DN1}) from the S-matrix for all energies. 
We have thus proven the following theorem.

\begin{theorem}
The conductivity of the periodic hexagonal lattice is determined by the S-matrix of all energies.
\end{theorem}

\subsubsection{Defects}
\label{SubsectionDefects}
By \textit{defects} we mean to delete some edges and to remove isolated vertices.
Consider, for example, the simplest case in which only one edge, say the edge with end points $a$ and $d$ in Figure \ref{Probing}, is removed.
Our main idea of detecting defects consists in using the solution $\widehat u$ given in Lemmas \ref{alphareconstruct} and \ref{betareconstruct} with $\widehat Q = - \lambda$. As will be discussed in the proof of Lemma \ref{LambdaHdeffk=LambdaHofkLemma}, one can detect the defect by observing the D-N map, or equivalently the S-matrix. Before going to the general case, we prepare a little more notion.

\subsubsection{Convex polygon}

\begin{figure}[htbp]
 \begin{minipage}{0.49\hsize}
  \begin{center}
 \includegraphics[width=60mm, bb=0 0 254 202]{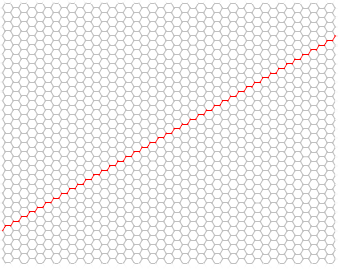}
  \end{center}
  \caption{Straight line 1 in the hexagonal lattice}
  \label{fig:line1}
 \end{minipage}
 \begin{minipage}{0.49\hsize}
  \begin{center}
  \includegraphics[width=60mm, bb=0 0 256 202]{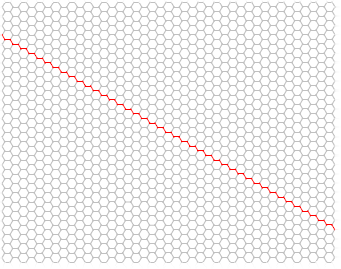}
  \end{center}
  \caption{Straight line 2 in the hexagonal lattice}
  \label{fig:line2}
 \end{minipage}
\end{figure}

\begin{figure}[htbp]
 \begin{minipage}{0.49\hsize}
  \begin{center}
\includegraphics[width=60mm, bb=0 0 256 202]{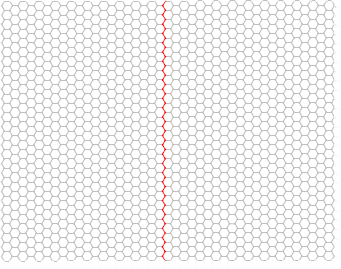}
  \end{center}
  \caption{Straight line 3 in the hexagonal lattice}
  \label{fig:line3}
 \end{minipage}
 \begin{minipage}{0.49\hsize}
  \begin{center}
\includegraphics[width=60mm, bb=0 0 256 202]{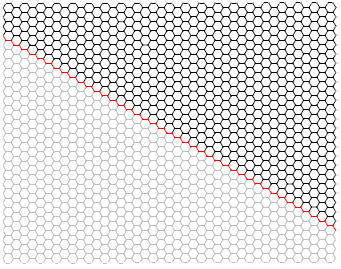}
  \end{center}
  \caption{Half space in the hexagonal lattice}
  \label{fig:half_space}
 \end{minipage}
\end{figure}

There are 3 kinds of straight lines in the hexagonal lattice (Figure \ref{fig:line1}, \ref{fig:line2}, \ref{fig:line3}). 
%
%
%
%
%
%
%
%
Then a half-space is defined as in Figure \ref{fig:half_space}.
%

%
%
\begin{figure}[htbp]
 \begin{minipage}{0.49\hsize}
  \begin{center}
\includegraphics[width=60mm, bb=0 0 251 170]{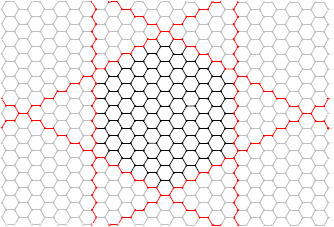}
  \end{center}
  \caption{Hexagonal honeycomb in the hexagonal lattice}
  \label{fig:hexagon}
 \end{minipage}
 \begin{minipage}{0.49\hsize}
  \begin{center}
 \includegraphics[width=42mm, bb=0 0 173 214]{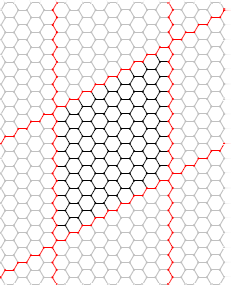}
  \end{center}
  \caption{Hexagonal parallelogram in the hexagonal lattice}
  \label{fig:parallel}
 \end{minipage}
\end{figure}

Let
\begin{equation}
\omega = e^{\pi i/3},
\end{equation}
and put
\begin{equation}
\left\{
\begin{split}
& v_1 = 1 + \omega = (3 + \sqrt3 i)/2,\\
& v_2 = \omega(1 + \omega) = \sqrt 3 i, \\
& v_3 = \omega^2(1 + \omega) = (-3 + \sqrt 3 i)/2.
\end{split}
\right.
\end{equation}
Let $U_h$ be the unit hexagon, which is defined to be a hexagon with vertices $\omega^j, j = 0, 1, \cdots, 5$. For $i,j = 1, 2, 3$, $i \neq j$, and $k \in {\bf Z}$, define the half-space $H_{i,j,k}^{(\pm)}$ by
\begin{equation}
H^{(+)}_{i,j,k} ={\mathop\cup_{l = -\infty}^{\infty}} \, {\mathop\cup_{m\geq k }}\left\{ m v_i + l v_j+ U_h\right\} ,
\end{equation}
\begin{equation}
H^{(-)}_{i,j,k} ={\mathop\cup_{l = -\infty}^{\infty}} \, {\mathop\cup_{m\leq k }}\left\{ m v_i + l v_j+ U_h\right\} .
\end{equation}

By a {\it convex polygon}, we mean an intersection of finite number of half-spaces. In the following, we consider only finite convex polygons. See e.g. Figures \ref{fig:hexagon} and \ref{fig:parallel}.
As typical examples of convex polygon, we consider the following two types of domain : hexagonal honeycomb and hexagonal parallelogram. 
Figure \ref{fig:hexagon} suggests the former, 
and Figure \ref{fig:parallel} the latter.

\begin{figure}[hbtp]
\includegraphics[width=8cm, bb=0 0 350 353]{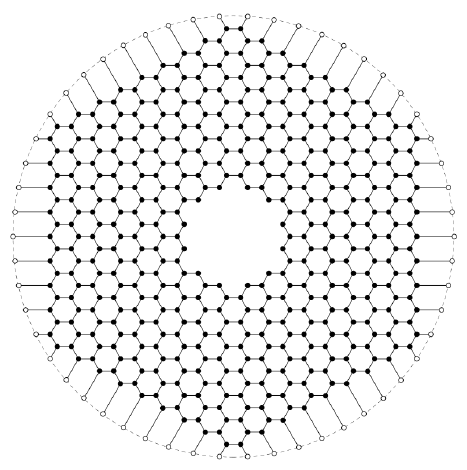}
\caption{Hexagonal honeycomb}
\label{S5HexLattIntDomain}
\end{figure}

 Taking $n$ large enough, we construct the 0th vertical block
$$
B_0 = {\mathop\cup_{k=-n}^n} \left(U_h + k\sqrt{3}i\right).
$$
We next put $U_h' = U_h + \frac{\sqrt3}{2}i$, and make the 1st block
$$
B_{\pm 1} = {\mathop\cup_{k=-n}^{n-1}}\left(U_h' + k\sqrt{3}i \pm \frac{3}{2}\right).
$$
The 2nd block $B_{\pm 2}$ consists of the $2n-1$ number of translated $U_h$'s. 
Repeating this procedure $n$ times, we obtain the {\it hexagonal honeycomb}. Look at  Figure \ref{S5HexLattIntDomain} and imagine the case without hole inside. It is the hexagonal honeycomb.
We attach edges to the vertices with degree 2 on it and also the new vertices on the end points of these edges (white dots in  Figure \ref{S5HexLattIntDomain}), we define the {\it hexagonal honeycomb with boundary}.

We next define another block
$$
B'_0 = {\mathop \cup_{k=0}^m }\left(U_h + k\sqrt 3i\right),
$$
and translate $B'_0$ by $\ell (1 + \omega)$ :
$$
B'_{\ell} = B'_0 + \ell(1 + \omega).
$$
We define
$$
P_{N} = {\mathop \cup_{\ell = - N}^{N}}B'_{\ell},
$$ 
and call it a {\it hexagonal parallelogram} (cf. Figure \ref{S6HexaParallel}). When $N = \infty$, it is called 
{\it graphene nanoribbon} (see \cite{KoKut}).

Letting $e_1, e_2$ be the edges such that $o(e_1) = \omega^2, t(e_1) = \omega$, $o(e_2) = \omega, t(e_2) = 2\omega$, we put
$$
L_{12,k} = e_1\cup e_2 + k \sqrt{3} i, \quad -1 \leq k \leq m,
$$
which are the horizontal edges of $B'_k$. 
We put
$$
Z_k = {\mathop\sum_{\ell=-\infty}^{\infty}}\left(L_{12,k} + \ell(1 + \omega)\right),\quad -1 \leq k \leq m,
$$
and call it a {\it  parallel line in $P_{\infty}$}.


\begin{lemma}\label{Lemmaheaxahoneycritical}
A convex polygon with boundary is a critical graph.
\end{lemma}

Proof. Since the proof is similar in all cases, we give the proof for the hexagonal honeycomb. Letting $\mathcal H$ be a hexagonal honeycomb with boundary, we remove an edge from $\mathcal{H}$. 
By rotating $\mathcal{H}$, we can assume that the removed edge $e_r$ is horizontal. 
Let $B$ be the block of $\mathcal{H}$ containing $e_r$. By translation, we can assume that the bottom of the block $B$ is the edge with vertices $\omega^4$ and $\omega^5$. By translating $B$ to the directions $\pm (1 + \omega)$, we obtain an infinite hexagonal parallelogram $P_{\infty}$, and the associated parallel lines $Z_k$, $ -1 \leq k \leq m$. Let $p_k, q_k$ be the intersection of $Z_k$ with $\partial \mathcal{H}$, where $p_k$ is the left point of intersection, and $q_k$ the right point of intersection. Then $(P,Q)$, where $P = (p_{-1},\cdots,p_m)$ and $Q = (q_{-1},\cdots,q_m)$, is an 
 $(m+2)$-connection of the graph $\mathcal{H}$. Then, if we remove one horizontal edge from $B$, it is no longer a connection, since $B$ has only $m+1$ horizontal edges. \qed
 
\begin{figure}[htbp]
 \begin{minipage}{0.49\hsize}
  \begin{center}
\includegraphics[width=60mm, bb=0 0 269 185]{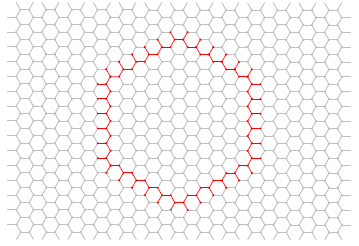}
  \end{center}
  \caption{Outer wall of the hexagonal honeycomb with boundary}
  \label{fig:honeycomb outer wall}
 \end{minipage}
 \begin{minipage}{0.49\hsize}
  \begin{center}
\includegraphics[width=42mm, bb=0 0 175 218]{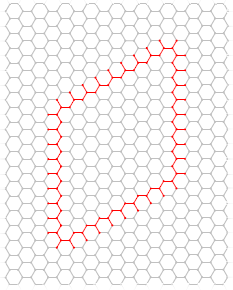}
  \end{center}
  \caption{Outer wall of the hexagonal parallelogram with boundary}
  \label{fig:parallelogram outer wall}
 \end{minipage}
\end{figure}

\medskip
We define the outer wall of a convex polygon taking the hexagonal honeycomb as an example. From the hexagonal honeycomb, we remove all the edges inside and leave only the edges on the periphery. We attach the edge to the vertices with inner angle $2\pi/3$ and a vertex at its end point. Let us call the resulting graph {\it outer wall of the hexagonal honeycomb with boundary} (Figure \ref{fig:honeycomb outer wall}). We can also define, for example, {\it outer wall of the hexagonal parallelogram with boundary} (Figure \ref{fig:parallelogram outer wall}). 
It is another critical case.


\begin{lemma}
\label{LemmaOuterwallhexahoney}
The outer wall of convex polygon with boundary is a critical graph.
\end{lemma}

Proof. As above, we give the proof for the hexagonal honeycomb.
Take two vertices $p_1, q_1$ on the boundary top of the hexagonal honeycomb, $p_1$ being the right to $q_1$. Take $p_2, q_2$ on the bottom, $p_2$ being left to $q_2$. Then, 
$P = \{p_1,p_2\}$, $Q = \{q_1,q_2\}$ are 2-connections. Let $a, b$ be the end points of the edges emanating from $p_1, q_1$, respectively. Then, if we delete the edge $e_{ab}$, $P$ and $Q$ are no longer 2-connected. \qed


\medskip

 In order to detect several defects, we restrict ourselves to the case in which the defects are of the shape of hexagonal honeycomb of one connected component and every component is separated from each other.


\begin{theorem}
\label{Theorem7.5DefectProbe}
Let the defect $\mathcal D$ be of the form $\mathcal D = \cup_{i=1}^N\mathcal D_i$, where $\mathcal D_i \cap \mathcal D_j = \emptyset$ if $i \neq j$ and one of $\mathcal D_i$'s is a convex polygon. 
Assume that the unique continuation property holds on the exterior domain of $ \stackrel{\circ}{\mathcal{D}} $.
Then the set $\{\lambda \in \sigma_e(\widehat H) \, ; \, S(\lambda) = I\}$ is of measure $0$.
\end{theorem}

Proof. Let $\mathcal D_{ow} = \cup_{i=1}^N\mathcal D_{i,ow}$, where $\mathcal D_{i,ow}$ is the outer wall of $\mathcal D_i$. In the assumptions (B-1) $\sim$ (B-4), we take 
$\mathcal V_{int} = \mathcal D_{ow}$,
and $\mathcal V_{ext}$ to be the domain exterior to $\stackrel{\circ}{\mathcal D}$. Then, $\mathcal V = \mathcal V_{ext}\cup\mathcal V_{int}$.  
Note that if $\mathcal D_{ow}$ is replaced by $\mathcal D$, the associated S-matrix is the identity. 
Note that we can apply results in \S 4 to the D-N maps for $\mathcal{D}_{ow} $ and $ \mathcal{D}$.

Suppose there exists a set of positive measure $E \subset \sigma_e (\widehat{H} )$ such that $S(\lambda) = I$ for $\lambda \in E$. 
By taking $\mathcal V_{int}$ to be $\mathcal D_{ow}$ and $\mathcal D$, we see that the D-N map for $\mathcal D_{ow}$, which is the product of each D-N map for $\mathcal D_{i,ow}$,  coincides with that for $\mathcal D$. Suppose $\mathcal D_1$ is a convex polygon. 
Since $\mathcal D_{1,ow}$ and $\mathcal D_1$ are critical,  Theorem \ref{TheoremGraphInv} and Lemmas \ref{Lemmaheaxahoneycritical} and \ref{LemmaOuterwallhexahoney} imply that they coincide as a planar graph. In particular, they must have the same number of arcs, which is not true. 
This proves the theorem.
\qed

\medskip

This theorem asserts that one can detect the  existence of defects from the knowledge of the S-matrix for all energies, however it does not tell us its location. In the next subsection, we find it by employing a different idea.


\subsection{Probing waves}
 \label{subsectionprovingwaves}

Let $\stackrel{\circ}{\mathcal D}$ be a convex polygon. Take a sufficiently large hexagonal parallelogram $\mathcal V_{0,int}$ which contains $D$, and put
\begin{equation}
\mathcal V_{def} = \mathcal V_{0,int}\setminus \stackrel{\circ}{\mathcal D}. 
\end{equation}

In the following, $\Sigma $ denotes $\partial \mathcal{V}_{0,int} = \partial \mathcal{V}_{def}$ and $(\Sigma )_T $, $(\Sigma )_B$, $(\Sigma )_R$ and $(\Sigma )_L$ denote the top, bottom, right and left side of $\Sigma$, respectively. 
We consider the following problem ${\mathcal H}_{def}$ on the region with defects
\begin{equation}
\left\{
\begin{split}
&( -\widehat{\Delta}_{\Gamma} - \lambda) \widehat{u}=0 \quad \text{in} \quad \stackrel{\circ}{\mathcal{V}_{def}} , \\
&\widehat{u} = \widehat{f}   \quad \text{on} \quad \Sigma ,
 \end{split}
\right.
\label{S7_eq_perturbed}
\end{equation}
and the problem $\mathcal H_0$ on the region without defects
\begin{equation}
\left\{
\begin{split}
&( -\widehat{\Delta}_{\Gamma_0}  - \lambda) \widehat{u}_0 =0 \quad \text{in} \quad \stackrel{\circ}{\mathcal{V}_{0,int}} , \\
&\widehat{u} _0 = \widehat{f}   \quad \text{on} \quad \Sigma .
 \end{split}
\right.
\label{S7_eq_nonperturbed}
\end{equation}
We assume :

\medskip
\noindent
{\bf (E)} \ {\it  The number $\lambda$ is not equal to 0, and also neither an eigenvalue of the Dirichlet problem (\ref{S7_eq_perturbed}) nor that of (\ref{S7_eq_nonperturbed}).}

\medskip
Then, the boundary value problems (\ref{S7_eq_perturbed}) and (\ref{S7_eq_nonperturbed}) can be solved uniquely for any data $\widehat{f}$. 
Let $\Lambda(\mathcal H_0)$ and $\Lambda({\mathcal H}_{def})$ be the D-N maps for $\mathcal H_0$ and $\mathcal H_{def}$, respectively.
As is seen in Figure \ref{S6HexaParallel}, there are two types of vertices on $\partial \mathcal{D} $.
One is the  vertex $v \in \partial \mathcal{D}$ with $ \mathrm{deg}_{\mathcal{V}_{def}} (v)=3 $ and inner angle $2\pi /3$,
 and  the other is the vertex  $v' \in \partial \mathcal{D}$ with $ \mathrm{deg}_{\mathcal{V}_{def}} (v' )=2 $ and inner angle $4\pi /3$.

Let  $A_k $, $A'_k $ be the lines  in Figure \ref{S6LineAkprime} in \S 6.2,  and take $ \alpha_{k,0} \in A_k \cap (\Sigma )_T $, $\alpha _{k,\ell'} \in A_k \cap (\Sigma )_R $.
Let $\widehat{u}_{0,k} $ be the solution of (\ref{S7_eq_nonperturbed}) with partial Dirichlet data $\widehat{f}$ such that 
\begin{equation}
\left\{
\begin{split}
&\widehat{f} ( \alpha_{k,0} )=1 ,  \\
&\widehat{f}  (\alpha _{k,\ell'} )= (-1)^{\ell'} , \\
&\widehat{f} (v)=0 \quad \text{for} \quad v\in \Sigma \setminus ( (\Sigma )_R \cup \{ \alpha_{k,0} \} ) ,
\end{split}
\right.
\label{Boundarydataf}
\end{equation}
and partial Neumann data  vanishing on $(\Sigma )_L$.
Lemma \ref{alphareconstruct} implies that $\widehat{u}_{k,0} $ exists uniquely on $\mathcal{V}_{0,int}$ and satisfies 
\begin{equation}
\widehat{u}_{0,k} ( \alpha_{k,\ell} )=
\left\{
\begin{split}
& (-1)^{\ell} \quad {\rm  for} \quad \ell = 0,1,\cdots , \ell', \\ 
&\widehat{u}_{0,k} (x_1 + i x_2 )=0 \quad {\rm  for} \quad x_1 + \sqrt{3} x_2 < a_k. 
\end{split}
\right.
\end{equation}
This solution $ \widehat{u} _{0,k} $ is an analogue of the exponentially growing solution introduced in \cite{Fa66},  \cite{SyUh87}.
We put
\begin{equation}
\widehat{f}_k = \widehat{u}_{0,k} \big| _{\Sigma},
\label{Ddfinefk}
\end{equation}
 and let $\widehat{u}_k $ be the solution of (\ref{S7_eq_perturbed}) with $\widehat{f} = \widehat{f}_k $.


\begin{lemma}
\label{u0k=uklemma}
We take $N$ large enough, and starting from  $k=N$, let $k$ vary downwards. Let $m$ be the largest $k$ such that $A_k $ meets $\mathcal{D}$.
Then, $\widehat{u}_{0,k} = \widehat{u}_k$ on $\mathcal{V}_{def}$ for $k\geq m$.

\label{S7_lem_equi}
\end{lemma}

Proof.
Note that $A_m$ passes through only
 vertices $v\in \partial \mathcal{D}$ with $\mathrm{deg}_{\mathcal{V}_{def}} (v) =3$ and inner angle $2\pi /3$. 
Then, for  any function $\widehat{u}$ on $\mathcal V_{0,int}$,
 we have $((- \widehat{\Delta}_{\Gamma} +1) \widehat{u} )(x_1+ix_2 )=((- \widehat{\Delta}_{\Gamma_0} +1) \widehat{u} )(x_1+ix_2 )$ for any $x_1 +ix_2 \in \stackrel{\circ}{\mathcal{V}_{def}} $ with $x_1 +\sqrt{3} x_2 \geq a_m$. 
By virtue of this equality and $\widehat{u}_{0,k} (x_1 +ix_2 )=0$ for $x_1 + \sqrt{3} x_2 < a_k$, 
  $\widehat{u}_{0,k} $ is a solution of (\ref{S7_eq_perturbed}) with $\widehat{f}= \widehat{f}_k $ if $k \geq m$.
Since (\ref{S7_eq_perturbed}) is uniquely solvable, we have $\widehat u_k = \widehat u_{0,k}$.
\qed

\medskip

We have now arrived at the following probing algorithm for the defects $\mathcal{D}$ of convex polygon in the hexagonal lattice.

\begin{lemma}
\label{LambdaHdeffk=LambdaHofkLemma}
Let $f_k$ be defined by (\ref{Ddfinefk}), and $m$ the number defined in Lemma \ref{S7_lem_equi}. Then, we have 
\begin{equation}
\Lambda ( \mathcal{H}_{def} )\widehat{f}_k = \Lambda (\mathcal{H}_0 )\widehat{f}_k , \quad k \geq m ,
\label{S7_eq_kleqk'}
\end{equation}
and
\begin{equation}
\Lambda ( \mathcal{H}_{def} )\widehat{f}_k \not= \Lambda (\mathcal{H}_0 )\widehat{f}_k , \quad k = m-1 .
\label{S7_eq_kgk'}
\end{equation}

\label{S7_lem_reconstdefect}
\end{lemma}

\begin{figure}[hbtp]
\includegraphics[width=7cm, bb=0 0 316 231]{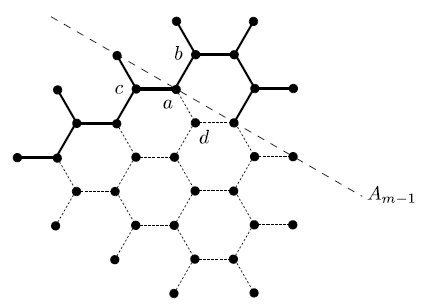}
\caption{Probing for defects}
\label{Probing}
\end{figure}

Proof.
For $k\geq m$, (\ref{S7_eq_kleqk'}) is a direct consequence of Lemma \ref{S7_lem_equi}.
Let us show (\ref{S7_eq_kgk'}).
As  is illustrated in 
Figure \ref{Probing},  $A_{m-1}$ passes through a vertex $a \in \partial \mathcal{D}$ with $ \mathrm{deg} _{\mathcal{V}_{def}} (a)=2 $ and inner angle $4\pi /3$. 
 Let $b,e \in \stackrel{\circ}{\mathcal{V}_{def}} $ be the adjacent vertices of $a$.

Assume that $ \Lambda ( \mathcal{H}_{def} ) \widehat{f}_{m-1} = \Lambda ( \mathcal{H}_0 ) \widehat{f}_{m-1} $. Then, 
by the same reasoning as in the proof of Lemma \ref{u0k=uklemma},  we have $\widehat{u}_{m-1} = \widehat{u}_{0,m-1}$ on $\mathcal{V}_{def}$.

Computing the equation $ (-\widehat{\Delta}_{\Gamma} - \lambda) \widehat{u}_{m-1} =0$ at the above vertex $a \in A_{m-1} \cap \partial \mathcal{D}$, we have
\begin{equation}
- \frac{1}{2} \left( \widehat{u}_{m-1}  (b) + \widehat{u}_{m-1} (c) \right) - \lambda\widehat{u} _{m-1} (a) =0 .
\label{S7_eq_contra1}
\end{equation}
Similarly, the equation $(-\widehat{\Delta}_{\Gamma_0}  +1) \widehat{u}_{0,m-1} =0$ at $a \in A_{m-1} \cap \partial \mathcal{D}$ is 
\begin{equation}
-\frac{1}{3} \left( \widehat{u}_{0,m-1} (b)+ \widehat{u} _{0,m-1} (c) + \widehat{u}_{0,m-1} (d)   \right) -\lambda \widehat{u} _{0,m-1} (a)=0 ,
\label{S7_eq_contra2}
\end{equation}
where $d=a+\omega ^5$. 
Putting $\alpha _{m-1,\ell} =a$ in Figure \ref{S6LineAkprime}, we have $\widehat{u}_{0,m-1} (a)=(-1)^{\ell} $ and $\widehat{u}_{0,m-1} (d)=0$.
However, (\ref{S7_eq_contra1}) and (\ref{S7_eq_contra2}) imply $\widehat{u} _{m-1} (a)= \widehat{u}_{0,m-1} (a)=0$, since  $\widehat{u} _{m-1} = \widehat{u} _{0,m-1} $ on $\mathcal{V}_{def}$.
This is a contradiction. 
\qed

\medskip

Let us pay attention to a relation between Lemma \ref{LambdaHdeffk=LambdaHofkLemma} and 
the partial data problem for the Laplacian on $\mathcal{V}_{def}$. 
In fact, the pair $( (\Sigma )_L , (\Sigma )_R) $ is an $( N+1) $-connection for $\mathcal{H}_0$ and is broken for $\mathcal{H}_{def}$.
Then $\mathcal{H}_0 $ is critical.
It follows that the submatrix of $\Lambda (\mathcal{H}_{def} )$ mapping from $(\Sigma )_R $ to $(\Sigma )_L$ (in the sense of Lemma \ref{S6partialDNdata}) is singular. Therefore, the partial data  problem on $\mathcal{V}_{def} $ in the sense of Lemma \ref{S6partialDNdata} is overdetermined, hence is ill-posed.

\medskip
We consider the  probing for the defects $\mathcal{D} = \bigcup _{j=1}^s \mathcal{D}_j $ where each $\mathcal{D}_j $ is a convex polygon such that $\mathcal{D}_j \cap \mathcal{D}_k = \emptyset $ if $j\not= k$. 
Our method can be applied also to this case.


\begin{lemma}
\label{ProbeLemmaHoneycombDefectslambda}
For $\mathcal{D} = \bigcup _{j=1}^s \mathcal{D}_j $, the assertion of Lemma \ref{S7_lem_reconstdefect} holds.

\label{S7_lem_multihoneycomb}
\end{lemma}

Proof.
For $k\geq m$, the proof is completely the same.
For $k=m-1$,  we have $\widehat{u}_{0,k} = \widehat{u}_k $ in $\mathcal{V}_{def} \setminus \, \stackrel{\circ}{C(\mathcal{D} )} $ where $C(\mathcal{D} )$ is the  hexagonal convex hull of $\mathcal D$.
It then leads to a contradiction at a vertex $v\in \partial C(\mathcal{D} ) \cap \partial \mathcal{D} $ with degree $2$ as in the proof of Lemma \ref{S7_lem_reconstdefect}.
\qed

\medskip
A similar probing procedure is possible by using $B_{\ell}$ in Figure \ref{S6LineBell}. We have only to rotate the domains $\mathcal{H}_0 $ and $\mathcal{H}_{def} $ in the above arguments. Then,  one can enclose the region with defects by the convex hull $C(\mathcal{D} )$.

In Lemma \ref{LemmaS4A2-A1=int2-int1}, we take $\mathcal V_{int,1}$ to be the interior domain without defects, and $\mathcal V_{int,2}$ with defects. 
Using (\ref{LemmaAlambdatoBSigmalambda}), we define $\phi_k \in M_{\lambda}$ by
\begin{equation}
\phi_k = \big(\widehat J^{(-)}(\lambda)^{\ast}\big)^{-1}B^{(+)}_{\Sigma,0}
(\lambda)\widehat f_k.
\end{equation}
Then we have
\begin{equation}
\widehat f_k = \big(B^{(+)}_{\Sigma,0}(\lambda)\big)^{-1}\widehat I^{(-)}(\lambda)^{\ast}\phi_k,
\label{S8widehatfkformula}
\end{equation}
where $B^{(+)}_{\Sigma,0}$ is defined by (\ref{DefineBpmSigmalambda}) for $\mathcal V_{int,1}$. 
Letting $A(\lambda)$ be the scattering amplitude for the lattice with defects, and recalling that the scattering amplitude vanishes for the case without defects,   we have, in view of
Lemmas \ref{LemmaS4A2-A1=int2-int1} and \ref{ProbeLemmaHoneycombDefectslambda},

\begin{equation}
\left\{
\begin{split}
& A(\lambda)\phi_k = 0, \quad {\rm if} \quad k \geq m, \\
& A(\lambda)\phi_k \neq 0, \quad {\rm if} \quad k = m-1.
\end{split}
\right.
\end{equation}
 
We have thus obtained the following theorem.

\begin{theorem}
\label{ProbeTheoremHoneycombDefectslambda}
If the set $\mathcal{D}$ of defects consists of a finite number of convex polygons, its convex hull $C(\mathcal{D} )$ can be computed from the S-matrix $S(\lambda)$ for an arbitrarily fixed energy $\lambda \in \sigma_e(\widehat H)\setminus{\mathcal T}$ satisfying the assumption (E).
\label{ProbeTheoremHoneycombDefects}
\end{theorem}


\begin{thebibliography}{99}
\bibitem{Agmon75}
S. Agmon, \textit{Spectral properties of Schr{\"o}dinger operators and scattering theory}, Ann. Scuola Norm. Sup. Pisa \textbf{2} (1975), 151-218.

\bibitem{AgHo76}
S. Agmon and L. H\"ormander, \textit{Asymptotic properties of solutions of 
differential equations with simple characteristics}, J. d'Anal. Math. 
\textbf{30} (1976), 1-38.



\bibitem{Ando13}
K. Ando, \textit{Inverse scattering theory for discrete Schr{\"o}dinger operators on the hexagonal lattice}, Ann. Henri Poincar{\'e}, \textbf{14} (2013), 347-383.

\bibitem{AndIsoMor}
K. Ando, H. Isozaki and H. Morioka, \textit{Spectral properties of Schr{\"o}dinger operators on perturbed lattices}, Ann. Henri Poincar{\'e} \textbf{17} (2016), 2103-2171.


\bibitem{BuIvKu13}
D. Burago, S. Ivanov and Y. Kurylev, \textit{A graph discretization of the Laplace-Beltrami operator}, J. Spectr. Theory \textbf{4} (2014), 675-714.


\bibitem{Ca80}
A. P. Calder{\'o}n, \textit{On an inverse boundary value problem}, Seminar on Numerical Analysis and its Applications to Continuum Physics, Soc. Brazileira deMathematica, Rio de Janeiro (1980), 65-73.

\bibitem{ChaColPaiRun}
K. Chadan, D. Colton, L. P{\"a}iv{\"a}rinta and W. Rundell, 
\textit{An Introduction to Inverse Scattering and Inverse Boundary Value Porblems}, SIAM, Philadelphia, (1997).

\bibitem{Chung}
R. K. Chung, \textit{Spectral Graph Theory}, Amer. Math. Soc., Providence, Rhode-Island (1997).

\bibitem{Col94} 
Y. Colin de Verdi{\`e}re, \textit{R{\'e}seaux {\'e}lectriques planaires I}, Commentarii Math. Helv., \textbf{69} (1994), 351-374.

\bibitem{Col98}
Y. Colin de Verdi{\`e}re, \textit{Cours Sp{\'e}cialis{\'e}s 4, Spectre de Graphes}, Soc. Math. de France (1998).

\bibitem{ColFra13}
Y. Colin de Verdi{\`e}re and T. Fran{\c c}oise, \textit{Scattering theory for graphs isomorphic to a regular tree at infinity}, J. Math. Phys. \textbf{54} (2013), 063502.

\bibitem{ColGit96} 
Y. Colin de Verdi{\`e}re, I. de Gitler and D. Vertigan, \textit{R{\'e}seaux {\'e}lectriques planaires II}, Commentarii Math. Helv., \textbf{71} (1996), 144-167.

\bibitem{CuSie}
J. C. Cuenin and H. Siedentop, \textit{Dipoles in graphene have infinitely many bound states}, J. Math. Phys. \textbf{55} (1204),  122304.

\bibitem{CuMo90}
E. B. Curtis and J. A. Morrow, \textit{Determining the resistors in a network}, SIAM J. Appl. Math.  \textbf{50} (1990), 918-930.

\bibitem{CuMo91}
E. B. Curtis and J. A. Morrow, \textit{The Dirichlet to Neumann map for a resistor network}, SIAM J. Appl. Math.  \textbf{51} (1991), 1011-1029.

\bibitem{CuMo00}
E. B. Curtis and J. A. Morrow, \textit{Inverse Problems for Electrical Networks}, World Scientific, Singapore-New Jersey-London-Hong Kong (2000).


\bibitem{CuMoMo94}
E. B. Curtis, E. Mooers and J. A. Morrow, \textit{Finding the conductors in circular networks}, Math. Modeling Numer. Anal.  \textbf{28} (1994), 781-813.




\bibitem{CuInMo98}
E. B. Curtis, D. Ingerman and J. A. Morrow, \textit{Circular planar graphs and resisitor networks}, Lin. Alg. and its Appl. \textbf{283} (1998), 115-150.


\bibitem{Do84}
J. Dodziuk, \textit{Difference equations, isoperimetric inequality and transience of certain random walks}, Trans. Amer. Math. Soc., \textbf{284} (1984), 787-794.


\bibitem{Es} M. S. Eskina, \textit{The direct and the inverse scattering problem for a partial difference equation}, Soviet Math. Doklady, \textbf{7} (1966), 193-197. 

\bibitem{Fa56}
L. D. Faddeev, \textit{Uniqueness of the inverse scattering problem}, Vestnik Leningrad Univ. 
\textbf{11} (1956), 126-130.

\bibitem{Fa66}
L. D. Faddeev, \textit{Increasing solutions of the Schr{\"o}dinger eqations}, 
Sov. Phys. Dokl. \textbf{10} (1966), 1033-1035.

\bibitem{Fa76}
L. D. Faddeev, \textit{Inverse problem of quantum scattering theory}, J. Sov. Math. \textbf{5}  (1976), 334-396.

\bibitem{GeNi98}
C. G{\'e}rard and F. Nier, \textit{The Mourre theory for analytically fibred operators}, J. Funct. Anal. \textbf{152} (1989), 202-219.

\bibitem{GoGuVo}
J. Gonz{\'a}lez, F. Guiner and M. A. H. Vozmediano, \textit{The electronic spectrum of fullerens from the Dirac equation}, Nuclear Phys. B \textbf{406} (1993), 771-794.

\bibitem{HiShi04}
Y. Higuchi and T. Shirai, \textit{Some spectral and geometric properties for infinite graphs}, Contemp. Math. \textbf{347} (2004), 29-56.

\bibitem{HiNo09}
Y. Higuchi and Y. Nomura, \textit{Spectral structure of the Laplacian on a covering graph}, Euro. J. of Combinatrics, \textbf{30} (2009), 570-585.

\bibitem{HSSS12}
F. Hiroshima, I. Sakai, T. Shirai and A. Suzuki, \textit{Note on the spectrum of discrete Schr{\"o}dinger operators}, J. Math-for-Industry \textbf{4} (2012), 105-108.


\bibitem{HoVol3}
L. H{\"o}rmander, \textit{The Analysis of Linear Partial Differential Operators III, Pseudodifferential operators}, Springer Verlag, Berlin-Heidelberg-New York-Tokyo (1985).

\bibitem{Ikehata99}
M. Ikehata, \textit{Reconstruction of obstacles from boundary measurements}, Wave Motion \textbf{30} (1999), 205-223.

\bibitem{Ikehata04}
M. Ikehata, \textit{Inverse scattering problems and the enclosure method}, Inverse Problems \textbf{20} (2004), 533-551.

\bibitem{Is03}
H. Isozaki, \textit{Inverse spectral theory}, in {\it Topics In The Theory of Schr{\"o}dinger Operators}, eds. H. Araki, H. Ezawa, World Scientific (2003), pp. 93-143.


\bibitem{IsKo} H. Isozaki and E. Korotyaev, \textit{Inverse problems, trace formulae for discrete Schr\"{o}dinger operators},   Ann.  Henri Poincar{\'e}, \textbf{13} (2012), 751-788.


\bibitem{IsMo1}
H. Isozaki and H. Morioka, \textit{A Rellich type theorem for discrete Schr{\"o}dinger operators}, Inverse Problems and Imaging, \textbf{8} (2014), 475-489.

\bibitem{IsMo2}
H. Isozaki and H. Morioka, \textit{Inverse scattering at a fixed energy for discrete Schr{\"o}dinger operators on the square lattice}, Ann. Inst. Fourier \textbf{65} (2015), 1153-1200.

\bibitem{Kuroda}
S. T. Kuroda, \textit{Scattering theory for differential operators, I, II}, J. Math. Soc. Japan \textbf{25} (1973), 75-104, 222-234.

\bibitem{KhNo87}
G. M. Khenkin and R. G. Novikov, \textit{The $\overline\partial$-equation in the multi-dimensional inverse scattering problem}, Russian Math. Surveys \textbf{42} (1087), 109-180.

\bibitem{KoOnoSuna}
T. Kobayashi, K. Ono and T. Sunada, \textit{Periodic Schr{\"o}dinger operators  on a manifold}, Forum Math. \textbf{1} (1989), 69-79.


\bibitem{KCSSSHSOTTN12}
T. Kondo, S. Casolo, T. Suzuki, T. Shikano, M. Sakurai, Y. Harada, M. Saito, M. Oshima, M. Trioni, G. Tantardini and J. Nakamura, 
\textit{Atomic-scale  characterization of nitrogen-dopted graphite : Effects of dopant nitrogen on the local electronic structure of the surrounding carbon stoms}, Phyisical Review B \textbf{86}, 035436 (2012).

\bibitem{KoKut}
E. Korotyaev and A. Kutsenko, \textit{Zigzag nanoribbons in external electric fields}, Asymptotic Anal. \textbf{66} (2010), 187-206.

\bibitem{KoSa13}
E. Korotyaev and N. Saburova, \textit{Schr{\"o}dinger operators on periodic graphs}, J. Math. Anal. Appl. \textbf{420} (2014), 576-611.


\bibitem{KoShiSu98} M. Kotani, T. Shirai, and T. Sunada, \textit{Asymptotic behavior of the transition probability of a random walk on an infinite graph}, J. Funct. Anal. \textbf{159} (1998), 664-689.

\bibitem{KuPo}
P. Kuchment and O. Post, \textit{On the spectra of carbon nano-structures}, 
Commun. Math. Phys. \textbf{275} (2007), 805-826.

\bibitem{MoWo}
B. Mohar and W. Woess, \textit{A survey on spactra of infinite graphs}, 
Bull. London Math. Soc. \textbf{21} (1989), 209-234.

\bibitem{Muho77}
	R. G. Muhometov, \textit{The problem of recovery of a two-dimensional Riemannian metric and integral geometry}, Soviet Math. Dokl. \textbf{18} (1977), 27-31.


\bibitem{Na88}
A. Nachman, \textit{Reconstruction from boundary measurements}, Ann. of Math. \textbf{128} (1988), 531-576.


\bibitem{Na96}
A. Nachman, \textit{Global uniqueness for a two-dimensional inverse boundary value problem}, Ann. of Math. \textbf{143} (1996), 71-96.

\bibitem{Na14}
S. Nakamura, \textit{Modified wave operators for discrete Schr{\"o}dinger operators with long-range perturbations}, J. Math. Phys. \textbf{55} (2014), 112101.

\bibitem{NeGuPeNov}
A. H. C. Neto, F. Guiner, N. M. R. Peres, K. S. Novoselov and A. K. Geim, 
\textit{The electronic properties of graphene}, Rev. Mod. Phys. \textbf{81} (2009), 109-162.


\bibitem{No88}
R. G. Novikov, \textit{A multidimensional inverse spectral problem for the equation 
$- \Delta\psi + (v(x)-E)\psi=0$}, Funct. Anal. Appl. \textbf{22} (1988), 263-272.

\bibitem{Ober00} 
R. Oberlin, \textit{Discrete inverse problems for Schr\"{o}dinger and resistor networks}, Research archive of Research Experiences for Undergraduates program at Univ. of Washington, (2000).

\bibitem{PesUhl}
L. Pestov and G. Uhlmann, \textit{Two dimensional compact simple Riemannian manifolds are boundary distance rigid}, Ann. of Math. (2) \textbf{161} (2005), 1093-1110.

\bibitem{Seme}
G. W. Semenoff, \textit{Condensed-matter simulation of a three-dimensional anomaly}, Phys. Rev. Lett. \textbf{53} (1984), 244-2452.

\bibitem{Sha} W. Shaban and B. Vainberg, \textit{Radiation conditions for the difference Schr\"{o}dinger operators}, Applicable Analysis, \textbf{80} (2001), 525-556.



\bibitem{Ship} S. P. Shipman, \textit{Eigenfunctions of unbounded support for 
embedded eigenvalues of locally perturbed periodic graph operators}, Commun. Math. Phys. \textbf{332} (2014), 605-626.

\bibitem{S99} T. Shirai,  \textit{The spectrum of infinite regular line graphs}, Trans. Amer. Math. Soc., \textbf{352} (1999), 115-132.

\bibitem{Suna90}
T. Sunada, \textit{A periodic Schr{\"o}dinger operator on abelian cover}, 
J. Fac. Sci. Univ. Tokyo, Sect. IA, Math. \textbf{37} (1990), 575-583.

\bibitem{Su13}
A. Suzuki, \textit{Spectrum of the Laplacian on a covering graph with pendant 
edges : The one-dimensional lattice and beyond}, Lin. Alg. and its Appl. \textbf{439} (2013), 3464-3489.


\bibitem{SyUh87}
J. Sylvester and G. Uhlmann, \textit{A global uniqueness theorem for an inverse boundary value problem}, Ann. of Math. \textbf{125} (1987), 153-169.

\bibitem{Uhl09}
G. Uhlmann, \textit{Electrical impedance tomography and Calder{\'o}n's problem}, Inverse Problems \textbf{25} (2009), 123011.

\bibitem{Yaf}
D. Yafaev, \textit{Mathematical Scattering Theory}, Amer. Math. Soc., Providence, Rhode Island,  (2009).









\end{thebibliography}
\end{document}